\documentclass[aos]{imsart}
\usepackage[utf8]{inputenc}
\usepackage[numbers]{natbib}
\usepackage{amsmath}
\usepackage{amsfonts}
\usepackage{amssymb}
\usepackage{amsthm}
\usepackage{multirow}
\usepackage{enumerate}
\usepackage{bbm}
\usepackage{bm}
\usepackage{graphicx}
 \usepackage{mathtools}
\usepackage{xcolor}
\usepackage{tikz}
\usetikzlibrary{datavisualization}
\usetikzlibrary{datavisualization.formats.functions}
\usepackage[colorlinks,citecolor=blue,urlcolor=blue]{hyperref}

\arxiv{2506.22588}
\startlocaldefs
\def\journalname#1{\def\journal@name{#1}}

\newcommand*{\iidsim}{\stackrel{\mathrm{i.i.d.}}{\sim}}

\newcommand*{\indicator}[1]{\mathbf{1}\left\{#1\right\}}

\newcommand*{\bracks}[1]{\left\{#1\right\}}
\newcommand*{\sqbrack}[1]{\left[#1\right]}
\newcommand*{\paren}[1]{\left(#1\right)}

\newcommand*{\rme}{\mathrm{e}}
\newcommand*{\rmd}{\mathrm{d}}
\newcommand*{\KL}{\mathrm{KL}}

\DeclareMathOperator{\var}{Var}

\DeclareBoldMathCommand{\e}{e}
\DeclareBoldMathCommand{\u}{u}
\DeclareBoldMathCommand{\v}{v}
\DeclareBoldMathCommand{\w}{w}
\DeclareBoldMathCommand{\q}{q}
\DeclareBoldMathCommand{\p}{p}
\DeclareBoldMathCommand{\R}{R}
\DeclareBoldMathCommand{\valpha}{\alpha}
\DeclareBoldMathCommand{\vbeta}{\beta}
\DeclareBoldMathCommand{\vmu}{\mu}
\DeclareBoldMathCommand{\vDelta}{\Delta}
\DeclareBoldMathCommand{\vdelta}{\delta}
\DeclareBoldMathCommand{\vgamma}{\gamma}
\DeclareBoldMathCommand{\H}{H}
\DeclareBoldMathCommand{\X}{X}
\DeclareBoldMathCommand{\Y}{Y}
\DeclareBoldMathCommand{\Z}{Z}
\DeclareBoldMathCommand{\x}{x}
\DeclareBoldMathCommand{\y}{y}
\DeclareBoldMathCommand{\ones}{1}
\DeclareBoldMathCommand{\veta}{\eta}
\DeclareBoldMathCommand{\vpi}{\pi}

\theoremstyle{definition}
\newtheorem{theorem}{Theorem}[section]
\newtheorem{corollary}[theorem]{Corollary}
\newtheorem{proposition}[theorem]{Proposition}

\newtheorem{remark}[theorem]{Remark}

\newtheorem{definition}[theorem]{Definition}

\newtheorem{lemma}[theorem]{Lemma}
\newtheorem{claim}[theorem]{Claim}

\newcommand{\eq}{=}
\definecolor{red}{RGB}{213,94,0}
\definecolor{green}{RGB}{0,158,115}

\endlocaldefs

\journalurl{}
\journalname{}

\begin{document}

\begin{frontmatter}
  \title{Anytime-Valid Tests for Sparse Anomalies}
  \runtitle{Anytime-Valid Tests for Sparse Anomalies}

  \begin{aug}
    \author[A]{\fnms{Muriel F.}~\snm{Pérez-Ortiz}\ead[label=e1]{m.f.perez.ortiz@tue.nl}},
    \author[A]{\fnms{Rui
        M.}~\snm{Castro}\ead[label=e2]{r.m.pires.da.silva.castro@tue.nl}}
    \and
    \author[A]{\fnms{Ivo V.}~\snm{Stoepker}\ead[label=e3]{i.v.stoepker@tue.nl}}
    \address[A]{Department of Mathematics and Computer Science,
      Eindhoven University of Technology\printead[presep={,\ }]{e1,e2,e3}}
    \runauthor{M.F. Pérez-Ortiz et al.}
  \end{aug}
\begin{abstract}
  We consider the problem of testing sequentially for the presence of sparse
  anomalies among a large number of data streams. To this end, we design and
  analyze Anytime-Valid (AV) tests, which retain type-I error control at
  arbitrary stopping times. Existing results address exclusively the
  nonsequential case, which exhibits a subtle phase transition between two
  regimes where tests are either powerless or powerful. In our sequential
  setting, we argue, two challenges arise: (1) the standard analysis of AV tests
  cannot be executed in the relevant sample-size regime; and (2) standard
  constructions of parameter-adaptive AV tests are either analytically
  intractable or computationally unfeasible. This work addresses these
  challenges. Borrowing insights from the nonsequential literature, we propose a
  framework to analyze AV tests and their shortest possible sample sizes. Under
  this framework, we show that, in the Gaussian location setting, the oracle AV
  test has a delicate threshold behavior that is related to---but not implied
  by---the phase transition observed in optimal nonsequential tests. Our main
  results include a computationally efficient, parameter-adaptive AV test; we
  show that it achieves the same threshold behavior as the oracle AV test.
  Numerical simulations illustrate these theoretical findings.
\end{abstract}


\end{frontmatter}

\section{Introduction}
\label{sec:introduction}

In essence, the statistical challenge behind many scientific questions is the
problem of anomaly detection. For example, in detecting disease outbreaks, the
statistical problem is to detect, among a large number of places, the existence
of places with anomalously high number of disease
cases~\citep{kulldorff2005stp}; in monitoring production processes, it is to
detect, among a large number of the product's characteristics, anomalous
deviations from existing quality standards~\citep{Bersimis:2007,zagar2022}. In
these applications, it is desirable to make decisions sequentially in time while
monitoring data continuously as it is collected. Indeed, in disease
surveillance, early detection of disease outbreaks is critical to curb further
spread; in production process monitoring, early detection of process instability
prevents defects from escalating. Inspired by these challenges, we investigate
the detection of anomalies that manifest themselves as sparse effects within a
large set of data streams.

We consider $K\in\mathbf{N}$ data streams, which are monitored continuously over
time. For each stream $i \in [K] := \{1, \dots, K\}$, a real-valued sequence
$(X_{i, t})_{t \geq 1}$ is observed. We assume that each stream is either
anomalous or not, and that anomalies affect streams in their entirety, starting
from $t=1$. Our primary interest is in detecting sparse anomalies; that is, when
anomalies might only affect a small fraction of the streams. Our goal is to
detect the presence of anomalous streams sequentially (as the data comes in and
the streams grow in length) as soon as possible.

The data streams and the potential anomalies are modeled through a probabilistic
contamination model that is commonly used both for sparse-signal detection and
multiple testing~\citep{donoho_higher_2004, ingster_adaptive_2002}. Despite the
existence of fixed-sample methods for sparse anomaly detection that are valid
and powerful under very general distributional assumptions (within such
contamination
models)~\citep{arias-castro_detection_2019,donoho_higher_2015,stoepker_sparse_2023,stoepker_anomaly_2022,cai_optimal_2014,cai_estimation_2007,cai_optimal_2011},
translating these insights to design sequential tests turns out to be, as we
will see, surprisingly challenging. These challenges are already apparent in
what is arguably the simplest version of the problem that shows a rich
statistical behavior: a Gaussian location model. Under this model,
\begin{equation}\label{eq:gaussian-alternative}
  \begin{aligned}
    &A_{1}, \dots, A_{K} \iidsim \mathrm{Bernoulli}(\varepsilon) \text{ and, for
      }  t\geq 1,\\
    &X_{i, 1}, \dots, X_{i, t} \ | \ A_i
      \iidsim
      \mathrm{Normal}(\delta A_i, 1)
      \text{ independently for each } i \in [K].
  \end{aligned}
\end{equation}
When $A_i=1$, we call stream $i$ anomalous; when $A_i=0$, we refer to it as
normal. Note that the anomaly indicators $(A_i)_{i\in[K]}$ are latent
variables---they are not observed---and the only observables up to time $t$ are
the data streams $(X_{i, s})_{i\in[K],s\leq t}$. The parameter
$\varepsilon \in [0, 1]$ controls the ``sparsity'' of the anomalies as it is the
average proportion of anomalous streams. Our  main interest is in the case that
$\varepsilon$ is  close to zero. In the model, the distribution of the
observations in anomalous and normal streams differ by a mean shift of size
$\delta$. Furthermore, conditionally on the indicators $(A_i)_{i\in[K]}$, the
streams are probabilistically independent. The model is stationary in the sense
that, for $t \geq 1$, the distribution of each data stream remains unchanged.

We study anomaly detection as a hypothesis testing problem: the null hypothesis
represents the absence of anomalies; the alternative hypothesis, their presence.
Within the probabilistic model~(\ref{eq:gaussian-alternative}), the no-anomaly
null hypothesis $\mathcal{H}_0$ and its alternative $\mathcal{H}_1$ are,
respectively,
\begin{equation}\label{eq:generic-alternative}
  \mathcal{H}_0: \varepsilon = 0,
  \qquad\text{ and }\qquad
  \mathcal{H}_1: \varepsilon > 0 \text{ and } \delta > 0.
\end{equation}
Intuitively, there is a trade-off between the sparsity $\varepsilon$ and the
signal strength $\delta$ that controls the statistical difficulty of the testing
problem~\eqref{eq:generic-alternative}. Indeed, for smaller values of
$\varepsilon$ (fewer anomalies), larger values of $\delta$ (stronger signal) will
be required to reject $\mathcal{H}_0$ confidently---results in this direction
are discussed in Section~\ref{sec:relat-work-prel}. In order to simplify the
exposition, we focus on the one-sided alternative that $\delta>0$, but all the
methods presented apply to the two-sided, $\delta\neq 0$ alternative as well.

Our goal is to test the hypotheses~\eqref{eq:generic-alternative} sequentially.
In particular, we design and evaluate Anytime-Valid (AV) tests for the
hypothesis testing problem~(\ref{eq:generic-alternative}). AV tests are designed
to solve the problem of multiple testing over time. Indeed, the defining
characteristic of AV tests is that if experimentation is stopped using any
stopping rule, they retain type-I error guarantees (we expand on this in
Section~\ref{sec:relat-work-prel}).

A natural starting point is to construct AV tests that build on known
fixed-sample tests; for example, by applying them at each timepoint $t$ and
using a Bonferroni-type correction to account for multiple testing over time.
Although this Bonferroni-type approach could be applied under more general
distributional assumptions than~(\ref{eq:gaussian-alternative}), it performs
poorly~(see Appendix~\ref{app:comparison-bonferroni}). This fact motivates
investigating AV test constructions tailored specifically to our sparse
detection setting, rather than extending known fixed-sample approaches. Not only
is this necessary to obtain effective AV tests, but, since AV tests can be used
to formulate fixed-sample tests, this investigation also yields powerful
fixed-sample tests as an interesting by-product.

Designing AV tests for~(\ref{eq:generic-alternative}) connects and makes
contributions to two branches of literature.
On the one hand, the model~(\ref{eq:gaussian-alternative}) and the sparse
alternative in \eqref{eq:generic-alternative} have been studied in the
literature on high-dimensional, sparse signal detection (e.g.,
\citep{donoho_higher_2004, ingster_adaptive_2002}). Of particular interest in
this line of work is the design of adaptive tests; that is, tests that are
optimal across all sparsity regimes $\varepsilon > 0$ without requiring any
knowledge about the value of $\varepsilon$. This work also designs adaptive
tests for this problem, but the main difference with existing work is that the
proposed AV tests have strong guarantees that hold without a prespecified time
(sample size) horizon, whereas previous work is on fixed-sample tests.
On the other hand, the sequential aspect of this work connects with the AV
testing literature (e.g., \cite{ramdas_game-theoretic_2023,johari_always_2022}),
born of renewed interest in power-one tests (e.g.,
\cite{lai_power-one_1977,robbins_expected_1974}). The main difference with this
line of literature is that the sparse alternatives considered
in~(\ref{eq:generic-alternative}) call for a different type of analysis than the
existing ones. Indeed, the analysis of existing AV tests revolves exclusively
along the time dimension, while, as we will see, the testing problem
\eqref{eq:generic-alternative} is amenable to an asymptotic analysis along the
$K$ dimension (a large number of streams) for each fixed time $t$. We further
discuss related work, along with necessary preliminaries in
Section~\ref{sec:relat-work-prel}.

\paragraph*{Informal summary of contributions} Informally, our contributions are
threefold. (A detailed summary of the contributions of this article is found in
Section~\ref{sec:contribution}.)
\begin{enumerate}
\item A framework for analysis: the nuances of the high-dimensional testing
  problem~(\ref{eq:generic-alternative}) render the analysis of optimality
  measures delicate. We provide, for the problem of testing sparse anomalies, an
  asymptotic, large-$K$ framework to analyze the optimality of AV tests. This
  framework makes it possible to quantify, under ${\cal H}_1$
  in~(\ref{eq:generic-alternative}), both the potential growth of evidence
  against ${\cal H}_0$, and to quantify the sample size at which a decision can
  be made at a prespecified level.
\item Analysis of oracle AV tests: under the proposed framework, we analyze an
  oracle test that has perfect knowledge of the alternative
  in~(\ref{eq:generic-alternative}). Loosely speaking, we find that there is a
  threshold phenomenon: there exists a specific time until which both evidence
  accumulation and decision (at a given level) are impossible; after this
  specific time, achieving both objectives becomes unavoidable.
\item An optimal adaptive AV test: we develop an adaptive AV test for
  \eqref{eq:generic-alternative}---it does not require knowledge of the sparsity
  $\varepsilon$ or the signal strength $\delta$. We show that the proposed test
  is optimal under our proposed framework. The construction of an adaptive AV
  tests is a challenge on its own: a naive adaptation of known fixed-sample
  tests provably fails (see Appendix~\ref{app:e-max-test}), and standard
  adaptive constructions from the AV literature either exhibit unsatisfactory
  performance or are computationally unfeasible (see Section~\ref{sec:prior}).
\end{enumerate}
Incidentally, since any AV test can be used to produce a fixed-sample test, a
by-product of this work is an adaptive, fixed-sample test
for~(\ref{eq:generic-alternative}). This test exhibits the optimal power
behavior shown by existing powerful adaptive tests from the fixed-sample
literature~\citep{donoho_higher_2004,porter_beyond_2020,li_higher_2015}.
Although the main focus of this article is in AV testing, this is not a minor
fact, and is discussed, among other interesting emerging questions, in
Section~\ref{sec:discussion}.

\paragraph*{Outline of the rest of the introduction} The remainder of this
section is structured as follows: Section~\ref{sec:relat-work-prel} discusses
related work and preliminaries in both sparse-anomaly detection and AV tests,
building up to the asymptotic analysis framework in Section~\ref{sec:framework};
Section~\ref{sec:contribution} contains a formal summary of the contributions;
Section~\ref{sec:notation} presents the notation that is used in the remainder
the article; and Section~\ref{sec:outline} outlines the rest of the article.

\subsection{Preliminaries and related work}
\label{sec:relat-work-prel}
This section discusses related work and preliminary results in both the
sparse-signal (Section~\ref{sec:state-things}) and the AV literature
(Section~\ref{sec:explain-av}). We highlight nuances encountered in both lines
of work: together these motivate our framework to judge optimality of AV tests
for the problem of testing sparse signals---this framework is introduced in
Section~\ref{sec:framework}.

\subsubsection{Detection and identification of sparse
  anomalies}\label{sec:state-things}

In the context of the model~\eqref{eq:gaussian-alternative} and related
settings, research has focused on two goals: studying the power of fixed-sample
tests under ``sparse'' alternatives \eqref{eq:generic-alternative} and
identifying which streams are anomalous. We call the first objective detection;
the second, identification. Given the Gaussian structure
of~(\ref{eq:gaussian-alternative}), sufficiency considerations simplify the
exposition of existing (fixed-$t$) results. Indeed, the vector of standardized
stream sums $(Z_{1, t}, \dots, Z_{K, t})$ with
$Z_{i, t} = (1 / \sqrt{t})\sum_{s\leq t} X_{i,s}$ is a sufficient statistic
for $(\varepsilon, \delta)$, and
$\mathbf{E}_{\varepsilon, \delta}[Z_{i,t}|A_i]=\delta \sqrt{t} A_i$ under a
specific alternative given by $(\varepsilon,\delta)$. Hence, to the effect of
fixed-$t$ procedures, it is enough to present the case that $t=1$ because the
case where $t\geq 1$ can be recovered by rescaling $\delta$. We return to
$t\geq 1$ at the end of the section.

For detection, the hypothesis problem~(\ref{eq:generic-alternative}), a
nonasymptotic characterization of the $\varepsilon$-$\delta$ trade-off is
challenging despite the simple appearance of~(\ref{eq:gaussian-alternative}). It
has been useful to take an asymptotic perspective, where the fixed-sample power
properties of tests are studied as the model dimension $K$ grows. In sparse
regimes, where $\varepsilon = \varepsilon(K) \to 0$ as $K\to\infty$, a first
question is to identify the sequences $\delta = \delta(K)$ against which
powerful tests exist. Such sequences have been identified in the Gaussian
version of the problem (and in more general
settings)~\citep{ingster_minimax_1994}, and parameter-adaptive versions (when
the parameters $\delta$ and $\varepsilon$ in the alternative hypothesis are
unknown) have been proposed
\citep{ingster_adaptive_2002,donoho_higher_2004,porter_beyond_2020}. Numerous
extensions of this setting have been considered in the literature as
well~\citep{arias-castro_detection_2019,donoho_higher_2015,stoepker_sparse_2023,stoepker_anomaly_2022}.
All existing results focus on fixed-sample (fixed-$t$) hypothesis testing. This
article departs from this line of work in that it studies a sequential AV
version of the sparse anomaly detection problem, an extension motivated by
sequential analysis \citep{tartakovsky_sequential_2014} and its applications;
for example, in syndromic surveillance \citep{Zamba2021}.

Consider the problem of detection, the hypothesis testing problem
\eqref{eq:generic-alternative}, in the case where the true parameters
$\varepsilon^\star>0$ and $\delta^\star > 0$ in the alternative hypothesis are
specified. In this case, the alternative hypothesis becomes simple, and the
likelihood ratio test is the most powerful fixed-sample test by
Neyman--Pearson's Theorem. Despite the most powerful test being known, the
sparsity--signal-strength trade-off remains nonobvious for fixed $K$. This is
where the asymptotic, large-$K$ point of view becomes
valuable~\citep{ingster_minimax_1994,donoho_higher_2004}. The insightful
parametrization $\varepsilon^\star = \varepsilon^\star(K)$ in terms of the
number of streams $K$ turns out to be $\varepsilon^\star = K^{-\beta^\star}$ for
$\beta^\star\in [0, 1]$. When $\beta^\star<1/2$, detection is easy even when
$\delta^\star$ is close to zero. To see this, consider the standardized grand
mean of the observations across streams
$T_K = (1 / \sqrt{K})\sum_{i=1}^K X_{i,1}$, a Gaussian random variable with
unit variance under~(\ref{eq:gaussian-alternative}). The mean of $T_K$ is zero
under the null hypothesis and, under the alternative hypothesis, $T_K$ has mean
$\delta^\star \varepsilon^\star \sqrt{K} = \delta^\star K^{1/2 - \beta^\star}$.
Therefore, when $\beta^\star<1/2$ the two distributions can be distinguished
with a test based on $T_K$ even if $\delta^\star$ converges to zero at a
slow-enough rate. The more interesting regime is when $\beta^\star>1/2$. In that
case, the insightful parametrization $\delta^\star = \delta^\star(K)$ turns out to
be $\delta^\star=\sqrt{2r^\star\ln K}$ for $r^\star > 0$. Here,
$\delta^\star\to\infty$ as $K\to\infty$ and, surprisingly, it is no longer
possible to detect sequences $\delta = \delta(K)$ such that $\delta\to 0$ in the
same limit. The problem exhibits a phase transition that has been characterized
sharply~\citep{porter_beyond_2020} and is best described in terms of the pair of
parameters $(\beta^\star, r^\star)$. To this end, define
\begin{equation}
  \label{eq:detection_boundary}
  \rho(\beta)
  :=
  \begin{cases}
    \beta - \frac{1}{2} & \text{ if } \beta \in (1/2, 3/4);\\
    (1 - \sqrt{1 - \beta})^2 & \text{ if } \beta\in [3/4, 1].
  \end{cases}
\end{equation}
The function $\rho$ plays the role of a ``detection boundary''; it divides the
parameter space in two regions. On the one hand, when
$r^\star<\rho(\beta^\star)$, all tests (with significance level strictly smaller
than one) become powerless as $K\to\infty$. On the other hand, when
$r^\star>\rho(\beta^\star)$, the likelihood ratio test (with significance level
$\alpha\in(0,1)$) has power converging to one as
$K\to\infty$~\citep{ingster_adaptive_2002,ingster_minimax_1994}.

In the same $(r^\star, \beta^\star)$ parametrization, the problem of identifying
which streams are anomalous has also been considered. In order to ensure that
the number of errors is small (how many normal streams are labeled as
anomalous), then necessarily $r^\star>1$. This is due to the fact that, under
the null model, $\max_{i\in[K]} X_{i,1}$ concentrates around $\sqrt{2\ln K}$,
and the stream-normalized means must be above that level to be identified as
anomalous with confidence. More lenient error metrics, such as False Discovery
Rate (FDR) and Non-Discovery Rate (NDR), can be considered. Nevertheless, to
ensure that both FDR and NDR converge to zero as $K\to\infty$, it is necessary
and sufficient that $r^\star>\beta^\star$ for $\beta^\star\in(0,1]$
\citep{donoho_higher_2004,DS_Haupt_2011,Arias-Castro:2017}. In contrast with
other statistical problems where testing and estimation are two sides of the
same coin, these results show a very different picture in this problem: there is
a large region in the parameter space where the presence of anomalies can be
detected, but their location cannot be identified. The identification--detection
gap is particularly large when $\beta^\star$ is close to $1/2$.
Figure~\ref{fig:regions} shows both the ``detection boundary'' given
by~(\ref{eq:detection_boundary}) and region of the $\beta$-$r$~plane where there
is a separation between detection and identification (in the sense of asymptotic
FDR and NDR control). As we will see, the gap between testing and identification
has implications in the anytime-valid setting as well.

\begin{figure}[h]
  \centering
    \begin{tikzpicture}[baseline]
      \datavisualization [ scientific axes,
      y axis={label = {$r$}},
      x axis={label = {$\beta$}},
      legend = right,
      visualize as smooth line/.list={ro,ex},
      style sheet=strong colors,
      style sheet=vary dashing,
      ex={label in legend={text=$r\eq\beta$}},
      ro={label in legend={text=$r\eq\rho(\beta)$}},
      data/format=function ]
      data [set=ro] {
        var x : interval [0.5:0.75];
        func y = (\value x - 1/2); }
      data [set=ro, samples = 500] {
        var x : interval [0.75:1];
        func y = (1 - sqrt(1 - \value x))^2;
      }
      data [set=ex] {
        var x : interval [0.5:1];
        func y = (\value x);
      }
      info {
        \draw [green] (visualization cs: x=0.67, y=0.82)
        node [above] {Identification \& Detection};
      }
      info {
        \draw [green] (visualization cs: x=0.75, y=0.4)
        node [above] {Detection};
      }
      info {
        \draw [green] (visualization cs: x=0.9, y=0.1)
        node [above] {Neither};
      };
    \end{tikzpicture}
    \caption{\label{fig:regions} What is possible? Regions in the $\beta$-$r$~plane where
      identification or detection of anomalies is possible in an asymptotic
      sense. In the middle region, only detection is possible and $\rho$ is
      given by~(\ref{eq:detection_boundary})~(see
      Section~\ref{sec:state-things}). }
\end{figure}

For the composite alternative hypothesis~\eqref{eq:generic-alternative}, there
exist parameter-adaptive, fixed-sample tests that achieve the same ``detection
boundary'' as an oracle likelihood ratio test, a test that uses the (unknown)
true parameters. An example of such a parameter-adaptive test is the test based
on the Higher Criticism (HC) statistic~\citep{donoho_higher_2004}, the
null-standardized empirical cumulative distribution of the standardized stream
averages (see Section~\ref{sec:simulations}). Even though other tests (the
Generalized Likelihood Ratio Test (GLRT) and its associated score tests) also
achieve this goal, there is no known test that dominates the others in an
uniform sense~\citep{porter_beyond_2020,li_higher_2015}.

The goal of this work is to construct AV adaptive tests that are optimal in a
well-defined sense, as described in the next section.

\subsubsection{Anytime-valid tests and optimality }
\label{sec:explain-av}

AV tests are a specific type of sequential tests, decision procedures that make
use of the information available at each time point. This information is modeled
by the data filtration
$\mathbf{F} := (\sigma(X_{i, s} : s\leq t, i\in [K]))_{t\geq 1}$. A sequential
test is a pair $(\varphi, \tau)$, where $\tau$ is an $\mathbf{F}$-stopping time,
and $\varphi = (\varphi_t)_{t\in \mathbf{N}}$ is an $\mathbf{F}$-adapted
sequence of decision functions; that is, $\varphi_t\in \{0, 1\}$ for each
$t\in\mathbf{N}$. The stopping time $\tau$ signifies the (random) sample size of
the test and $\varphi_\tau$ signifies the decision made at that sample size. In
this case, $\varphi_\tau = 1$ means rejecting the null hypothesis;
$\varphi_\tau = 0$, accepting it. For the applications that we have in mind,
there is no reason to stop in the absence of anomalies. Hence, our primary
interest is in sequential tests $(\varphi, \tau)$ for which $\varphi_\tau = 1$;
in other words, sequential tests that only stop to reject the null hypothesis.

The anytime validity of a sequential test is actually a property of the decision
sequence alone. Indeed, we say that the decision sequence $\varphi$ is anytime valid at
level $\alpha$ if, for any $\mathbf{F}$-stopping time $\tau$, the type-I error
$\alpha_0(\varphi, \tau) := \mathbf{P}_0\{\varphi_\tau = 1\}$ of the test
$(\varphi, \tau)$ is smaller or equal than $\alpha$. Typically, for a fixed
decision sequence $\varphi$ and two different $\mathbf{F}$-stopping times
$\tau\neq \tau'$, the type-I errors of the sequential tests $(\varphi, \tau)$
and $(\varphi, \tau')$ are very different; namely,
$\alpha_0(\varphi, \tau) \neq \alpha_0(\varphi, \tau')$. Intuitively, for fixed
$\varphi$, the stopping rule that would inflate the most the type-I error
probability is $\tau = \inf\{t \in \mathbf{N} : \varphi_t = 1\}$, the rule that
stops as soon as possible. We say that a test is anytime valid if it has a
type-I error guarantee precisely under this worst-case rule.
\begin{definition}[Level-$\alpha$, Anytime-Valid (AV) test]
  \label{def:av-test}
  Let $\alpha\in (0,1)$, let $\varphi = (\varphi_t)_{t\in \mathbf{N}}$ be an
  $\mathbf{F}$-adapted sequence of decision functions, and define the stopping
  time $\tau = \inf \{t\in \mathbf{N}: \varphi_t = 1\}$. We say that the
  sequential test $(\varphi, \tau)$ is anytime valid at level $\alpha$ if
  $\alpha_0(\varphi, \tau)=\mathbf{P}_0\{\varphi_\tau = 1\} \leq \alpha$.
\end{definition}
The claim that, given a fixed decision sequence $\varphi$, the stopping rule
$\tau = \inf\{t: \varphi_t = 1\}$ is the worst case in the type-I error sense is
now verified. Indeed, if $(\varphi, \tau)$ is anytime valid at level $\alpha$,
then $\sup_{\omega}\alpha_0(\varphi, \omega)\leq \alpha$, where the supremum is
over all stopping times~\citep[][Lemma~1]{ramdas_admissible_2020} (justifying
the anytime-valid name). This uniform type-I error guarantee over stopping rules implies that a single decision sequence can accommodate the requirements of potentially multiple users with different preferences about power and sample
sizes~\citep{johari_always_2022}, it allows for reliable post-hoc inference when
the stopping rule is unknown~\citep{grunwald_safe_2024}, or when its data
dependency is complex~\citep{howard_uniform_2018}. Our focus is on AV tests and
the properties of the worst-case stopping time $\tau$.

By far, the most common and, in a sense, the only admissible way to build AV
tests is monitoring test martingales in
time~\citep{ramdas_admissible_2020}.
\begin{definition}[Test (super)martingale]
  An $\mathbf{F}$-adapted $\mathbf{P}_0$-(super)martingale $(E_t)_{t\geq 1}$ is
  a test (super)martingale if it is nonnegative and it starts at one; that is,
  $E_1 = 1$.
\end{definition}

Given both a test martingale $(E_t)_{t\geq 1}$ and a desired level
$\alpha \in (0, 1)$, we can construct an AV test $(\varphi, \tau)$ as
$\varphi_t = \indicator{E_t \geq 1 / \alpha}$ and
$\tau = \inf\{t : E_t \geq 1/\alpha\}$. By Ville's
inequality~\citep{ville_etude_1939}, an instance of Doob's martingale
inequality, this test is anytime valid with level at most $\alpha$; that is,
$ \alpha(\varphi, \tau) = \mathbf{P}_0\{\sup_{t\geq 1}E_t \geq 1/\alpha\} \leq
\alpha$.

In the case that the null hypotheses $\mathcal{H}_0$ is simple, the previous
test-martingale construction is admissible in the following sense: for every
level-$\alpha$ AV test, a test that stops more often and has the same level can
be built using this construction~\citep[][Theorem~18]{ramdas_admissible_2020}.
However, for a given test martingale $(E_t)_{t\geq 1}$, the associated AV test
$(\varphi, \tau)$ may be conservative (its type-I error might be strictly
smaller than $\alpha$) because Ville's inequality is not an equality in general.
To this effect, there are two possible methods to ``calibrate'' an AV test
rejecting ${\cal H}_0$ based on large values of $E_t$. The first method is to
use a smaller threshold to ensure that the level is exactly $\alpha$. This
method is particularly well suited for our case, where the null hypothesis
in~(\ref{eq:generic-alternative}) is simple and the threshold needed for exact
$\alpha$ type-I error can be estimated through Monte Carlo simulation (see
Section~\ref{sec:simulations}). The second method is to modify the
(super)martingale to guarantee that the overshoot $E_\tau - 1/\alpha$ is zero
$\mathbf{P}_0$-almost surely~\citep{fischer_improving_2025} in which case
Ville's inequality is an equality and the AV test has exact level $\alpha$. This
method is discussed in more detail in Section~\ref{sec:discussion}. For
simplicity, the generic threshold $1/\alpha$ is used when presenting the
theoretical results of Section~\ref{sec:results} on the power properties of the
AV tests, even if the ensuing tests might be slightly conservative. (In fact, as
we will see, those results are independent of the threshold.)

The analysis of the power properties of AV tests is divided in two parts: a
stopping-time--independent analysis that focuses on properties of the underlying
test martingales, and a stopping-time--specific analysis for the worst-case
stopping rule.

In order to perform a stopping-time--independent evaluation of the power
properties of admissible AV tests, the prevalent approaches in the literature
characterize the growth properties of the underlying test
martingales~\citep{ramdas_game-theoretic_2023}. The rationale is that if data
are generated under an element
$\mathbf{P}_{\varepsilon^\star, \delta^\star}=: \mathbf{P}^\star$ of the
alternative hypothesis $\mathcal{H}_1$ in (\ref{eq:generic-alternative}), the
$\mathbf{P}_0$-test martingale should lose the martingale property and it
typically grows exponentially fast as time passes. For a given
$\mathbf{P}_0$-test martingale $(E_t)_{t\in\mathbf{N}}$, the growth rate under
the alternative is quantified by $\mathbf{E}^\star[\ln E_t]$, where
$\mathbf{E}^\star$ denotes the expectation under $\mathbf{P}^\star$. Hence, the
presence of the logarithmic function is unrelated to ``the utility'' (which
could be could be measured using another function) provided by $E_t$; instead,
it is related to its potential exponential
growth~\citep[see][]{kelly_jr_new_1956}. We use the alternative-expected
logarithmic value as a measure of optimality for test martingales; we call
log-optimal a process that maximizes it.
\begin{definition}[Log-optimal test martingale]\label{def:log-opt-mart} Let ${\cal E}$ be the set of all
  $\mathbf{F}$-test martingales. A test martingale $(E^\star_t)_{t\geq 1}$ is
  log-optimal under $\mathbf{P}^\star$ if, for any fixed stopping time $\sigma$,
  it maximizes $(E_t)_{t \geq 1}\mapsto \mathbf{E}^\star[\ln E_\sigma]$ over
  test martingales; that is, if
  $\mathbf{E}^\star[\ln E^\star_\sigma] = \sup_{E \in {\cal
      E}}\mathbf{E}^\star[\ln E_\sigma]$.
\end{definition}
In the case that $\delta^\star$ and $\varepsilon^\star$ are specified, the
alternative hypothesis is simple and, since $\mathcal{H}_0$ is also simple, the log-optimal test martingale
$(E^\star_t)_{t\geq 1}$ is the likelihood ratio $E^\star_t := p_{\varepsilon^\star,
    \delta^\star}(X_1, \dots, X_t) / p_0(X_1, \dots,
  X_t)$, where $p_{\varepsilon^\star,
  \delta^\star}$ is the density of $\mathbf{P}^\star$ and
$p_0$ is the density of
$\mathbf{P}_0$~\citep[][Theorem~8]{koolen_log-optimal_2022}. When the property
defining log-optimality is required at every fixed
$t$ instead of every stopping time
$\sigma$, the result is classic~\citep{kelly_jr_new_1956,breiman_optimal_1961}.
In both forms, the result that the likelihood ratio is log-optimal is an
analogue of Neyman--Pearson's theorem; it implies that the AV test with decision
sequence determined by $\varphi^\star_t = \indicator{E^\star_t \geq
  1/\alpha}$, which uses the log-optimal test martingale
$E^\star_t$, is actually an (approximate) one-sided---also called open-ended or
power-one---Sequential Probability Ratio Test
(SPRT)~\citep{lai_power-one_1977,robbins_expected_1974}.
The random sample size of this test is
$\tau^\star = \inf\{t : E^\star_t \geq 1 / \alpha\}$, and its expected value
$\mathbf{E}^\star[\tau^\star]$ is, in a certain sense, minimal over all
sequential tests against $\mathbf{P}^\star$ with the same type-I
error~\citep[][Section~3.4]{tartakovsky_sequential_2014}. Note that since the alternative
${\cal H}_1$ in~(\ref{eq:generic-alternative}) is not i.i.d., the Wald--Wolfowitz
optimality result for the SPRT~\citep{wald_optimum_1948} does not apply.

Our primary goal is to design AV tests that does not require knowledge of the
specific element $\mathbf{P}^\star$ that generated the data and to assess its
optimality. Given an AV test $(\varphi, \tau)$ that uses an underlying test
martingale $(E_t)_{t\geq 1}$ and no knowledge of $\mathbf{P}^\star$, we
divide this assessment in two parts: a stopping-time--independent part that only
looks at the log-optimality properties of the test martingale $(E_t)_{t\geq 1}$,
and a stopping-time--specific analysis of
$\tau_E = \inf\{t : E_t \geq 1 / \alpha\}$. The task is therefore to design
$(E_t)_{t\geq 1}$ so that $\mathbf{E}^\star[\ln E_{t}]$ and $\tau_E$ are as
close to possible to the theoretically optimal
$\mathbf{E}^\star[\ln E_{t}^\star]$ and $\tau^\star$. To achieve this, it is
necessary to evaluate the latter, oracle quantities.

\subsection{Evaluating AV tests for sparse alternatives}
\label{sec:framework}

Despite the uncomplicated appearance of~(\ref{eq:generic-alternative}),
quantifying the optimal values $\mathbf{E}^\star[\ln E_{t}^\star]$ and
$\tau^\star$ turns out to be a delicate task. Firstly, existing literature on AV
tests focuses, for a given test martingale $(E_t)_{t\geq 1}$, on either the
finite-time regret
${\cal R}_t := \mathbf{E}^\star[\ln E_{t}^\star] - \mathbf{E}^\star[\ln E_{t}]$
or its large-time rate
$\lim_{t\to \infty}{\cal R}_t /
t$~\citep{ramdas_game-theoretic_2023,wang_anytime-valid_2025,grunwald_safe_2024,vovk_nonparametric_2024}.
In our case, since the observations are not independent across streams under
${\cal H}_1$ in~(\ref{eq:gaussian-alternative}), the value
$\mathbf{E}^\star[\ln E_{t}^\star]$ has no interpretable closed-form expression,
making it hard to evaluate ${\cal R}_t$ for each fixed $t$. Furthermore, taking
the limit $\lim_{t\to \infty}{\cal R}_t / t$ is not satisfactory either: as $t$
grows larger, the signal available on the anomalous streams grows as
$\sqrt{t}\delta$ and the hypotheses become well separated eventually
(see~Section~\ref{sec:state-things}). For this reason, it is not hard to design
test martingales $(E_t)_{t\geq 1}$ so that
$\lim_{t\to \infty}{\cal R}_t / t = 0$. Regarding the evaluation of the stopping
time $\tau^\star$, many existing results are not applicable to our problem
because they focus on the i.i.d. case. When data are not i.i.d., to the best of
the author's knowledge, all existing approximations for distributional
properties of $\tau^\star$, like its expected value
$\mathbf{E}^\star[\tau^\star]$ (or higher-order moments), rely on taking the
limit of small $\alpha$, which, a posteriori, is also a large-sample
approximation. For instance, if
$\lim_{t\to\infty} t^{-1} \mathbf{E}^\star[\ln E^\star_t] = r$ for a
positive $r$, it is known that
$\mathbf{E}^\star[\tau^\star] = \ln(1 / \alpha) / r + o(\ln(1 / \alpha))$ as
$\alpha\to 0$~\citep[][Section~3.4]{tartakovsky_sequential_2014}, an
approximation for large average sample sizes.

In a nutshell, for our setting, a conventional analysis cannot be performed at
the relevant sample sizes: for finite $t$, the quantities involved are not
tractable analytically, and taking a large-time, asymptotic approach would be
outside of the sample-size regime that is relevant to the problem. In order to
circumvent these difficulties, the strategies used in the fixed-sample
literature strongly suggest that an asymptotic, large-$K$ stance might be
productive (see Section~\ref{sec:state-things}).

In a large-$K$ analysis, a tantalizing question is whether the fixed-sample
power transition observed for sparse alternatives with
$\varepsilon^\star =\varepsilon^\star(K) \to 0$ and
$\delta^\star = \delta^\star(K) \to \infty$ (see Section~\ref{sec:state-things})
also has a sequential analogue. The existing results strongly suggest that
sequences $\varepsilon^\star = K^{-\beta^\star}$ and
$\delta^\star \propto \sqrt{2\ln K}$ are candidates for evaluating both
$(\mathbf{E}^\star[\ln E_t^\star])_{t\geq 1}$ and $\tau^\star$ (defined in
Section~\ref{sec:explain-av}) in a large-$K$ analysis. Note that the sparsity
parameter $\varepsilon^\star$ is unaffected by the passage of time, but this is
not the case for the signal strength: the signal available in the anomalous
streams, should there be any, accumulates as $\sqrt{t}\delta^\star$. Hence, the
fixed-sample question of \emph{whether} testing is possible becomes a question
of \emph{when} it becomes possible. This suggests that, in order to fix a time
scale for the sequential problem, it is useful to parametrize $\delta^\star$ in
terms of a time scale $T^\star$ (unknown to us) such that
$\delta^\star \propto 1 / \sqrt{T^\star}$. An educated guess for a large-$K$
analysis is hence of the form $\delta^\star = \sqrt{2(1 / T^\star)\ln K}$, where
$T^\star$ is an arbitrary reference time, and the fraction $t/T^\star$ in
$\sqrt{t}\delta^\star = \sqrt{2(t/T^\star)\ln K}$ plays the role of an
``information fraction'' that controls the easiness of the problem, a role
played by $r^\star$ in Section~\ref{sec:state-things}.

Although not clear in principle---this is discussed as the results are presented
in full in Section~\ref{sec:results}---, this parametrization turns out to be
the adequate one. Furthermore, since every AV test can be used to build a
fixed-sample test, in informal terms, a reasonable conjecture is that stopping
using an AV test should not be possible (in this large-$K$ sense) for
$t < t^\star := T^\star \rho(\beta^\star)$, because, in that case,
$r^\star = t / T^\star$ is below the fixed-sample ``detection boundary''
$\rho(\beta^\star)$ from~(\ref{eq:detection_boundary}) in
Section~\ref{sec:state-things}. However, for $t \in (t^\star, \sqrt{2\ln K})$,
when the signal is above the fixed-sample ``detection boundary'' but below the
identification boundary (see the discussion around Figure~\ref{fig:regions}), the
existence of powerful fixed-sample tests in the large-$K$ regime has, in
principle, no direct implications for AV testing (see Section~\ref{sec:results}
for an in-depth discussion). In this regime lies the main contribution of this
work, where both the large-$K$ growth properties of test martingales and
distributional properties of type-I--worst-case stopping times are
evaluated.

\subsection{Summary of the main results}
\label{sec:contribution}

We study the problem of testing the null hypothesis ${\cal H}_0$ that
$\varepsilon = 0$ in (\ref{eq:generic-alternative}) against the alternative
${\cal H}_1$ of sparse (small-$\varepsilon$) anomalies. When considering a
simple alternative hypothesis specified by a set of parameters
$(\varepsilon^\star,\delta^\star)$, the log-optimal test martingale is the
likelihood ratio process $(E^\star_t)_{t\geq 1}$. In terms of the joint density
of a standard normal distribution with mean $\delta$ and unit variance,
$p_{\delta}(x_1,\dots, x_t) = \prod_{s\leq t}\{ (1/\sqrt{2\pi})\exp\paren{-(x_s
  - \delta)^2/2}\}$, the likelihood ratio $E^\star_t$ has the form
\begin{equation}\label{eq:optimal_estat_simple_alt}
  E^\star_t
  :=
  \prod_{i \in [K]}\bracks{
    (1 - \varepsilon^\star) + \varepsilon^\star\frac{p_{\delta^\star}(X_{i,1},\dots,
      X_{i,t})}{p_{0}(X_{i,1},\dots, X_{i,t})}
  }.
\end{equation}
While the tests that we study are valid at any fixed $t$ and $K$, we take a
fixed-$t$ and large-$K$ stance in the analysis (see
Section~\ref{sec:framework}). Our results show that if
$\varepsilon^\star = K^{-\beta^\star}$ for $\beta^\star\in (1/2, 1]$, then, for
each sequence $\delta^\star = \sqrt{2(1 / T^\star)\ln K}$ (see
Section~\ref{sec:state-things}), there is indeed a dichotomy at
$t^\star := T^\star \rho(\beta^\star) = (2\rho(\beta^\star)\ln K) /
\delta^{\star 2}$ for the oracle test martingale: asymptotically as
$K\to\infty$, the maximum possible logarithmic growth transitions from trivial
to infinite. Moreover, the event that the ensuing AV test stops goes from
impossible to unavoidable. This work characterizes this dichotomous behavior
around the time point $t^\star$ and its contributions are threefold.
\begin{enumerate}
\item Theorem~\ref{thm:lr_achieves_boundary} shows that if $t < t^\star$, then
  no test martingale can have nontrivial logarithmic growth rate (this is stated
  in Corollary~\ref{cor:log-impossible}), while if $t > t^\star$, the log-optimal
  test martingale, the likelihood ratio~(\ref{eq:optimal_estat_simple_alt}), has
  logarithmic growth tending to infinity as $K\to\infty$.
\item Theorem~\ref{thm:stop-on-time} shows that, under the alternative
  hypothesis, the one-sided SPRT using the stopping rule
  $\tau^\star = \inf\{t : E^\star_t \geq 1/\alpha \}$ stops at around time
  $t^\star$ with $\mathbf{P}^\star$-probability tending to one.
\item We build a test martingale---a mixture likelihood ratio---that can be
  computed efficiently in Section~\ref{sec:prior}. The construction of this
  martingale requires no knowledge of the specific $\mathbf{P}^\star$ that
  generates the data; it is the base of our adaptive AV test
  for~(\ref{eq:generic-alternative}). Theorem~\ref{thm:prior-works} shows that
  the logarithmic growth rate of the proposed test martingale exhibits the same
  behavior as that of the oracle likelihood ratio.
  Theorem~\ref{thm:prior-stops-on-time} shows that a mixture SPRT that
  uses the proposed statistic also stops at around $t^\star$ with
  $\mathbf{P}^\star$-probability tending to one.
\end{enumerate}
Additionally, the finite-$K$ performance of our methods is demonstrated through
numerical simulations in Section~\ref{sec:simulations}.

\begin{figure}[ht!]
  \centering
  \includegraphics{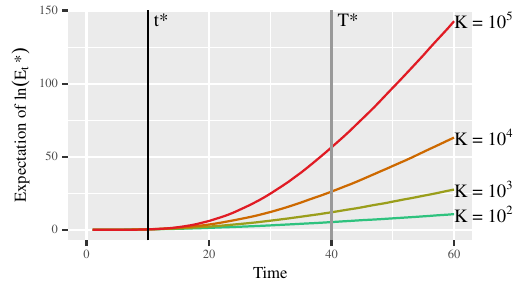}
  \caption{Alternative-expected logarithmic value of the likelihood ratio under
    $(\varepsilon^\star(K), \delta^\star(K))$. For increasing K, and fixed
    $\beta^\star= 0.75$ and $T^\star = 40$. See Section~\ref{sec:contribution} and
    Remark~\ref{rmk:expected_logE}.}
  \label{fig:expected_log_E}
\end{figure}

\begin{remark}[Behavior in time]\label{rmk:expected_logE}
  Figure~\ref{fig:expected_log_E} illustrates the expected log-growth
  $\mathbf{E}^\star[\ln E_{t}^\star]$ as a function of time $t$.
  The figure underscores the subtlety of the anomaly detection problem in the
  time period between $t^\star$ and $T^\star$, when $r = t/ T^\star$ is in the
  Detection Region in Figure~\ref{fig:regions} and anomalies are strong enough
  for detection but not for identification (see also
  Section~\ref{sec:state-things}).
  Arguably, this is the most interesting time interval. In this interval, the
  growth of $\mathbf{E}^\star[\ln E_{t}^\star]$ in time is nonlinear and slower
  than it would be if the anomalous streams were known; however, our results
  show that it is still enough to enable early stopping in the large-$K$ limit.
  For $t>T^\star$, the expected value $\mathbf{E}^\star[\ln E_{t}^\star]$ grows
  approximately linearly in $t$ with slope close to
  $K\varepsilon^\star \delta^{\star 2}/2$, the value that it would have if we
  already knew where the anomalies are. This is expected when $t>T^\star$
  because, in that case, the anomalous streams can be identified reliably.
\end{remark}

\subsection{Notation}
\label{sec:notation}
This subsection describes notational choices used in the main body of the
article.
\paragraph*{Distributions} We use $\mathbf{P}^\star$ when referring to the
distribution $\mathbf{P}_{\varepsilon^\star, \delta^\star}$ of the data under a
pair of true alternative parameters $\varepsilon^\star$ and $\delta^\star$, and
$\mathbf{P}_0$ for the null hypothesis, that is any distribution in
(\ref{eq:gaussian-alternative}) with either $\varepsilon = 0$ or $\delta=0$. We
call $\mathbf{E}^\star$ the expectation operator under $\mathbf{P}^\star$; and
$\mathbf{E}_0$, under $\mathbf{P}_0$.

\paragraph*{Asymptotic notation} Let $f$ and $g$ be two real-valued functions,
$f,g: \mathbf{N}\to \mathbf{R}$. We write that $f = O(g)$ as $n\to \infty$ and
say that $f$ is in the big order of $g$ if there is $n_0\in \mathbf{N}$ and
$C>0$ such that $|f(n)| \leq C|g(n)|$ anytime that $n\geq n_0$.

\subsection{Outline}
\label{sec:outline}
The rest of this article is organized as follows. Section~\ref{sec:results}
contains the main theoretical results summarized in
Section~\ref{sec:contribution}. In Section~\ref{sec:simulations}, the finite-$K$
performance of the ensuing tests is demonstrated through computer simulations
and it is compared to what would be achieved with fixed-sample tests. In the
last section, Section~\ref{sec:discussion}, we recount the main results and
discuss open problems; in particular, in relation both the ensuing fixed-sample
test and to the computational complexity of the methods that are presented here.

Four appendices accompany this article: Appendix~\ref{app:e-max-test} shows an
instructive failed attempt that the authors performed while trying to derive the
adaptive AV tests contained in this article, Appendix~\ref{app:proofs} contains
the proofs that are omitted from the main text,
Appendix~\ref{app:comparison-bonferroni} presents additional simulations that show
the superiority of this work's approach over existing standard techniques to
construct AV tests, and Appendix~\ref{app:simulation-details} gives further
details on the simulations of both Section~\ref{sec:simulations} and
Appendix~\ref{app:comparison-bonferroni}.

\section{Results}
\label{sec:results}

This section contains the main theoretical results in this article.
Section~\ref{sec:oracle} contains Theorem~\ref{thm:lr_achieves_boundary}, characterizing the logarithmic growth of the optimal test martingale $E_t^\star$; and
Theorem~\ref{thm:stop-on-time}, about the distributional properties of the
associated one-sided SPRT,~ $\tau^\star$. Section~\ref{sec:prior} contains a
construction of an adaptive test martingale (a mixture likelihood ratio) and its characterization, specifically Theorem~\ref{thm:prior-works}, on the log-optimality properties of the construction; and Theorem~\ref{thm:prior-stops-on-time}, on the sample size of the associated mixture SPRT.

\subsection{Oracle test martingale}
\label{sec:oracle}

If $\varepsilon^\star$ and $\delta^\star$ are known, the log-optimal test
martingale is the likelihood ratio $E^\star_t$ in
(\ref{eq:optimal_estat_simple_alt}). Thus, analyzing the logarithmic expected
value of $E^\star_t$ under well-specified $\mathbf{P}^\star$, we find the
barrier of possibility for test martingales with oracle knowledge of the true
alternative. The next result, proven in
Section~\ref{sec:proof-lr_achieves_boundary}, provides a zeroth order
characterization.
\begin{theorem}\label{thm:lr_achieves_boundary}
  Fix $t\geq 1$. Let $\beta^\star\in (1/2, 1]$, and let $T^\star > 0$. Use the
  parametrization
  $\delta^\star = \delta^\star(K) = \sqrt{2(1 / T^\star)\ln K}$ and
  $\varepsilon^\star = \varepsilon^\star(K) = K^{-\beta^\star}$. If
  $E^\star_t = E^\star_t(\varepsilon^\star, \delta^\star) $ is the likelihood
  ratio in~(\ref{eq:optimal_estat_simple_alt}) and $\rho$ is the ``detection
  boundary'' from~(\ref{eq:detection_boundary}) and
  $t^\star = T^\star\rho(\beta^\star)$, then, as $K\to\infty$,
  \begin{equation*}
    \mathbf{E}^\star[\ln E^\star_t]
    \to
    \begin{cases}
      0 &\text{ if } \   t < t^\star;\\
      \infty &\text{ if } \  t > t^\star,
    \end{cases}
  \end{equation*}
  where $\mathbf{E}^\star$ is the expectation under $\mathbf{P}^\star$, the
  distribution with true alternative parameters $\varepsilon^\star$ and
  $\delta^\star$.
\end{theorem}
In practice, when $K$ is finite,
Theorem~\ref{thm:lr_achieves_boundary} (jointly with
Theorem~\ref{thm:stop-on-time} below) helps to connect a given signal-strength--sparsity
pair $(\delta^\star, \varepsilon^\star)$ to a sample size $t^\star$ from which
detection is possible. This is useful, for instance, if a particular target pair
$(\delta^\star, \varepsilon^\star)$ is of interest, and we would like to
determine the necessary time elapsed until detection becomes possible under that
alternative. Indeed, if $\delta^\star > 0$ and $\varepsilon^\star > 0$ are
fixed, we can first define both $T^\star = (2\ln K)/\delta^{\star 2}$ and
$\beta^\star = \ln(1 / \varepsilon^\star)$, and use
Theorem~\ref{thm:lr_achieves_boundary} to conclude that, approximately, it would
take up to time $t^\star \approx (2\rho(\beta^\star)\ln K) / \delta^{\star 2}$
for any test martingale to have nontrivial logarithmic growth. As is intuitive,
$t^\star$ becomes larger either when the signal strength becomes smaller (if
$\delta^\star$ decreases), and it becomes larger if there are less anomalous
streams on average (if $\beta^\star$ increases).

Theorem~\ref{thm:lr_achieves_boundary} is related but not implied by existing
fixed-sample results. The proof of Theorem~\ref{thm:lr_achieves_boundary} is, at
its core, an analysis of the Kullback-Leibler (KL) divergence between the null
and the alternative hypotheses (the KL divergence is the alternative-hypothesis
expected logarithmic value of the likelihood ratio process). Because of known
relations that exist between statistical risk and the KL
divergence~\citep[Chapter~2]{tsybakov_introduction_2009}, it is possible to
conclude, in a certain sense (see the proof of
Theorem~\ref{thm:lr_achieves_boundary}) that if asymptotic--power-one,
fixed-sample tests exist, then the logarithmic growth of the likelihood ratio
must tend to infinity. However, the fixed-sample impossibility results that
preclude powerful testing when the signals are too weak (under the ``detection
boundary'') do not have, a priori, any direct implication on the existence
of test martingales with nontrivial asymptotic expected logarithmic growth; this
turns out to be the case a posteriori.

Additionally, known tools of analysis developed for fixed-sample
likelihood-ratio tests such as ``the second-moment
method''~\citep{addario-berry_combinatorial_2010} are not directly applicable in
our situation and additional technical work is needed.

Theorem~\ref{thm:lr_achieves_boundary} implies an impossibility result for any
other test martingale. Indeed, since the likelihood ratio is log-optimal in the
sense of Definition~\ref{def:log-opt-mart},
Theorem~\ref{thm:lr_achieves_boundary} implies that no test martingale can have
a nontrivial asymptotic growth rate before the ``detection moment'' $t^\star$.
\begin{corollary}\label{cor:log-impossible}
  Fix $t\geq 1$ and let $(E_t)_{t\geq 1}$ be a test martingale
  for~(\ref{eq:generic-alternative}). Let
  $\delta^\star, \varepsilon^\star, T^\star$ and $\beta^\star$ as in
  Theorem~\ref{thm:lr_achieves_boundary}. Then, for
  $t^\star = T^\star\rho(\beta^\star)$,
  \begin{equation*}
    \mathbf{E}^\star[\ln E_t] \to 0 \text{ as $K\to\infty$ if } t < t^{\star}.
  \end{equation*}
\end{corollary}
\begin{proof}
The claim follows from Theorem~\ref{thm:lr_achieves_boundary} because
$\sup_{E \in {\cal
      E}}\mathbf{E}^\star[\ln E_t] = \mathbf{E}^\star[\ln E_t^\star]$, where
$E_t^\star$ is the likelihood ratio (\ref{eq:optimal_estat_simple_alt}) and the
sup is over test supermartingales (see the discussion after
Definition~\ref{def:log-opt-mart}).
\end{proof}

If we interpret the value of a test martingale as a measure of ``evidence''
against $\mathcal{H}_0$, Corollary~\ref{cor:log-impossible} shows that gathering
this evidence is impossible in this precise log-optimality asymptotic sense
when it is ``too soon'' (measured by the relation of $t$ with respect to
$t^\star = T^\star\rho(\beta^\star)$). The question of what this implies for
sequential decision making at a given level remains unanswered based on the
results presented so far; an answer is given next in
Theorem~\ref{thm:stop-on-time}, which studies the sample size of the one-sided
SPRT~\citep{wald_sequential_1947,siegmund_sequential_1985}, the test that stops
at $\tau^\star = \min\{t: E_t^\star \geq 1 / \alpha\}$ (see
Section~\ref{sec:explain-av}).

The following theorem, proven in Section~\ref{sec:proof-stop-on-time}, shows
that, with $\mathbf{P}^\star$-probability tending to one, the one-sided SPRT has
sample size $\tau^\star \leq \lfloor t^\star\rfloor +1$, where $t^\star$ is as
in Theorem~\ref{thm:lr_achieves_boundary} and
$\lfloor x \rfloor := \max\{n\in\mathbf{Z}: n \leq x\}$ for $x\in\mathbf{R}$ is
the floor function. Hence, the moment at which
the test martingale $E_t^\star$ has infinite asymptotic logarithmic growth
coincides in the same limit with the moment at which an SPRT that monitors
$E^\star_t$ would stop with certainty.
\begin{theorem}\label{thm:stop-on-time}
  Adopt the assumptions of Theorem~\ref{thm:lr_achieves_boundary}. Let
  $\alpha\in(0,1)$, define
  $\tau^\star := \inf\{t : E_t^\star \geq 1 / \alpha\}$, and let
  $t^\star = T^\star\rho(\beta^\star)$. Then,
  \begin{equation*}
    \text{ as } K\to \infty, \
    \mathbf{P}^\star\{\tau^\star \leq t\} \to
    \begin{cases}
      0 & \text{ if } \  t < t^\star;\\
      1 & \text{ if }  \ t > t^\star.
    \end{cases}
  \end{equation*}
\end{theorem}
Despite the similar threshold behavior shown by both results, there is no
direct relation between Theorem~\ref{thm:lr_achieves_boundary} and
Theorem~\ref{thm:stop-on-time}. Indeed, for a fixed $t > t^\star$, showing that
$\ln E^\star_t$ diverges to infinity in probability as $K\to\infty$ would imply
Theorem~\ref{thm:stop-on-time}, but Theorem~\ref{thm:lr_achieves_boundary}
provides this information in expectation, not in probability. Hence, the
first-moment information given by Theorem~\ref{thm:lr_achieves_boundary} is, on
its own, not enough to prove Theorem~\ref{thm:stop-on-time}; it is also
necessary to study the concentration properties of $\ln E^\star_t$ around its
expected value. The proof of Theorem~\ref{thm:stop-on-time}, shown in
Section~\ref{sec:proof-stop-on-time}, achieves this goal through a precise
asymptotic analysis of the second moment of $\ln E^\star_t$ under
$\mathbf{P}^\star$.

Studying additional distributional properties of $\tau^\star$, such as its
moments, remains an interesting open direction.

\subsection{Adaptive methods and misspecification}
\label{sec:prior}

Two main families of methods exist to test the composite alternative
in~\eqref{eq:generic-alternative}; that is,  when the true value of the parameters
$(\varepsilon^\star, \delta^\star)$ is unknown: plug-in methods and mixture
methods. In this section, we outline the idea behind these two families of
methods. As we will wee, the plug-in approach has unclear optimality properties, whereas the mixture methods (when naively implemented) lead to statistics that are computationally intractable.
Fortunately, a specific instantiation of the mixture
method gives a statistic that is both tractable analytically and can be computed
efficiently. An empirical comparison between the two is provided in
Appendix~\ref{app:comparison-bonferroni}.

\paragraph*{Plug-in methods} The idea behind sequential plug-in methods is to
use predictable (with respect to the natural data filtration) point estimates of
the unknown parameters and to ``plug'' them into the likelihood ratio
statistic~\citep{waudby-smith_estimating_2023,grunwald_safe_2024,ramdas_game-theoretic_2023}.
While monitoring plug-in test martingales yields valid AV tests, it is not clear
whether strong power properties can be assured because of the assumed sparsity
of the alternative in~(\ref{eq:generic-alternative}). The main ingredient of
plug-in methods is a sequence
$((\hat{\varepsilon}_t, \hat{\delta}_t))_{t \geq 1}$ of $\mathbf{F}$-predictable
estimators for $(\varepsilon^\star, \delta^\star)$. A standard choice for
$(\hat{\varepsilon}_t, \hat{\delta}_t)$ is a smoothed maximum likelihood
estimator based on the data up to $t - 1$. Given this sequence, a plug-in test
martingale $(M_t)_{t\geq 1}$ is
\begin{equation}\label{eq:preq-plug-in-mart}
  M_t = \prod_{s\leq t} \frac{\hat{p}_s(X_s \ | \ X_1, \dots, X_{s-1})
  }{
    p_0(X_s \ | \ X_1, \dots, X_{s-1})
  },
\end{equation}
where, for each $t$,
$\hat{p}_t(x_t \ | \ x_1, \dots, x_{t-1}) := p_{\hat{\varepsilon}_t,
  \hat{\delta}_t}(x_t \ | \ x_1, \dots, x_{t-1})$ is the conditional density of
$X_t$ given $X_1, \dots, X_{t - 1}$ under the alternative
in~(\ref{eq:generic-alternative}) evaluated at the predictable estimators
$(\hat{\varepsilon}_t, \hat{\delta}_t)$. While fixed-$K$, large-$t$ properties
of plug-in test martingales are well
known~\citep{grunwald_minimum_2019,grunwald_minimum_2007}, the same cannot be
said in our large-$K$, fixed-$t$ sparse regime. The reason for this is that in
our sparse setting, where $\beta^\star \geq 1/2$ (and $\varepsilon^\star$ is of
lower order than $1 / \sqrt{K}$), to the best of the authors' knowledge, it is
not known whether consistent (or at least accurate enough) joint estimation of
$\varepsilon^\star$ and $\delta^\star$ is possible below the ``detection
boundary'' (see Section~\ref{sec:state-things}); only estimators for
$\varepsilon^\star$ are known in this
regime~\citep{meinshausen_estimating_2006,cai_estimation_2007} (and above the
detection boundary). For the case $\beta < 1/2$, which is outside of our current
interest, the picture is very different and joint estimation of
$\varepsilon^\star$ and $\delta^\star$ is possible in a minimax
sense~\citep{cai_optimal_2010} (see also Section~\ref{sec:state-things}).
Despite the challenges surrounding estimation of
$(\varepsilon^\star, \delta^\star)$, a test based
on~(\ref{eq:preq-plug-in-mart}) is anytime valid. We evaluate numerically a
plug-in test martingale based on a smoothed maximum likelihood estimator in
Appendix~\ref{app:comparison-bonferroni}, where it is shown that the plug-in
method consistently underperforms the upcoming adaptive mixture-based
construction.

\paragraph*{Mixture methods} The second main family of methods to test the
composite alternative in~\eqref{eq:generic-alternative} is the family of mixture
methods~\citep{robbins_statistical_1970}; the adaptive tests that we propose use
this method. Consider the likelihood ratio $E_t(\varepsilon,
\delta)$ 
against a generic alternative set of parameters $(\varepsilon, \delta)$, that
is,
\begin{equation}
  \label{eq:generic-lr}
  E_t(\varepsilon, \delta) := \prod_{i \in [K]}\bracks{ (1 - \varepsilon) +
    \varepsilon\frac{p_{\delta}(X_{i,1},\dots, X_{i,t})}{p_{0}(X_{i,1},\dots,
      X_{i,t})} }.
\end{equation}
For a given choice of a probability distribution $\mathbf{\Pi}$ on the
alternative-hypothesis parameters, called a prior or mixture distribution, define
the mixture likelihood ratio process $(E_t(\mathbf{\Pi}))_t$ as the mixture
$E_t(\mathbf{\Pi}) := \int E_t(\varepsilon, \delta)\rmd
\mathbf{\Pi}(\varepsilon, \delta)$. Since the process
$(E_t(\varepsilon,\delta))_{t\geq 1}$ is a test martingale with respect to the
natural data filtration for each pair $(\varepsilon,\delta)$, the mixture
$(E_t(\mathbf{\Pi}))_{t\geq 1}$ is a test martingale as well and it can be used
to build AV tests. Notice that even though we refer to $\mathbf{\Pi}$ as a
``prior'', no Bayesian posterior computations are made in this work.

The fact that $(E_t(\mathbf{\Pi}))_{t\geq 1}$ is a test martingale for any prior
$\mathbf{\Pi}$ guarantees anytime validity for the ensuing tests, but its power
properties depend crucially on $\mathbf{\Pi}$. In order to evaluate a choice of
prior $\mathbf{\Pi}$, we take the following fixed-$t$ and large-$K$ asymptotic
stance. Under the $\beta$-$r$ parametrization of
Theorem~\ref{thm:lr_achieves_boundary}, the true parameters
$\varepsilon^\star = \varepsilon^\star(K)$ and $\delta^\star = \delta^\star(K)$
correspond to a sequence of true distributions
$\mathbf{P}^{\star} = \mathbf{P}^{\star}(K)$; in other words,
$\mathbf{P}^\star(K)$ corresponds to the distribution
$\mathbf{P}_{\varepsilon^\star, \delta^\star}$ with parameters
$\varepsilon^\star = \varepsilon^\star(K) = K^{-\beta^\star}$ and
$\delta^\star =\delta^\star(K) = \sqrt{2(1 / T^\star)\ln K}$, as in
Section~\ref{sec:oracle}. We build a sequence of priors
$\mathbf{\Pi} = \mathbf{\Pi}(K)$ and evaluate the asymptotic logarithmic value
$\mathbf{E}^\star[\ln E_t(\mathbf{\Pi})]$ for fixed $t$ as $K\to\infty$, where,
recall, both $\mathbf{P}^\star$ and $\mathbf{\Pi}$ depend on $K$. From
Corollary~\ref{cor:log-impossible}, it follows that if
$t < t^\star$ (under the ``detection boundary''), then necessarily
$\mathbf{E}^\star[\ln E_t(\mathbf{\Pi})] \to 0$ as $K \to \infty$ for any choice
$\mathbf{\Pi} = \mathbf{\Pi}(K)$. Therefore, the goal is to construct mixture
distributions $\mathbf{\Pi} =\mathbf{\Pi}(K) $ such that, for
$\tau^\mathbf{\Pi} :=\inf \{E_t(\mathbf{\Pi}) \geq 1 / \alpha\}$, both
$\mathbf{E}^\star[\ln E_t(\mathbf{\Pi})] \to \infty$ and
$\mathbf{P}^\star\{\tau^{\mathbf{\Pi}} \leq t\} \to 1$ whenever
$t > t^\star$---just as having known the true parameters $\varepsilon^\star$ and
$\delta^\star$ (see Theorem~\ref{thm:lr_achieves_boundary} and
Theorem~\ref{thm:stop-on-time}).

Despite the existence of default (conjugate) priors for this problem, due to
computational considerations, we choose to use discrete priors instead. Indeed,
using, for example, a Gaussian distribution on $\delta^\star$ and a Beta distribution on
$\varepsilon^\star$ yields closed-form expressions for $E_t(\mathbf{\Pi})$, but
the resulting statistics cannot be computed efficiently. This inefficiency stems
from the fact that computing such mixture martingale requires expanding the
product that defines $E_t$ in (\ref{eq:generic-lr}) and, after integration with
respect to the mixture distribution, computing the resulting statistic would
need an exponential (in $K$) number of operations. This exponentially complex
computation is, as we will see, avoided using discrete mixture distributions.

In order to decide on which parameters $(\varepsilon,\delta)$ to place prior
mass, it is necessary to quantify how much mass needs to be placed close to
$(\varepsilon^\star, \delta^\star)$. The next proposition, which guides our choice of
prior, shows that if $(\varepsilon,\delta)$ and
$(\varepsilon^\star, \delta^\star)$ are ``close enough'' in an asymptotic sense
as $K\to\infty$, then $E_t(\varepsilon, \delta)$ has the same threshold behavior
as $E_t^\star$ at $t^\star$ (see Theorem~\ref{thm:lr_achieves_boundary} and
Theorem~\ref{thm:stop-on-time}). The
proof of the Proposition~\ref{prop:misspecification-cost} is given in
Section~\ref{sec:proofs-prior}. In proving
Proposition~\ref{prop:misspecification-cost}, we give explicit rates at which
$\mathbf{E}^\star[\ln E_t(\varepsilon,\delta)]$ grows to infinity as
$K\to\infty$; the main technical tool is Proposition~\ref{prop:master-rates}.
\begin{proposition}\label{prop:misspecification-cost}
  Let $\varepsilon^\star = \varepsilon^\star(K) = K^{-\beta^\star}$ and
  $\delta^\star = \delta^\star(K) = \sqrt{2(1/T^\star)\ln K}$ be the true
  parameters in $\mathbf{P}^\star$ corresponding to a fixed pair
  $(T^\star,\beta^\star)$. Let $((\beta_K, T_K))_{K\geq 1}$ be a sequence of
  parameters defining $\varepsilon = \varepsilon(K) = K^{-\beta_K}$ and
  $\delta = \delta(K) = \sqrt{2(1 / T_K)\ln K}$. Define
  $d_{\beta} = d_{\beta}(K) = \beta_K - \beta^\star$ and
  $d_{T} = d_{T}(K) = 1 - \sqrt{T_K / T^\star}$. Then, if
  $t^\star = T^\star\rho(\beta^\star)$,
  \begin{equation*}
    \text{ as } K\to\infty, \ \
    \mathbf{E}^\star[\ln E_t(\varepsilon, \delta)] \to
    \begin{cases}
      \infty &\text{ if } t > t^\star;\\
      0 &\text{ if } t < t^\star\\
    \end{cases}
  \end{equation*}
  whenever both $d_T(K)$ and $d_{\beta}(K)$ satisfy,
  respectively,
  \begin{equation}\label{eq:misspecification-cost}
    |d_{T}|K^{1 - \beta_K} = O(1)
    \text{ and }
    |d_{\beta}| \ln K  \to 0
    \text{ as }
    K\to\infty
    .
  \end{equation}
\end{proposition}
Written in terms of $\delta$ and $\varepsilon$, the quantities $d_{T}$ and
$d_{\beta}$ become $d_{T} = (\delta - \delta^\star) / \delta$, the relative
difference between $\delta$ and $\delta^\star$; and
$d_{\beta} = \ln(\varepsilon^\star / \varepsilon)$, the logarithmic ratio
between $\varepsilon^\star$ and $\varepsilon$. Then,
(\ref{eq:misspecification-cost}) reads
\begin{equation}\label{eq:distance-requirement}
  \left|\frac{\delta - \delta^\star}{\delta}\right| \varepsilon K \ln K = O(1)
  \text{ and }
  \left|\ln\paren{\frac{\varepsilon^\star}{\varepsilon}}\right| \ln K \to 0
  \text{ as }
  K\to\infty.
\end{equation}
Thus, both the relative difference in means multiplied by the expected number of
anomalies $\varepsilon K$, and the logarithmic ratio
$\ln(\varepsilon^\star /\varepsilon)$ must tend to zero sufficiently fast. Our
analysis shows that the quantities involved in (\ref{eq:misspecification-cost})
are the correct quantities to control. Indeed, our proofs could be adapted to
give not only sufficient, but also necessary conditions on the quantity involved
in (\ref{eq:misspecification-cost}) for $E_t(\varepsilon,\delta)$ to achieve the
detection moment $t^\star$ for any choice $\varepsilon = \varepsilon(K)$ and
$\delta = \delta(K)$. However, pursuing this would not advance our present
purposes since such a necessary condition would depend on the true parameters
$(\varepsilon^\star, \delta^\star)$, which are assumed to be unknown.

Proposition~\ref{prop:misspecification-cost} motivates the choice of a uniform
prior $\mathbf{\Pi} = \mathbf{\Pi}(K)$ on a grid
$\mathcal{G}_C = \mathcal{G}_C(K) \subset (0, 1]\times\mathbf{R}_+$ in the
parameters $(\varepsilon, \delta)$. Note that both $\mathbf{\Pi}$ and
${\cal G}_C$ depend on the number of streams $K$, but, to keep notation simple,
we keep this $K$ dependence implicit. The grid ${\cal G}_C$ is built so that
there is always a point $(\varepsilon, \delta)\in \mathcal{G}_C$ close enough to
$(\varepsilon^\star, \delta^\star)$ as described by
Proposition~\ref{prop:misspecification-cost}. In showing performance guarantees
for $(E_t(\mathbf{\Pi}))_t$ and the ensuing tests, it is assumed that there
exists $C > 0$ such that $T^\star < C$. This is not a stringent assumption as
$C$ might be determined from the maximum length that the experiment might take,
and, in the absence of such knowledge, a nonuniform prior can be placed on
arbitrarily small values of $\delta^\star$ (which corresponds to arbitrarily
large values of $T^\star$) and our analysis still applies with minor
modifications. We build first an exponential grid for $\varepsilon$ of size
$\lceil\ln^2 K\rceil$,
\begin{equation*}
  \mathcal{G}_\varepsilon
  =
  \bracks{
    K^{-\beta_i}:
    \beta_i =
    \frac{1}{2}
    +
    \frac{i}{2\lceil\ln^2 K\rceil}
    \ , \
    i = 1, \dots, \lceil \ln^2 K \rceil
  },
\end{equation*}
followed by a exponentially spaced $\varepsilon$-dependent grid
$\mathcal{G}_\delta(\varepsilon)$ on $\delta$ of size
$\lceil \varepsilon K\rceil \lceil \ln C \rceil$,
\begin{equation*}
  \mathcal{G}_{\delta, C}(\varepsilon)
  =
  \bracks{
    \sqrt{\frac{2\ln K}{\exp\paren{i/\lceil \varepsilon K \rceil}}}
     \ : \
    i = 1,  \dots, \lceil \varepsilon K\rceil\lceil \ln C \rceil
  }.
\end{equation*}
With $\mathcal{G}_\varepsilon$ and $\mathcal{G}_{\delta, C}$ at hand, consider
the set ${\cal G}_C$ given by
\begin{equation}\label{eq:grid}
  \mathcal{G}_C
  =
  \bigcup_{\varepsilon' \in \mathcal{G}_\varepsilon}
  \{(\varepsilon', \delta') : \delta' \in \mathcal{G}_{\delta, C}(\varepsilon')\}.
\end{equation}
We use $\mathbf{\Pi} = \mathrm{Uniform}({\cal G}_C)$ as prior on
$(\varepsilon,\delta)$. The nonuniform grid ${\cal G}_C = {\cal G}_C(K)$ (whose
size and shape depends on $K$) is designed so that, for any sequence of true
parameters $(\varepsilon^\star(K), \delta^\star(K))$, there is always a sequence
$((\varepsilon^\circ(K),\delta^\circ(K)))_{K\geq 1}$ of pairs
$(\varepsilon^\circ(K),\delta^\circ(K)) \in {\cal G}_C(K)$ that satisfies
(\ref{eq:distance-requirement}).

Using this grid construction avoids the need of exponential-time computations;
the statistic $E_t(\mathbf{\Pi})$ can be computed, up to logarithmic factors, in
$O(K^{3/2})$ steps at each $t$ as $K\to\infty$. A crude bound on the size of the sets already
reveals this fact. Indeed, notice that each
$\varepsilon\in {\cal G}_{\varepsilon}$ is of the form
$\varepsilon = K^{-\beta}$ with $\beta \in (1/2, 1]$ and consequently the size
of the $\varepsilon$-dependent grid $\mathcal{G}_{\delta, C}(\varepsilon')$ is
$O(K^{1/2})$. Since $|{\cal G}_{\varepsilon}| = \lceil \ln^2 K\rceil$, then,
$|\mathcal{G}_C| = \sum_{\varepsilon' \in
  \mathcal{G}_\varepsilon}|\mathcal{G}_{\delta, C}(\varepsilon')| =
O(K^{1/2}\ln^2 K)$ as $K\to\infty$. A more careful analysis (involving the sum of the actual sizes of each
$|\mathcal{G}_{\delta, C}(\varepsilon')|$) improves this bound by a logarithmic
factor; that is, $|\mathcal{G}_C| = O(K^{1/2}\ln K)$ as $K\to\infty$. Since computing
$E_t(\varepsilon, \delta)$ takes $O(K)$ operations for large $K$, computing
$E_t(\mathbf{\Pi})$ takes $O(K^{3/2}\ln K)$ operations in the same limit. Additionally, the
dependence of the computational complexity on the bound $C$ is logarithmic
because $|{\cal G}_C|\propto \lceil \ln C\rceil$, which we consider benign. For
ease of future reference, we summarize the proposed construction in a
definition.

We refer to the resulting test martingale as an adaptive test martingale for reasons that will become
apparent soon.
\begin{definition}[Adaptive test martingale]
  \label{def:adaptive-test-mart}
  Let $C > 0$ and let $t \geq 1$. Let ${\cal G}_C$ be the parameter grid
  from~\eqref{eq:grid} and let $\mathbf{\Pi} = \mathrm{Uniform}(\mathcal{G}_C)$.
  We define the adaptive test martingale $E_t(\mathbf{\Pi})$ as
  \begin{equation*}
    E_t(\mathbf{\Pi}) := \int E_t(\varepsilon,
    \delta)\rmd\mathbf{\Pi}(\varepsilon, \delta) =
    \frac{1}{|\mathcal{G}_C|} \sum_{(\varepsilon, \delta)\in \mathcal{G}_C}
    E_t(\varepsilon, \delta),
  \end{equation*}
  where $E_t(\varepsilon,\delta)$ is the likelihood ratio \eqref{eq:generic-lr}.
\end{definition}

For this construction, the following guarantee
holds (its proof is given in Section~\ref{sec:proofs-prior}).
\begin{theorem}\label{thm:prior-works}
  Let $C>1$, and let
  $\delta^\star = \delta^\star(K) = \sqrt{2(1 / T^\star)\ln K}$ and
  $\varepsilon^\star = \varepsilon^\star(K) = K^{-\beta^\star}$. Assume
  $T^\star < C$. Then, as $K\to\infty$ and for $1/2 <\beta^\star < 1$, the
  adaptive test martingale $(E_t(\mathbf{\Pi}))_{t\geq 1}$ from
  Definition~\ref{def:adaptive-test-mart} satisfies
  \begin{equation*}
    \mathbf{E}^\star[\ln E_t(\mathbf{\Pi})]
    \to
    \begin{cases}
      \infty &\text{ if } t > t^\star; \\
      0 &\text{ if }  t < t^\star,
    \end{cases}
  \end{equation*}
  where $t^\star = T^\star \rho(\beta^\star)$.
\end{theorem}
This implies that, through this zeroth order lens, using the adaptive test
martingale $(E_t(\mathbf{\Pi}))_{t\geq 1}$ from
Definition~\ref{def:adaptive-test-mart} guarantees the same performance as
having known the true parameters in advance (compare with
Theorem~\ref{thm:lr_achieves_boundary}). This justifies calling
$E_t(\mathbf{\Pi})$ adaptive. Whether the same performance guarantee can be
achieved using a statistic with lower computational complexity remains an open
question (see the discussion in Section~\ref{sec:discussion}). The proof of
Theorem~\ref{thm:prior-works} relies on the fact that, for any pair
$(\varepsilon,\delta) \in \mathcal{G}_C$,
\begin{equation}\label{eq:prior-balance}
  \mathbf{E}^\star[\ln E_t(\mathbf{\Pi})] \geq \mathbf{E}^\star[\ln
  E_t(\varepsilon,\delta)] - \ln(|\mathcal{G}_C|).
\end{equation}
The right hand side of~(\ref{eq:prior-balance}) shows that there is a trade-off
between the number of elements $|\mathcal{G}_C|$ in the support of the prior and
the ability to ensure at least one element is close enough to
$(\varepsilon^\star, \delta^\star)$. This balance is achieved by using the
proposed prior construction and Proposition~\ref{prop:misspecification-cost},
which gives the explicit rates at which
$\mathbf{E}^\star[\ln E_t(\varepsilon,\delta)]$ diverges for a fixed
$t> t^\star$ as $K\to\infty$ when the pair $(\varepsilon,\delta)$ is close
enough to $(\varepsilon^\star, \delta^\star)$.

Given that $E_t(\mathbf{\Pi})$ is a mixture likelihood ratio, the sequential
test defined by the stopping rule
$\tau^{\mathbf{\Pi}} := \inf\{t : E_t(\mathbf{\Pi}) \geq 1 / \alpha\}$ and the
decision to reject the null hypothesis if $\tau^{\mathbf{\Pi}} < \infty$, is an
(with level approximately $\alpha$) mixture SPRT~\citep{robbins_statistical_1970}. The next
proposition, whose proof is in Section~\ref{sec:proof-prior-stops}, shows that
this particular mixture SPRT stops asymptotically with certainty before
$t^\star$, just as the well specified (oracle) SPRT from
Section~\ref{sec:oracle}.
\begin{theorem}\label{thm:prior-stops-on-time}
  Adopt the assumptions of Theorem~\ref{thm:prior-works}. Let $\alpha\in(0,1]$,
  define $\tau^{\mathbf{\Pi}} := \inf\{t : E_t(\mathbf{\Pi}) \geq 1 / \alpha\}$,
  and let $t^\star = T^\star\rho(\beta^\star)$. Then,
  \begin{equation*}
    \text{ as } K\to \infty, \
    \mathbf{P}^\star\{\tau^{\mathbf{\Pi}} \leq t\} \to
    \begin{cases}
      0 & \text{ if } \  t < t^\star;\\
      1 & \text{ if } \ t > t^\star.
    \end{cases}
  \end{equation*}
\end{theorem}
The proof of Theorem~\ref{thm:prior-stops-on-time} is significantly more
involved than that of Theorem~\ref{thm:stop-on-time} and it implies the latter
result. As the two arguments use fundamentally different techniques, we keep both in the interest of the reader.

\section{Simulation study}
\label{sec:simulations}
Using computer simulation, we study the finite-sample performance of the AV
tests proposed in the previous sections. In
Appendix~\ref{app:comparison-bonferroni}, we compare these tests to prequential
plug-in test martingales, and a test based the higher-criticism statistic together with a Bonferroni-type correction. The assessment of our AV tests is done in two
ways: (1) by comparing the cumulative rejection rates of the AV tests against
the power of fixed-sample tests; and (2) by comparing other distributional
properties of the stopping times (their quantiles and their expected value).

The cumulative rejection rate can be interpreted as the power of a
time-truncated (mixture) SPRT, and it is sensible to compare it to the power of a
fixed-sample test. Indeed, for a given value $t$, we compare the rejection rate
of a (mixture) SPRT that is monitored continuously up to time $t$ to the power of a
fixed-sample test that is performed with the complete data at time $t$. We use
as fixed-sample benchmarks the (oracle) most powerful test and the
parameter-adaptive test based on the Higher Criticism ($\mathrm{HC}$) statistic,
a standard test for the fixed-sample problem. The cumulative rejection rate of
AV tests is expected to be lower (higher is better) than that of
fixed-sample tests because the latter do not account for multiple testing over
time.

Other distributional properties of the stopping time of the (mixture) SPRT
associated to our tests (see~Section~\ref{sec:explain-av}) are of interest as
well. Given a desired power, we compute the sample size required by fixed-sample
tests to achieve that power and compare it both to the expected stopping times of
the AV tests and their quantiles. We will see that more often than not, AV tests
stop before this sample size. In a nutshell, for a given power, the (mixture)
SPRT associated to the AV test for anomalies requires larger sample
sizes than the fixed-sample test in the worst case, but they stop earlier often.
This is in line with the findings in the
literature~\citep{tartakovsky_sequential_2014,siegmund_sequential_1985}.

In this section, we consider sample sizes required to obtain 80\% of the maximum
power achievable (see Section~\ref{sec:quantile-comparison}) for tests that
guarantee a type-I error rate of 5\%. Qualitatively, similar results hold for
other error rates.

\subsection{Cumulative rejection rate comparisons}
\label{sec:cumul-reject-rate}

We evaluate the cumulative rejection rate of the test that uses the log-optimal
statistic $E^\star_t$ from~(\ref{eq:optimal_estat_simple_alt}) and that of the
adaptive test martingale $E_t(\mathbf{\Pi})$ from
Definition~\ref{def:adaptive-test-mart}. Indeed, for $\alpha = 0.05$, and a pair
of true parameters $(\varepsilon^\star, \delta^\star)$, let
$\tau^\star = \inf\{t: E^\star_t\geq A^\star \}$ and let
$\tau^{\mathbf{\Pi}} = \inf\{t: E_t(\mathbf{\Pi})\geq A^{\mathbf{\Pi}} \}$,
where $A^\star$ and $A^{\mathbf{\Pi}}$ are the thresholds that guarantee that
the level of both tests equals $\alpha$; that is,
$\alpha = \mathbf{P}_0\{\sup_tE^\star_t \geq A^\star\} =
\mathbf{P}_0\{\sup_tE(\mathbf{\Pi}) \geq A^{\mathbf{\Pi}}\}$. Both $A^\star$ and
$A^{\mathbf{\Pi}}$ are estimated using Monte Carlo simulation (see the details
in Appendix~\ref{app:simulation-details}). Specifying $E_t(\mathbf{\Pi})$
requires the choice of a parameter $C$ (an upper bound on the duration of the
experiment), which is required to satisfy $C > 2\ln K / \delta^{\star 2}$.
Simulations (not presented) show that the impact of choosing a large $C$ (even
by several orders of magnitude above the theoretical requirement) has very small
effects in the observed results. We estimate the cumulative rejection rates,
defined as the cumulative distribution functions of $\tau^\star$ and
$\tau^{\mathbf{\Pi}}$; that is, the respective functions
$t\mapsto F_\tau^\star(t) := \mathbf{P}^\star\{\tau^\star \leq t\}$ and
$t\mapsto F_{\tau}^{\mathbf{\Pi}}(t) := \mathbf{P}^\star\{\tau^{\mathbf{\Pi}}
\leq t\}$. Given a fixed $t$, these cumulative distribution functions may be
interpreted as the power of a test that uses
$\tau^\star \wedge t := \min\{\tau^\star, t\}$ (or
$\tau^{\mathbf{\Pi}} \wedge t := \min\{\tau^{\mathbf{\Pi}}, t\}$) as a stopping
rule and rejects if $E_{\tau^\star \wedge t}^{\star} \geq A^\star$ (or if
$E_{\tau^{\mathbf{\Pi}} \wedge t}^{\mathbf{\Pi}} \geq A^{\mathbf{\Pi}}$,
respectively) . Because of this power interpretation, it is sensible to compare,
for a fixed $t$, the values of $F_\tau^\star(t)$ and
$F_{\tau}^{\mathbf{\Pi}}(t)$ to the power of fixed-sample test that analyzes a
snapshot of the data at time $t$. We include the most powerful test, the
likelihood ratio test based on $E^\star_t$, and the test that uses the higher
criticism statistic, a standard parameter-adaptive test in this situation.
Additionally to being parameter-adaptive in the sense of
Section~\ref{sec:state-things}, the higher criticism statistic $\mathrm{HC}_t$
is known to have good finite-sample power properties~\citep{donoho_higher_2015}.
The statistic $\mathrm{HC}_t$ can be computed as
\begin{equation}\label{eq:hc-statisic}
  \mathrm{HC}_{t} := \max_{i\in [K]}\frac{ \sqrt{K}(\hat{F}_{t}(U_{i,t}) - U_{i,
      t})}{\sqrt{U_{i,t}(1 - U_{i, t})}},
\end{equation}
where, for each $i\in [K]$, we write
$U_{i, t} = \mathbf{P}_0\{\mathrm{Normal}(0, 1) > Z_{i, t}\}$ for the Gaussian
complementary cumulative distribution function evaluated at the normalized stream sum
$Z_{i, t} = (1 / \sqrt{t})\sum_{s\leq t}X_{i, s}$---the null distribution of
$U_{i, t}$ is $\mathrm{Uniform}([0, 1])$---, and
$\hat{F}_{t}(u) := (1 / K)\sum_{i\in [K]}\indicator{ U_{i, t}\leq u}$ is the
empirical cumulative distribution function of $(U_{i, t})_{i\in [K]}$.
\begin{figure}[ht!]
  \centering
  \includegraphics{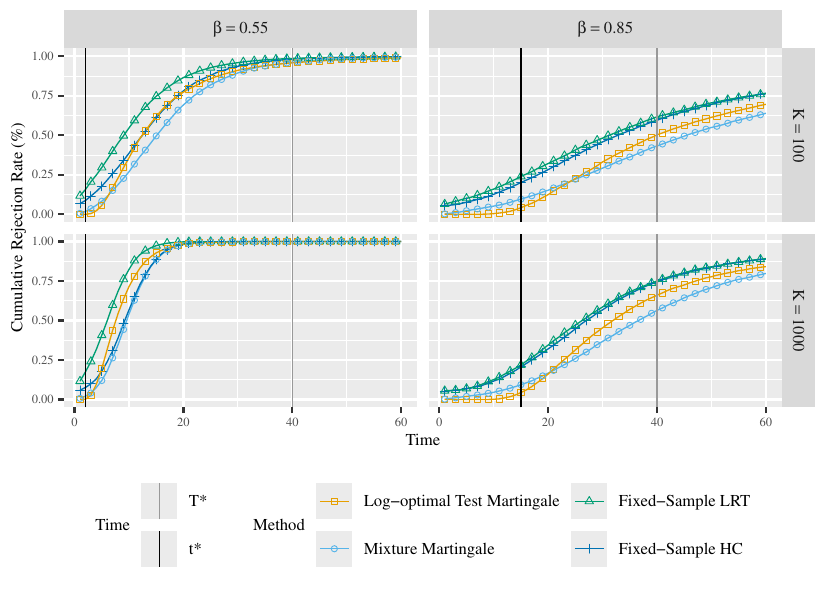}
  \caption{Cumulative rejection rate for AV methods are shown for the
    well-specified log-optimal test martingale
    (\ref{eq:optimal_estat_simple_alt}) and the adaptive construction from Definition~\ref{def:adaptive-test-mart} in
    Section~\ref{sec:prior}. The power is estimated for the likelihood ratio
    test (LRT) and a test based on the Higher Criticism (HC) statistic. Detection
    is possible between $T^\star$ and $t^\star$, when it is not always possible
    to estimate the anomalies. Only every other time point is marked to avoid
    overplotting. See more details in Section~\ref{sec:simulations}.}
  \label{fig:power}
\end{figure}
The result of the computer simulations is shown in Figure~\ref{fig:power}, and
the estimation details are found in Appendix~\ref{app:simulation-details}.
Estimates of the cumulative rejection rates of the (mixture) SPRT, given by
$t\mapsto F_\tau^\star(t)$ and $t\mapsto F_\tau^{\mathbf{\Pi}}(t)$ are shown for
the fraction of anomalies $\varepsilon = K^{-\beta^\star}$ with
$\beta^\star\in \{0.55, 0.85\}$, the number of streams $K\in \{100, 1000\}$, and
the signal strength $\delta^\star = \sqrt{2(1 / T^\star)\ln K}$ with
$T^\star = 40$. Changing the value of $T^\star$ only changes the scale of our
plots; we keep it constant. Furthermore, the corresponding moment $t^\star$,
where testing for anomalies in the asymptotic regime becomes possible; and
$T^\star$, when reliable identification is possible, are shown. The cumulative
rejection rate of the (mixture) SPRT is lower than the power of the fixed-sample
tests, it remains competitive for values of $\beta^\star$ close to $1 / 2$, and
deteriorates as $\beta^\star \to 1$. As $K$ grows larger, the low-to-high power
transition between $t^\star$ and $T^\star$ becomes more pronounced, supporting
the theoretical findings in Theorem~\ref{thm:lr_achieves_boundary} and
Theorem~\ref{thm:prior-works}.

\subsection{Quantile comparison}
\label{sec:quantile-comparison}

We analyze distributional properties of the sample sizes $\tau^\star$ and
$\tau^{\mathbf{\Pi}}$ for the SPRT and the mixture SPRT (the AV test associated to the adaptive test martingale from Definition~\ref{def:adaptive-test-mart}), respectively. Using the realizations of these stopping times, as described in the previous subsection
and in more detail in Appendix~\ref{app:simulation-details}, their quantiles are
readily estimated. Special attention is given to the sample size required by the
(mixture) SPRT to achieve $80\%$ of the maximum possible power, but qualitatively
similar results are obtained for other choices as well.

When sampling directly from the alternative hypothesis
in~(\ref{eq:generic-alternative}), the maximum power achievable by any test is
lower than one. This is because, under ${\cal H}_1$, there are zero anomalous
streams with probability $(1 - \varepsilon)^K$. Consequently, with nonzero
probability, the data are actually generated by the null hypothesis. If a test
for this problem had type-I error equal to $\alpha$ and rejected every time that
the number of anomalies is nonzero, its power, and therefore the maximum power
achievable by any level-$\alpha$ test for this problem, would be
$\gamma_{\max} = \alpha(1 - \varepsilon)^K + 1 - (1 - \varepsilon)^K$.

In each simulation scenario, the sample sizes required to achieve 80\% of the
maximum power $\gamma_{\max}$ are computed. We let $n^\star(0.8)$ be the lowest
value of $n$ such that $F^\star_\tau(n)\geq 0.8 \gamma_{\max}$. Then,
$n^\star(0.8)$ is the required maximum sample for the one-sided SPRT based on
$E^\star_t$ to attain $80\%$ of the maximum power. In other words, if the test
is stopped at $\tau^\star\wedge n^\star(0.8)$, then the resulting
(time-truncated) SPRT rejecting the null hypothesis when
$E^\star_{\tau^\star\wedge n^\star(0.8)} \geq A^\star$ has $80\%$ of the maximum
power $\gamma_{\max}$. Additionally,
$\mathbf{E}^\star[\tau^\star\wedge n^\star(0.8)]$ is the average sample of the
resulting test. Define, in the same vein, $n^{\mathbf{\Pi}}(0.8)$ to be the
smallest $n$ such that $F^{\mathbf{\Pi}}_\tau(n)\geq 0.8\gamma_{\max}$. In
Table~\ref{tab:average-st}, to simplify its presentation, we refer to
$n^\star(0.8)$ as ``Oracle'' $n_{\mathrm{AV}}(0.8)$ and to
$n^{\mathbf{\Pi}}(0.8)$ as ``Adaptive'' $n_{\mathrm{AV}}(0.8)$, and compare it
to the ``Oracle'' $n_{\mathrm{FS}}(0.8)$, which is the sample size required by
the fixed-sample likelihood ratio test to attain $80\%$ of the maximum power;
and to the ``Adaptive'' $n_{\mathrm{FS}}(0.8)$, the sample size required by the
HC test to achieve the same goal. We observe in Table~\ref{tab:average-st} that
although the maximum sample sizes $n_{\mathrm{AV}}(0.8)$ are larger, the average
sample sizes are smaller than those required by fixed-sample tests.
\begin{table}[ht]
  \label{tab:average-st}
  \caption{Comparison of sample sizes needed to attain $80\%$ power for
  fixed-sample tests and truncated (mixture) SPRT, and the average sample sizes for the
  (mixture) SPRT. ``Oracle'' $n_{\mathrm{FS}}(0.8)$ refers to the sample required by
  the Likelihood Ratio (fixed-sample) Test and $n_{\mathrm{AV}}(0.8)$ to the
  (well-specified) SPRT; ``Adaptive'' $n_{\mathrm{FS}}(0.8)$ to the HC test and
  $n_{\mathrm{AV}}(0.8)$ to the mixture martingale. While the worst case sample
  sizes needed to attain $80\%$ power are higher for AV tests than for
  fixed-sample tests, the average stopping time is lower. }
\centering
\begin{tabular}{lccccc}
  \hline
  \hline
  $K$ & $\beta^\star$ & Type & $\mathbf{E}^\star[\tau \wedge n_{\mathrm{AV}}(0.8)]$ & $n_{\mathrm{AV}}(0.8)$ & $n_{\mathrm{FS}}(0.8)$ \\
  \hline
  \hline
  \multirow{2}{*}{100}  & \multirow{2}{*}{0.55}
                      & Oracle   & 13.74 &  22 &  17 \\
      &  & Adaptive & 15.67 &  25 &  21 \\
  \hline
  \multirow{2}{*}{1000} & \multirow{2}{*}{0.55}
                      & Oracle   & 8.25 &  12 &  10 \\
      &  & Adaptive & 9.84 &  14 &  14 \\
  \hline
  \multirow{2}{*}{100}  & \multirow{2}{*}{0.85}
                      & Oracle   & 41.50 &  61 &  50 \\
      &   & Adaptive & 45.36 &  69 &  51 \\
  \hline
  \multirow{2}{*}{1000} & \multirow{2}{*}{0.85}
                      & Oracle  & 33.32 &  48 &  41 \\
      &  &  Adaptive & 36.79 &  55 &  41 \\
  \hline
\end{tabular}
\end{table}

\section{Conclusion and Discussion}
\label{sec:discussion}

This work addresses anytime-valid testing for anomalies in a large number of
data streams. Given the demand and flexibility of anytime-valid methods and the
nuances of the sparse anomaly detection problem, a fixed-$t$ and large-$K$ point
of view is taken. Theorem~\ref{thm:lr_achieves_boundary} in
Section~\ref{sec:results} shows that if the performance of test martingales is
measured using log-optimality criteria, the best possible test martingale---a
well-specified likelihood ratio---exhibits a threshold behavior as $K\to\infty$.
Under a specific parametrization with $\varepsilon^\star(K)$ tending to zero,
the sequences of $\delta^\star(K)$ against which any test martingale can have
nontrivial growth under the alternative are pinned down in
Theorem~\ref{thm:lr_achieves_boundary}. A similar threshold behavior is shown
for the stopping time of the associated AV test
(Theorem~\ref{thm:stop-on-time}). Section~\ref{sec:prior} shows a construction,
a mixture likelihood ratio summarized in
Definition~\ref{def:adaptive-test-mart}, that achieves the same threshold
behavior as the well-specified likelihood ratio without knowing the true
alternative parameters. The main results are the construction of
$E_t(\mathbf{\Pi})$ in Section~\ref{sec:prior}, and the characterization of both
its asymptotic growth properties in Theorem~\ref{thm:prior-works} and the
stopping time of the associated mixture SPRT in
Theorem~\ref{thm:prior-stops-on-time}. Section~\ref{sec:simulations} shows the
finite-$K$ power properties of these constructions through computer simulation.

In the following, we discuss a number of open problems and the relation between
this work and the existing literature.

\paragraph*{Computation} The parameter-adaptive statistic $E_{t}(\mathbf{\Pi})$
described in Section~\ref{sec:prior} can be computed, up to logarithmic factors
in $K$, in $O(K^{3 / 2})$ steps per time step. Furthermore, the analysis leading
to Theorem~\ref{thm:lr_achieves_boundary} shows that, even though there is
certain latitude in the choice of prior distribution $\mathbf{\Pi}$, this
computational complexity is not improvable for any discrete mixture of
likelihood ratios. The fact that, in the fixed-sample case, the higher criticism
statistic can be computed in $O(K)$ operations opens up an interesting question:
is there a test martingale that attains the log-optimality ``detection
boundary'' as in Theorem~\ref{thm:lr_achieves_boundary} and can be computed in
$O(K)$ steps? The answer to this question is far from trivial, and
Appendix~\ref{app:e-max-test} contains an example where an adaptation of ideas
from fixed-sample tests fails to achieve this goal.

\begin{figure}[ht!]
  \centering
  \includegraphics{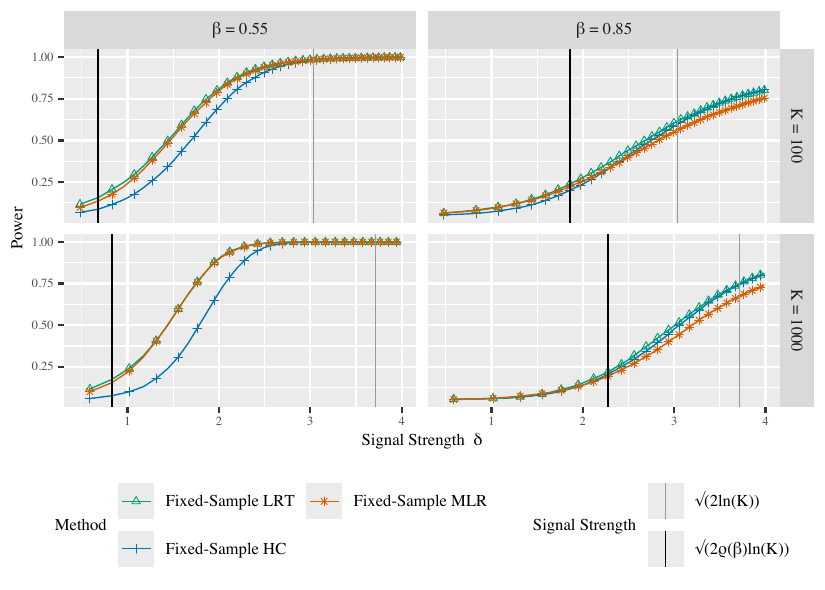}
  \caption{Power for fixed-sample methods. Comparison between tests based on the
    oracle Likelihood Ratio (LR), the Higher Criticism (HC) statistic and the
    Mixture Likelihood Ratio (MLR) from Section~\ref{sec:prior}. All tests are
    calibrated by Monte Carlo Simulation ($10^4$ runs) to attain 5\% type-I
    error. Power is estimated on $10^4$ runs. See more details in
    Section~\ref{sec:discussion}.}
  \label{fig:fspower}
\end{figure}

\paragraph*{E-statistics and fixed-sample tests} The adaptive AV tests that was
proposed also yields an adaptive fixed-sample test; namely,
$\indicator{E_t(\mathbf{\Pi}) \geq 1 / \alpha}$ for fixed $t$. This test is
based on the mixture likelihood ratio $E_t(\mathbf{\Pi})$, an example of an
e-statistic. An e-statistic is a nonnegative function of the data
$E_t = E_t(X_{i, s}: i\in [K], s\leq t)$ such that $\mathbf{E}_0[E_t] \leq 1$.
Using e-values, the realized values of e-statistics, has been proposed as a
basis for hypothesis
testing~\citep{ramdas_hypothesis_2025,grunwald_safe_2024,ramdas_game-theoretic_2023}
and fixed-sample log-optimality has been used as a way of quantifying power for
e-values. E-values play a pivotal role in multiple testing and allow for
post-hoc inference~\citep{ramdas_hypothesis_2025,grunwald_beyond_2024}. The
proof of Theorem~\ref{thm:stop-on-time} implies, under its assumptions, that the
fixed-sample test $\indicator{E_t(\mathbf{\Pi}) \geq 1 / \alpha}$ has
(conventional) power tending to one as $K\to\infty$ anytime that
$r^\star = t/T^\star$ is above the fixed-sample ``detection boundary''
$\rho(\beta^\star)$. This implies that the fixed-sample test
$\indicator{E_t(\mathbf{\Pi}) \geq 1 / \alpha}$ matches the power phase
transition of the (oracle) Likelihood Ratio Test and of the test based on the HC
statistic~(\ref{eq:hc-statisic})~\cite{donoho_higher_2004}.
Figure~\ref{fig:fspower} shows power plots for the $\alpha$-calibrated tests
(using Monte-Carlo simulation). In
some regimes, our mixture likelihood ratio test performs much better than the HC
test and nearly as well as the Likelihood Ratio Test; in other regimes, the
situation reverses (with a smaller margin). Although other adaptive tests are
known that achieve this phase transition, no known test dominates all the others
in a uniform sense~\citep{porter_beyond_2020,li_higher_2015}.

\paragraph*{The threshold $1/\alpha$ is approximate} Despite the usefulness of
having a universal threshold $1 / \alpha$, the type-I error of the AV test that
rejects when $E_t(\mathbf{\Pi})\geq 1/ \alpha$ is strictly smaller than
$\alpha$; that is, the resulting AV test is conservative. This conservativeness
has its origin in the fact that that the overshoot
$E_\tau(\mathbf{\Pi}) - 1/\alpha$ is positive at
$\tau = \inf\{t: E_t(\mathbf{\Pi}) \geq 1/\alpha \}$: if the overshoot is
exactly zero, the AV test has exact type-I error $\alpha$. In this work, the
following approach was taken: the value
$A^{\mathbf{\Pi}} = A^{\mathbf{\Pi}}(\alpha)$ such that the test
$\indicator{\sup_{t}E_{t}(\mathbf{\Pi}) > A^{\mathbf{\Pi}}}$ has level $\alpha$
was approximated through Monte Carlo simulation (see
Section~\ref{sec:simulations}). Approximating the value of the threshold $A$
analytically is challenging and existing analytical approximations (which
involve analyzing the distribution of the overshoot
$E_\tau(\mathbf{\Pi}) - 1/\alpha$) are, however, only applicable in simpler
scenarios than the sparse alternative that we consider
in~(\ref{eq:generic-alternative})~\citep{siegmund_sequential_1985,blier-wong_improved_2025}.
In order to reduce this conservativeness, additional methods have been developed
to ``boost'' $E_t(\mathbf{\Pi})$ and obtain a test martingale
$E _t^{*}(\mathbf{\Pi})$ that makes the same decision as $E _t(\mathbf{\Pi})$
and stops more often~\citep{fischer_improving_2025}. This ``boosting'' method
requires, however, computing explicitly
$\mathbf{E}_0[E_t(\mathbf{\Pi}) \ | \ {\cal F}_{t - 1}]$, which, although
possible numerically in our setting, would negate the computational advantages
of using the mixture distribution $\mathbf{\Pi}$ from Section~\ref{sec:prior}.

\paragraph*{Multiple testing} There has been a number of articles in the
literature unveiling a deep relation between e-statistics and multiple
testing~\citep{vovk_e-values_2021}, both in fixed-sample and in anytime-valid
scenarios~\citep{li_note_2025,ignatiadis_e-values_2024,xu_online_2024}. These
works, translated to our setup, focus on controlling control of False Discovery
Rate (FDR) and Non-Discovery Rate (NDR) when detecting the anomalies. Despite
their importance, these results are not applicable directly to our setting
because of the impossibility of controlling these quantities in the sample-size
regimes that are of our current interest (see Section~\ref{sec:state-things}).

\paragraph*{Bayesian Literature} The mixture distributions---priors---present in
the current work suggest a connection it to the Bayesian literature. In the Bayesian
literature a perhaps deceptively similar problem, the problem of finding
``needles in a haystack'' has been
treated~\citep{castillo_needles_2012,johnstone_needles_2004}. The problem is
similar in that there are $\varepsilon$-mixtures of Gaussian distributions, but
the focus is on the Bayesian identification of the anomalous observations while
our focus is on detection (see Section~\ref{sec:state-things}). Interestingly,
the best algorithms for the Bayesian problem also take $O(K^{3/2})$ operations
to compute and use discrete approximations to continuous
priors~\citep{erven_fast_2021}.

\begin{acks}[Acknowledgments]
  We thank Tyron Lardy and Wouter Koolen for fruitful discussions.
\end{acks}



\bibliographystyle{imsart-number}
\bibliography{criticism}

@misc{blier-wong_improved_2025,
	title = {Improved thresholds for e-values},
	url = {http://arxiv.org/abs/2408.11307},
	doi = {10.48550/arXiv.2408.11307},
	urldate = {2026-03-26},
	publisher = {arXiv},
	author = {Blier-Wong, Christopher and Wang, Ruodu},
	month = oct,
	year = {2025},
	note = {arXiv:2408.11307 [math]},
}

@article{grunwald_beyond_2024,
	title = {Beyond {Neyman}–{Pearson}: {E}-values enable hypothesis testing with a data-driven alpha},
	volume = {121},
	shorttitle = {Beyond {Neyman}–{Pearson}},
	url = {https://www.pnas.org/doi/10.1073/pnas.2302098121},
	doi = {10.1073/pnas.2302098121},
	number = {39},
	urldate = {2026-03-26},
	journal = {Proceedings of the National Academy of Sciences},
	publisher = {Proceedings of the National Academy of Sciences},
	author = {Grünwald, Peter D.},
	month = sep,
	year = {2024},
	pages = {e2302098121},
}

@article{ramdas_hypothesis_2025,
	title = {Hypothesis {Testing} with {E}-values},
	volume = {1},
	issn = {2978-4212},
	url = {https://doi.org/10.1561/3600000002},
	doi = {10.1561/3600000002},
	number = {1-2},
	urldate = {2026-03-26},
	journal = {Foundations and Trends in Statistics},
	author = {Ramdas, Aaditya and Wang, Ruodu},
	month = jul,
	year = {2025},
	pages = {1--390},
}

@article{li_higher_2015,
	title = {Higher criticism: $p$-values and criticism},
	volume = {43},
	issn = {0090-5364, 2168-8966},
	shorttitle = {Higher criticism},
	url = {https://projecteuclid.org/journals/annals-of-statistics/volume-43/issue-3/Higher-criticism-p-values-and-criticism/10.1214/15-AOS1312.full},
	doi = {10.1214/15-AOS1312},
	language = {en},
	number = {3},
	urldate = {2026-03-25},
	journal = {The Annals of Statistics},
	publisher = {Institute of Mathematical Statistics},
	author = {Li, Jian and Siegmund, David},
	month = jun,
	year = {2015},
	pages = {1323--1350},
}

@article{Zamba2021,
author = {Zamba, K. D.  and Tsiamyrtzis, P},
title = {Sequential detection framework for real-time biosurveillance based on Shiryaev-Roberts procedure with illustrations using COVID-19 incidence data},
journal = {Sequential Analysis},
volume = {40},
number = {2},
pages = {149--169},
year = {2021},
publisher = {Taylor \& Francis},
doi = {10.1080/07474946.2021.1912503},
URL = {https://doi.org/10.1080/07474946.2021.1912503}
}

@book{de_la_pena_self-normalized_2009,
	address = {Berlin, Heidelberg},
	series = {Probability and its {Applications}},
	title = {Self-{Normalized} {Processes}},
	copyright = {http://www.springer.com/tdm},
	isbn = {978-3-540-85635-1 978-3-540-85636-8},
	url = {http://link.springer.com/10.1007/978-3-540-85636-8},
	urldate = {2024-11-22},
	publisher = {Springer},
	author = {De La Peña, Victor H. and Lai, Tze Leung and Shao, Qi-Man},
	editor = {Gani, Joe and Heyde, C. C. and Jagers, P. and Kurtz, T. G.},
	year = {2009},
	doi = {10.1007/978-3-540-85636-8},
}

@book{hastie_elements_2009,
	address = {New York},
	edition = {2},
	series = {Springer {Series} in {Statistics}},
	title = {The {Elements} of {Statistical} {Learning}: {Data} {Mining}, {Inference}, and {Prediction}, {Second} {Edition}},
	isbn = {978-0-387-84857-0},
	shorttitle = {The {Elements} of {Statistical} {Learning}},
	url = {//www.springer.com/us/book/9780387848570},
	language = {en},
	urldate = {2018-07-19},
	publisher = {Springer-Verlag},
	author = {Hastie, Trevor and Tibshirani, Robert and Friedman, Jerome},
	year = {2009},
}

@article{dempster_maximum_1977,
	title = {Maximum {Likelihood} from {Incomplete} {Data} {Via} the {EM} {Algorithm}},
	volume = {39},
	issn = {0035-9246},
	url = {https://doi.org/10.1111/j.2517-6161.1977.tb01600.x},
	doi = {10.1111/j.2517-6161.1977.tb01600.x},
	number = {1},
	urldate = {2026-01-06},
	journal = {Journal of the Royal Statistical Society: Series B (Methodological)},
	author = {Dempster, A. P. and Laird, N. M. and Rubin, D. B.},
	month = sep,
	year = {1977},
	pages = {1--22},
}

@book{grunwald_minimum_2007,
	series = {Adaptive {Computation} and {Machine} {Learning} series},
	title = {The {Minimum} {Description} {Length} {Principle}},
	isbn = {978-0-262-52963-1},
	language = {Inglés},
	author = {Grünwald, Peter D.},
	month = mar,
	year = {2007},
}

@article{grunwald_minimum_2019,
	title = {Minimum description length revisited},
	volume = {11},
	issn = {2661-3352},
	url = {https://www.worldscientific.com/doi/10.1142/S2661335219300018},
	doi = {10.1142/S2661335219300018},
	number = {01},
	urldate = {2023-01-23},
	journal = {International Journal of Mathematics for Industry},
	author = {Grünwald, Peter D. and Roos, Teemu},
	month = dec,
	year = {2019},
	note = {Publisher: World Scientific Publishing Co.},
	pages = {1930001},
}

@article{Arias-Castro:2017,
author = {Arias-Castro, Ery and Chen, Shiyun},
title = {{Distribution-free multiple testing}},
volume = {11},
journal = {Electronic Journal of Statistics},
number = {1},
publisher = {Institute of Mathematical Statistics and Bernoulli Society},
pages = {1983 -- 2001},
year = {2017},
doi = {10.1214/17-EJS1277},
URL = {https://doi.org/10.1214/17-EJS1277}
}

@article{DS_Haupt_2011,
	author = {Haupt, Jarvis and Castro, Rui M and Nowak, Robert},
	journal = {IEEE Transactions on Information Theory},
	number = {9},
	pages = {6222--6235},
	publisher = {IEEE},
	title = {Distilled sensing: Adaptive sampling for sparse detection and estimation},
	volume = {57},
	year = {2011}}

@article{meinshausen_estimating_2006,
	title = {Estimating the proportion of false null hypotheses among a large number of independently tested hypotheses},
	volume = {34},
	issn = {0090-5364, 2168-8966},
	url = {https://projecteuclid.org/journals/annals-of-statistics/volume-34/issue-1/Estimating-the-proportion-of-false-null-hypotheses-among-a-large/10.1214/009053605000000741.full},
	doi = {10.1214/009053605000000741},
	number = {1},
	urldate = {2025-12-12},
	journal = {The Annals of Statistics},
	author = {Meinshausen, Nicolai and Rice, John},
	month = feb,
	year = {2006},
	note = {Publisher: Institute of Mathematical Statistics},
	pages = {373--393},
}

@article{cai_estimation_2007,
	title = {Estimation and confidence sets for sparse normal mixtures},
	volume = {35},
	issn = {0090-5364, 2168-8966},
	url = {https://projecteuclid.org/journals/annals-of-statistics/volume-35/issue-6/Estimation-and-confidence-sets-for-sparse-normal-mixtures/10.1214/009053607000000334.full},
	doi = {10.1214/009053607000000334},
	number = {6},
	urldate = {2025-12-15},
	journal = {The Annals of Statistics},
	author = {Cai, T. Tony and Jin, Jiashun and Low, Mark G.},
	month = dec,
	year = {2007},
	note = {Publisher: Institute of Mathematical Statistics},
	pages = {2421--2449},
}

@article{cai_optimal_2010,
	title = {Optimal rates of convergence for estimating the null density and proportion of nonnull effects in large-scale multiple testing},
	volume = {38},
	issn = {0090-5364, 2168-8966},
	url = {https://projecteuclid.org/journals/annals-of-statistics/volume-38/issue-1/Optimal-rates-of-convergence-for-estimating-the-null-density-and/10.1214/09-AOS696.full},
	doi = {10.1214/09-AOS696},
	number = {1},
	urldate = {2025-12-15},
	journal = {The Annals of Statistics},
	author = {Cai, T. Tony and Jin, Jiashun},
	month = feb,
	year = {2010},
	note = {Publisher: Institute of Mathematical Statistics},
	pages = {100--145},
}

@article{vovk_nonparametric_2024,
	title = {Nonparametric {E}-tests of {Symmetry}},
	volume = {2},
	issn = {2693-7166},
	url = {https://nejsds.nestat.org/journal/NEJSDS/article/69},
	doi = {10.51387/24-NEJSDS60},
	language = {en},
	number = {2},
	urldate = {2025-12-10},
	journal = {The New England Journal of Statistics in Data Science},
	author = {Vovk, Vladimir and Wang, Ruodu},
	month = feb,
	year = {2024},
	note = {Publisher: New England Statistical Society},
	pages = {261--270},
}

@misc{fischer_improving_2025,
	title = {Improving {Wald}'s (approximate) sequential probability ratio test by avoiding overshoot},
	url = {http://arxiv.org/abs/2410.16076},
	doi = {10.48550/arXiv.2410.16076},
	urldate = {2025-11-25},
	publisher = {arXiv},
	author = {Fischer, Lasse and Ramdas, Aaditya},
	month = jul,
	year = {2025},
	note = {arXiv:2410.16076 [stat]},
}

@article{koolen_log-optimal_2022,
	series = {Probability and {Statistics}: {Foundations} and {History}. {In} honor of {Glenn} {Shafer}},
	title = {Log-optimal anytime-valid {E}-values},
	volume = {141},
	issn = {0888-613X},
	url = {https://www.sciencedirect.com/science/article/pii/S0888613X21001523},
	doi = {10.1016/j.ijar.2021.09.010},
	urldate = {2025-11-20},
	journal = {International Journal of Approximate Reasoning},
	author = {Koolen, Wouter M. and Grünwald, Peter},
	month = feb,
	year = {2022},
	pages = {69--82},
}

@article{Bersimis:2007,
author = {Bersimis, S. and Psarakis, S. and Panaretos, J.},
title = {Multivariate statistical process control charts: an overview},
journal = {Quality and Reliability Engineering International},
volume = {23},
number = {5},
pages = {517-543},
doi = {https://doi.org/10.1002/qre.829},
year = {2007}
}

@article{Ingster1997,
	author = {Ingster, Yu. I.},
	journal = {Mathematical Methods of Statistics},
	number = {1},
	pages = {47--69},
	title = {{Some problems of hypothesis testing leading to infinitely divisible distributions}},
	url = {https://mathscinet.ams.org/mathscinet-getitem?mr=1456646},
	volume = {6},
	year = {1997}}

@inproceedings{xu_online_2024,
	author = {Xu, Ziyu and Ramdas, Aaditya},
	booktitle = {Proceedings of {The} 27th {International} {Conference} on {Artificial} {Intelligence} and {Statistics}},
	language = {en},
	month = apr,
	note = {ISSN: 2640-3498},
	pages = {3997--4005},
	publisher = {PMLR},
	title = {Online multiple testing with e-values},
	url = {https://proceedings.mlr.press/v238/xu24a.html},
	urldate = {2025-06-25},
	year = {2024}}

@article{ignatiadis_e-values_2024,
	author = {Ignatiadis, Nikolaos and Wang, Ruodu and Ramdas, Aaditya},
	doi = {10.1093/biomet/asad057},
	issn = {1464-3510},
	journal = {Biometrika},
	month = jun,
	number = {2},
	pages = {417--439},
	title = {E-values as unnormalized weights in multiple testing},
	url = {https://doi.org/10.1093/biomet/asad057},
	urldate = {2025-06-25},
	volume = {111},
	year = {2024}}

@article{li_note_2025,
	author = {Li, Guanxun and Zhang, Xianyang},
	doi = {10.1093/biomet/asae050},
	issn = {1464-3510},
	journal = {Biometrika},
	month = feb,
	number = {1},
	pages = {asae050},
	title = {A note on e-values and multiple testing},
	url = {https://doi.org/10.1093/biomet/asae050},
	urldate = {2025-06-25},
	volume = {112},
	year = {2025}}

@article{Zagar2022,
	author = {Zagar, Janja and Mihelic, Jurij},
	doi = {10.1038/s41597-022-01203-x},
	issn = {2052-4463},
	journal = {Scientific Data 2022 9:1},
	month = {mar},
	number = {1},
	pages = {1--11},
	pmid = {35322032},
	publisher = {Nature Publishing Group},
	title = {{Big data collection in pharmaceutical manufacturing and its use for product quality predictions}},
	url = {https://www.nature.com/articles/s41597-022-01203-x},
	volume = {9},
	year = {2022}}

@article{kulldorff2005stp,
	author = {Kulldorff, M. and Heffernan, R. and Hartman, J. and Assuncao, R. and Mostashari, F.},
	journal = {{PLOS} Medicine},
	number = {3},
	pages = {216},
	publisher = {PUBLIC LIBRARY OF SCIENCE},
	title = {A Space-Time Permutation Scan Statistic for Disease Outbreak Detection},
	volume = {2},
	year = {2005}}

@book{siegmund_sequential_1985,
	address = {New York, NY},
	author = {Siegmund, David},
	doi = {10.1007/978-1-4757-1862-1},
	isbn = {978-1-4419-3075-0 978-1-4757-1862-1},
	publisher = {Springer},
	series = {Springer {Series} in {Statistics}},
	title = {Sequential {Analysis}},
	url = {http://link.springer.com/10.1007/978-1-4757-1862-1},
	urldate = {2023-01-09},
	year = {1985}}

@book{wald_sequential_1947,
  author = {Wald, A.},
  publisher = {J. Wiley \& Sons, New York},
  title = {Sequential Analysis},
  year = 1947
}

@article{breiman_optimal_1961,
	author = {Breiman, L.},
	journal = {Proceedings of the Fourth Berkeley Symposium on Mathematical Statistics and Probability, Volume 1: Contributions to the Theory of Statistics},
	month = jan,
	note = {Publisher: University of California Press},
	pages = {65--79},
	title = {Optimal {Gambling} {Systems} for {Favorable} {Games}},
	url = {https://projecteuclid.org/ebooks/berkeley-symposium-on-mathematical-statistics-and-probability/Proceedings-of-the-Fourth-Berkeley-Symposium-on-Mathematical-Statistics-and/chapter/Optimal-Gambling-Systems-for-Favorable-Games/bsmsp/1200512159},
	urldate = {2022-08-26},
	volume = {4.1},
	year = {1961}}

@book{ville_etude_1939,
	address = {Paris},
	author = {Ville, Jean},
	number = {218},
	series = {Th{\`e}ses de l'entre-deux-guerres},
	title = {{\'E}tude {Critique} de la {Notion} de {Collectif}},
	url = {http://archive.numdam.org/item/THESE_1939__218__1_0/},
	year = {1939}}

@article{vovk_e-values_2021,
	author = {Vovk, Vladimir and Wang, Ruodu},
	doi = {10.1214/20-AOS2020},
	issn = {0090-5364, 2168-8966},
	journal = {The Annals of Statistics},
	month = jun,
	note = {Publisher: Institute of Mathematical Statistics},
	number = {3},
	pages = {1736--1754},
	shorttitle = {E-values},
	title = {E-values: {Calibration}, combination and applications},
	url = {https://projecteuclid.org/journals/annals-of-statistics/volume-49/issue-3/E-values-Calibration-combination-and-applications/10.1214/20-AOS2020.full},
	urldate = {2022-05-02},
	volume = {49},
	year = {2021}}

@article{ramdas_admissible_2020,
	author = {Ramdas, Aaditya and Ruf, Johannes and Larsson, Martin and Koolen, Wouter},
	journal = {arXiv:2009.03167 [math, stat]},
	month = sep,
	note = {arXiv: 2009.03167},
	title = {Admissible anytime-valid sequential inference must rely on nonnegative martingales},
	url = {http://arxiv.org/abs/2009.03167},
	urldate = {2020-11-27},
	year = {2020}}

@article{van_erven_renyi_2014,
	author = {van Erven, Tim and Harremo{\"e}s, Peter},
	doi = {10.1109/TIT.2014.2320500},
	issn = {0018-9448},
	journal = {IEEE Transactions on Information Theory},
	month = jul,
	number = {7},
	pages = {3797--3820},
	title = {R{\'e}nyi {Divergence} and {Kullback}-{Leibler} {Divergence}},
	volume = {60},
	year = {2014}}

@article{howard_uniform_2018,
	author = {Howard, Steven R. and Ramdas, Aaditya and McAuliffe, Jon and Sekhon, Jasjeet},
	journal = {arXiv:1810.08240 [math, stat]},
	month = oct,
	note = {arXiv: 1810.08240},
	title = {Uniform, nonparametric, non-asymptotic confidence sequences},
	url = {http://arxiv.org/abs/1810.08240},
	urldate = {2019-08-22},
	year = {2018}}

@article{wald_optimum_1948,
	author = {Wald, A. and Wolfowitz, J.},
	doi = {10.1214/aoms/1177730197},
	issn = {0003-4851, 2168-8990},
	journal = {The Annals of Mathematical Statistics},
	month = sep,
	note = {Publisher: Institute of Mathematical Statistics},
	number = {3},
	pages = {326--339},
	title = {Optimum {Character} of the {Sequential} {Probability} {Ratio} {Test}},
	url = {https://projecteuclid.org/journals/annals-of-mathematical-statistics/volume-19/issue-3/Optimum-Character-of-the-Sequential-Probability-Ratio-Test/10.1214/aoms/1177730197.full},
	urldate = {2023-02-07},
	volume = {19},
	year = {1948}}

@book{tartakovsky_sequential_2014,
	address = {New York},
	author = {Tartakovsky, Alexander and Nikiforov, Igor and Basseville, Michele},
	doi = {10.1201/b17279},
	isbn = {978-0-429-15234-4},
	month = aug,
	publisher = {Chapman and Hall/CRC},
	shorttitle = {Sequential {Analysis}},
	title = {Sequential {Analysis}: {Hypothesis} {Testing} and {Changepoint} {Detection}},
	year = {2014}}

@article{lai_power-one_1977,
	author = {Lai, Tze Leung},
	issn = {0090-5364},
	journal = {The Annals of Statistics},
	note = {Publisher: Institute of Mathematical Statistics},
	number = {5},
	pages = {866--880},
	title = {Power-{One} {Tests} {Based} on {Sample} {Sums}},
	url = {https://www.jstor.org/stable/2958514},
	urldate = {2023-02-09},
	volume = {5},
	year = {1977}}

@article{robbins_expected_1974,
	author = {Robbins, H. and Siegmund, D.},
	doi = {10.1214/aos/1176342704},
	issn = {0090-5364, 2168-8966},
	journal = {The Annals of Statistics},
	month = may,
	note = {Publisher: Institute of Mathematical Statistics},
	number = {3},
	pages = {415--436},
	title = {The {Expected} {Sample} {Size} of {Some} {Tests} of {Power} {One}},
	url = {https://projecteuclid.org/journals/annals-of-statistics/volume-2/issue-3/The-Expected-Sample-Size-of-Some-Tests-of-Power-One/10.1214/aos/1176342704.full},
	urldate = {2023-02-10},
	volume = {2},
	year = {1974}}

@article{waudby-smith_estimating_2023,
	author = {Waudby-Smith, Ian and Ramdas, Aaditya},
	doi = {10.1093/jrsssb/qkad009},
	issn = {1369-7412},
	journal = {Journal of the Royal Statistical Society Series B: Statistical Methodology},
	month = feb,
	pages = {qkad009},
	title = {Estimating means of bounded random variables by betting},
	url = {https://doi.org/10.1093/jrsssb/qkad009},
	urldate = {2023-02-21},
	year = {2023}}

@article{robbins_statistical_1970,
	author = {Robbins, Herbert},
	doi = {10.1214/aoms/1177696786},
	issn = {0003-4851, 2168-8990},
	journal = {The Annals of Mathematical Statistics},
	month = oct,
	note = {Publisher: Institute of Mathematical Statistics},
	number = {5},
	pages = {1397--1409},
	title = {Statistical {Methods} {Related} to the {Law} of the {Iterated} {Logarithm}},
	url = {https://projecteuclid.org/journals/annals-of-mathematical-statistics/volume-41/issue-5/Statistical-Methods-Related-to-the-Law-of-the-Iterated-Logarithm/10.1214/aoms/1177696786.full},
	urldate = {2023-02-20},
	volume = {41},
	year = {1970}}

@book{tsybakov_introduction_2009,
	address = {New York, NY},
	author = {Tsybakov, Alexandre B.},
	doi = {10.1007/b13794},
	isbn = {978-0-387-79051-0 978-0-387-79052-7},
	language = {en},
	publisher = {Springer},
	series = {Springer {Series} in {Statistics}},
	title = {Introduction to {Nonparametric} {Estimation}},
	url = {https://link.springer.com/10.1007/b13794},
	urldate = {2023-04-02},
	year = {2009}}

@article{kelly_jr_new_1956,
	author = {Kelly Jr., J. L.},
	doi = {10.1002/j.1538-7305.1956.tb03809.x},
	issn = {1538-7305},
	journal = {Bell System Technical Journal},
	language = {en},
	note = {\_eprint: https://onlinelibrary.wiley.com/doi/pdf/10.1002/j.1538-7305.1956.tb03809.x},
	number = {4},
	pages = {917--926},
	title = {A {New} {Interpretation} of {Information} {Rate}},
	url = {https://onlinelibrary.wiley.com/doi/abs/10.1002/j.1538-7305.1956.tb03809.x},
	urldate = {2023-03-29},
	volume = {35},
	year = {1956}}

@article{ingster_minimax_1994,
	author = {Ingster, Yu. I.},
	doi = {10.1007/BF01254275},
	issn = {1573-8795},
	journal = {Journal of Mathematical Sciences},
	language = {en},
	month = feb,
	number = {4},
	pages = {503--515},
	title = {Minimax detection of a signal in lp metrics},
	url = {https://doi.org/10.1007/BF01254275},
	urldate = {2024-04-08},
	volume = {68},
	year = {1994}}

@article{ramdas_game-theoretic_2023,
	author = {Ramdas, Aaditya and Gr{\"u}nwald, Peter and Vovk, Vladimir and Shafer, Glenn},
	doi = {10.1214/23-STS894},
	issn = {0883-4237, 2168-8745},
	journal = {Statistical Science},
	month = nov,
	note = {Publisher: Institute of Mathematical Statistics},
	number = {4},
	pages = {576--601},
	title = {Game-{Theoretic} {Statistics} and {Safe} {Anytime}-{Valid} {Inference}},
	url = {https://projecteuclid.org/journals/statistical-science/volume-38/issue-4/Game-Theoretic-Statistics-and-Safe-Anytime-Valid-Inference/10.1214/23-STS894.full},
	urldate = {2024-04-08},
	volume = {38},
	year = {2023}}

@article{castillo_needles_2012,
	author = {Castillo, Isma{\"e}l and Vaart, Aad van der},
	doi = {10.1214/12-AOS1029},
	issn = {0090-5364, 2168-8966},
	journal = {The Annals of Statistics},
	month = aug,
	note = {Publisher: Institute of Mathematical Statistics},
	number = {4},
	pages = {2069--2101},
	shorttitle = {Needles and {Straw} in a {Haystack}},
	title = {Needles and {Straw} in a {Haystack}: {Posterior} concentration for possibly sparse sequences},
	url = {https://projecteuclid.org/journals/annals-of-statistics/volume-40/issue-4/Needles-and-Straw-in-a-Haystack--Posterior-concentration-for/10.1214/12-AOS1029.full},
	urldate = {2024-03-26},
	volume = {40},
	year = {2012}}

@article{erven_fast_2021,
	author = {Erven, Tim van and Szab{\'o}, Botond},
	doi = {10.1214/20-BA1227},
	issn = {1936-0975, 1931-6690},
	journal = {Bayesian Analysis},
	month = sep,
	note = {Publisher: International Society for Bayesian Analysis},
	number = {3},
	pages = {933--960},
	title = {Fast {Exact} {Bayesian} {Inference} for {Sparse} {Signals} in the {Normal} {Sequence} {Model}},
	url = {https://projecteuclid.org/journals/bayesian-analysis/volume-16/issue-3/Fast-Exact-Bayesian-Inference-for-Sparse-Signals-in-the-Normal/10.1214/20-BA1227.full},
	urldate = {2024-03-26},
	volume = {16},
	year = {2021}}

@article{johnstone_needles_2004,
	author = {Johnstone, Iain M. and Silverman, Bernard W.},
	doi = {10.1214/009053604000000030},
	issn = {0090-5364, 2168-8966},
	journal = {The Annals of Statistics},
	month = aug,
	note = {Publisher: Institute of Mathematical Statistics},
	number = {4},
	pages = {1594--1649},
	shorttitle = {Needles and straw in haystacks},
	title = {Needles and straw in haystacks: {Empirical} {Bayes} estimates of possibly sparse sequences},
	url = {https://projecteuclid.org/journals/annals-of-statistics/volume-32/issue-4/Needles-and-straw-in-haystacks--Empirical-Bayes-estimates-of/10.1214/009053604000000030.full},
	urldate = {2024-03-26},
	volume = {32},
	year = {2004}}

@article{donoho_higher_2004,
	author = {Donoho, David and Jin, Jiashun},
	doi = {10.1214/009053604000000265},
	issn = {0090-5364, 2168-8966},
	journal = {The Annals of Statistics},
	month = jun,
	note = {Publisher: Institute of Mathematical Statistics},
	number = {3},
	pages = {962--994},
	title = {Higher criticism for detecting sparse heterogeneous mixtures},
	url = {https://projecteuclid.org/journals/annals-of-statistics/volume-32/issue-3/Higher-criticism-for-detecting-sparse-heterogeneous-mixtures/10.1214/009053604000000265.full},
	urldate = {2024-01-23},
	volume = {32},
	year = {2004}}

@article{arias-castro_detection_2019,
	author = {Arias-Castro, Ery and Ying, Andrew},
	doi = {10.1214/18-EJS1512},
	issn = {1935-7524, 1935-7524},
	journal = {Electronic Journal of Statistics},
	month = jan,
	note = {Publisher: Institute of Mathematical Statistics and Bernoulli Society},
	number = {1},
	pages = {208--230},
	shorttitle = {Detection of sparse mixtures},
	title = {Detection of sparse mixtures: higher criticism and scan statistic},
	url = {https://projecteuclid.org/journals/electronic-journal-of-statistics/volume-13/issue-1/Detection-of-sparse-mixtures-higher-criticism-and-scan-statistic/10.1214/18-EJS1512.full},
	urldate = {2024-01-23},
	volume = {13},
	year = {2019}}

@article{addario-berry_combinatorial_2010,
	author = {Addario-Berry, Louigi and Broutin, Nicolas and Devroye, Luc and Lugosi, G{\'a}bor},
	issn = {0090-5364},
	journal = {The Annals of Statistics},
	note = {Publisher: Institute of Mathematical Statistics},
	number = {5},
	pages = {3063--3092},
	title = {On {Combinatorial} {Testing} {Problems}},
	url = {https://www.jstor.org/stable/29765255},
	urldate = {2024-01-23},
	volume = {38},
	year = {2010}}

@article{stoepker_anomaly_2022,
author = {Stoepker, Ivo V. and Castro, Rui M. and Arias-Castro, Ery and van den Heuvel, Edwin},
title = {Anomaly Detection for a Large Number of Streams: A Permutation-Based Higher Criticism Approach},
journal = {Journal of the American Statistical Association},
volume = {119},
number = {545},
pages = {461--474},
year = {2024},
publisher = {Taylor \& Francis},
doi = {10.1080/01621459.2022.2126361},
eprint = {https://doi.org/10.1080/01621459.2022.2126361}
}

@article{grunwald_safe_2024,
	author = {Gr{\"u}nwald, Peter and de Heide, Rianne and Koolen, Wouter},
	doi = {10.1093/jrsssb/qkae011},
	issn = {1369-7412},
	journal = {Journal of the Royal Statistical Society Series B: Statistical Methodology},
	month = mar,
	pages = {qkae011},
	title = {Safe {Testing}},
	url = {https://doi.org/10.1093/jrsssb/qkae011},
	urldate = {2024-06-26},
	year = {2024}}

@article{stoepker_sparse_2023,
	author = {Ivo V. Stoepker and Rui M. Castro and Ery Arias-Castro},
	doi = {https://doi.org/10.1214/24-AOS2477},
	journal = {The Annals of Statistics},
	number = {2},
	pages = {676 -- 702},
	publisher = {Institute of Mathematical Statistics},
	title = {{Sparse anomaly detection across referentials: A rank-based higher criticism approach}},
	volume = {53},
	year = {2025}}

@article{porter_beyond_2020,
	author = {Porter, Thomas and Stewart, Michael},
	doi = {10.1214/19-AOS1885},
	issn = {0090-5364, 2168-8966},
	journal = {The Annals of Statistics},
	month = aug,
	note = {Publisher: Institute of Mathematical Statistics},
	number = {4},
	pages = {2230--2252},
	shorttitle = {Beyond {HC}},
	title = {Beyond {HC}: {More} sensitive tests for rare/weak alternatives},
	url = {https://projecteuclid.org/journals/annals-of-statistics/volume-48/issue-4/Beyond-HC-More-sensitive-tests-for-rareweak-alternatives/10.1214/19-AOS1885.full},
	urldate = {2025-03-14},
	volume = {48},
	year = {2020}}

@article{cai_optimal_2014,
	author = {Cai, Tony T. and Wu, Yihong},
	doi = {10.1109/TIT.2014.2304295},
	issn = {1557-9654},
	journal = {IEEE Transactions on Information Theory},
	month = apr,
	note = {Conference Name: IEEE Transactions on Information Theory},
	number = {4},
	pages = {2217--2232},
	title = {Optimal {Detection} of {Sparse} {Mixtures} {Against} a {Given} {Null} {Distribution}},
	url = {https://ieeexplore.ieee.org/document/6730948},
	urldate = {2025-03-14},
	volume = {60},
	year = {2014}}

@article{donoho_higher_2015,
	author = {Donoho, David and Jin, Jiashun},
	doi = {10.1214/14-STS506},
	issn = {0883-4237, 2168-8745},
	journal = {Statistical Science},
	month = feb,
	note = {Publisher: Institute of Mathematical Statistics},
	number = {1},
	pages = {1--25},
	title = {Higher {Criticism} for {Large}-{Scale} {Inference}, {Especially} for {Rare} and {Weak} {Effects}},
	url = {https://projecteuclid.org/journals/statistical-science/volume-30/issue-1/Higher-Criticism-for-Large-Scale-Inference-Especially-for-Rare-and/10.1214/14-STS506.full},
	urldate = {2025-03-14},
	volume = {30},
	year = {2015}}

@article{cai_optimal_2011,
	author = {Cai, Tony T. and Jessie Jeng, X. and Jin, Jiashun},
	doi = {10.1111/j.1467-9868.2011.00778.x},
	issn = {1369-7412},
	journal = {Journal of the Royal Statistical Society Series B: Statistical Methodology},
	month = nov,
	number = {5},
	pages = {629--662},
	title = {Optimal {Detection} of {Heterogeneous} and {Heteroscedastic} {Mixtures}},
	url = {https://doi.org/10.1111/j.1467-9868.2011.00778.x},
	urldate = {2025-03-14},
	volume = {73},
	year = {2011}}

@article{ingster_adaptive_2002,
	author = {Ingster, Y. I.},
	issn = {1066-5307},
	journal = {Mathematical Methods of Statistics},
	note = {OCLC: 196146366},
	number = {1},
	pages = {37--68},
	title = {Adaptive {Detection} of a {Signal} of {Growing} {Dimension}. {II}},
	volume = {11},
	year = {2002}}

@article{johari_always_2022,
	author = {Johari, Ramesh and Koomen, Pete and Pekelis, Leonid and Walsh, David},
	doi = {10.1287/opre.2021.2135},
	issn = {0030-364X},
	journal = {Operations Research},
	month = may,
	note = {Publisher: INFORMS},
	number = {3},
	pages = {1806--1821},
	shorttitle = {Always {Valid} {Inference}},
	title = {Always {Valid} {Inference}: {Continuous} {Monitoring} of {A}/{B} {Tests}},
	url = {https://pubsonline.informs.org/doi/abs/10.1287/opre.2021.2135},
	urldate = {2025-03-21},
	volume = {70},
	year = {2022}}

@article{wang_anytime-valid_2025,
author = {Wang, Hongjian and and Ramdas, Aaditya},
title = {Anytime-valid t-tests and confidence sequences for {Gaussian} means with unknown variance},
journal = {Sequential Analysis},
volume = {44},
number = {1},
pages = {56--110},
year = {2025},
publisher = {Taylor \& Francis},
doi = {10.1080/07474946.2024.2428245},
URL = {https://doi.org/10.1080/07474946.2024.2428245},
eprint = {https://doi.org/10.1080/07474946.2024.2428245}
}

@book{abramowitz_handbook_1965,
	author = {Abramowitz, Milton and Stegun, Irene A.},
	isbn = {978-0-486-61272-0},
	language = {en},
	month = jan,
	note = {Google-Books-ID: MtU8uP7XMvoC},
	publisher = {Courier Corporation},
	shorttitle = {Handbook of {Mathematical} {Functions}},
	title = {Handbook of {Mathematical} {Functions}: {With} {Formulas}, {Graphs}, and {Mathematical} {Tables}},
	year = {1965}}

\begin{appendix}

  \section{Lessons from an AV max test}
  \label{app:e-max-test}
  This section discusses the classic max test, an adaptive test
  for~(\ref{eq:generic-alternative}), and initial ideas that the authors had in using
  standard techniques to mimic its behavior with an AV test. These ideas provide a cautionary tale, suggesting that similar approaches on more sophisticated
  tests, such as a test based on the higher criticism statistic $\mathrm{HC}_t$
  from~(\ref{eq:hc-statisic}), might also not be fruitful.

  A classic adaptive fixed-sample test for~(\ref{eq:generic-alternative}) is the
  max test. For each $i\in [K]$ let $Z_{i, t} = (1/\sqrt{t})\sum_{s\leq t} X_{i, s}$ be the standardized stream sums. The fixed-sample max test rejects the null hypothesis when the maximum of these standardized stream sums, $\max_{i\in [K]} Z_{i, t}$, is large. The ``detection boundary'' for
  the max test is known and it differs from (it is worse than) $\rho$
  in~(\ref{eq:detection_boundary}), the ``detection boundary'' of the most
  powerful test~\citep{donoho_higher_2004}. Indeed, if judged under the same
  parametrization as in
  Section~\ref{sec:state-things}---$\varepsilon^\star = K^{-\beta^\star}$ and
  $\delta^\star = \sqrt{2r^\star\ln K}$---, the fixed-sample $\alpha$-calibrated
  max test is known to achieve a detection boundary
  $\rho_{\max}(\beta) = (1 - \sqrt{1 - \beta})^2$. In other words, if
  $r^\star > \rho_{\max}(\beta^\star)$, the max test has asymptotic power tending
  to one, while if $r^\star < \rho_{\max}(\beta^\star)$, its asymptotic power
  tends to zero~\citep{donoho_higher_2004}. This boundary coincides with the
  optimal detection boundary, $\rho$, for $\beta^\star\geq3/4$. Therefore, studying the max test provides a playground to
  try possible ideas for AV testing. We show an approach to mimic this test with an
  AV test.

  A common thread in the test-martingale literature is to approximate maxima of
  test martingales with mixtures of test martingales, which are interpreted soft
  maxima (a process also called pseudo-maximization) of test
  martingales~\citep[see][Chapter~11]{de_la_pena_self-normalized_2009}. In our
  context, a first idea might be to use this approach to mimic the fixed-sample
  max test using a soft maximum of test martingales. Indeed, given that
  $\max_{i\in [K]} Z_{i,t}$ can be approximated by a log-sum-exp (also called a
  soft-max), a first idea for an AV analogue of the max test is to use a mixture
  of the average of the single-stream likelihood ratios; that is, for a specific
  choice of $\mathbf{\Pi}$, to use $E^{\max}_t$ given by
  \begin{equation}\label{eq:e-max}
    E^{\max}_t
    =
    \int\bracks{
      (1 - \varepsilon)
      + \varepsilon\frac{1}{K} \sum_{i\in [K]}\rme^{\delta \sqrt{t}Z_{i, t} - \delta^2 t/ 2}
    }
    \rmd \mathbf{\Pi}(\varepsilon, \delta).
  \end{equation}
  Indeed, if $\mathbf{\Pi}$ puts mass on the maximizer of the last equation in
  $\delta$, then, as $t\to\infty$,~\eqref{eq:e-max} grows as a multiple of
  $\exp((\max_{i\in [K]}Z_{i, t})^2 / 2)$---this exponent is the
  maximum in $\delta$---, which, given the classic treatment of the problem,
  could point at a relationship with the fixed-sample max test.

  A question that arises is whether a test based on the statistic~\eqref{eq:e-max}
  achieves the same ``detection moment'' $t_{\max}^\star$ that would be expected
  from a max test AV upgrade,
  $t_{\max}^\star = T^\star \rho_{\max}(\beta^\star)$---in analogy to the result
  of Theorem~\ref{thm:lr_achieves_boundary}. Given that a similar approach could
  be used on $\mathrm{HC}_t$ from~(\ref{eq:hc-statisic}) to derive a test
  martingale, if a test based on~(\ref{eq:e-max}) achieves the same ``detection
  moment'' expected from an AV max test---for a properly chosen $\mathbf{\Pi}$---,
  there would be hope to derive, using this approach, an adaptive AV test that
  achieves the optimal ``detection moment'' $t^\star = \rho(\beta^\star)T^\star$
  (because the $\mathrm{HC}_t$ tests achieves the ``detection boundary'' $\rho$).
  This hope, however, turns out to be false: even if $\mathbf{\Pi}$ puts all of
  its mass on the ``correct'' values $(\varepsilon^\star, \delta^\star)$, the
  logarithmic growth of~\eqref{eq:e-max} does not achieve the ``detection
  boundary'' of the max test.
  \begin{lemma}\label{lem:e-max-log}
    Let
    $E_t^{\mathrm{max}} = (1 - \varepsilon^\star) + \varepsilon^\star(1 / K)
    \sum_{i\in [K]}\rme^{\delta^\star \sqrt{t}Z_{i, t} - \delta^{\star 2} t/ 2}$, and let
    $\varepsilon^\star = \varepsilon^\star(\beta^\star) = K^{-\beta^\star}$ and
    $\delta^\star = \delta^\star(T^\star) = \sqrt{2(1 / T^\star)\ln K}$. Then,
    anytime that $t / T^\star < \beta$, we have that
    $\mathbf{E}^\star[\ln E_t^{\mathrm{max}}] \to 0$ as $K \to \infty$.
  \end{lemma}
  \begin{proof}
    Use that $\ln x \leq x - 1$ for all $x > 0$ to deduce that
    $\mathbf{E}^\star[\ln E_t^{\mathrm{max}}]\leq
    \mathbf{E}^\star[E_t^{\mathrm{max}} - 1]$. A computation shows that
    $\mathbf{E}^\star[E_t^{\mathrm{max}}] = (1 - \varepsilon^\star) +
    \varepsilon^{\star 2}\rme^{\delta^{\star 2}t}$. Hence, setting
    $r^\star = t/T^\star$, we obtain that
    $\mathbf{E}^\star[\ln E_t^{\mathrm{max}}]\leq K^{2(r^\star - \beta^\star)} -
    \varepsilon^\star$, which implies the claim.
  \end{proof}
  Lemma~\ref{lem:e-max-log} shows that, in the expected logarithmic sense,
  the ``detection moment'' for the test martingale $E_t^{\mathrm{max}}$ is
  $t^\star = \beta^\star T^\star$, which is very far from the desired
  $t^\star = \rho_{\max}(\beta^\star) T^\star$ that would be expected from a
  max-test AV upgrade---in analogy with the results from
  Section~\ref{sec:results}. Similarly, the authors ran into similar issues when
  trying to apply this idea to the $\mathrm{HC}_t$ test (not shown).

  \section{Proofs}
  \label{app:proofs}

This appendix contains the proofs that
  were omitted from the main text. Section~\ref{sec:notation-proofs} introduces
  additional notation that was not used in the main body of the article but is
  necessary in the proofs; Section~\ref{sec:proof-lr_achieves_boundary} contains
  the proof of Theorem~\ref{thm:lr_achieves_boundary};
  Section~\ref{sec:proof-stop-on-time} contains the proof of
  Theorem~\ref{thm:stop-on-time}; Section~\ref{sec:proofs-prior} contains the
  proofs of both Proposition~\ref{prop:misspecification-cost} and
  Theorem~\ref{thm:prior-works}; Section~\ref{sec:proof-master-rates} contains
  the proof of Proposition~\ref{prop:master-rates} (which is the main result and
  technical tool in Section~\ref{sec:proofs-prior});
  Section~\ref{sec:technical-lemmas} contains technical lemmas related to
  Gaussian integrals; and Section~\ref{sec:proof-prior-stops} contains the proof
  of Theorem~\ref{thm:prior-stops-on-time}.

  \subsection{Notation}
  \label{sec:notation-proofs}

  Some additional notation is used in the proofs. To simplify the notation, the
  explicit dependency on $K$ is often dropped whenever it is clear from the
  context.

  \paragraph*{Probability distributions} For a single data stream $i\in [K]$, we
  use $\mathbf{P}_{\delta}$ to denote the distribution of i.i.d.
  $\mathrm{Normal}(\delta, 1)$ random variables. Crucially, under this notation,
  $\mathbf{P}^\star $ is the product distribution of $K$ copies of
  $(1 - \varepsilon^\star)\mathbf{P}_0 +
  \varepsilon^\star\mathbf{P}_{\delta^\star}$; the same relation exists between
  $\mathbf{P}_{\varepsilon,\delta}$ and
  $(1 - \varepsilon)\mathbf{P}_0 + \varepsilon\mathbf{P}_\delta$. We call
  $\mathbf{E}_\delta$ the expectation operator corresponding to
  $\mathbf{P}_\delta$.

  \paragraph*{Asymptotic notation}
  Additional asymptotic notation is used aside from the big-order notation used
  in the main text (see Section~\ref{sec:notation}). Let $f$ and $g$ be two
  real-valued functions, $f,g: \mathbf{N}\to \mathbf{R}$.
  We say that $f$ is in the small order of $g$ as $n\to\infty$, and write that
  $f = o(g)$ if $\lim_{n\to \infty} f(n)/g(n) = 0$.
  We say that $f \sim g$ as $n\to\infty$---$f$ is asymptotic to $g$---if
  $f(n) = g(n) (1 + o(1)) $ as $n\to\infty$.
  We say that $f\lesssim g$ as $n\to\infty $ if there is $\tilde{g}$ such that
  $f\leq \tilde{g}$ and $\tilde{g}\sim g$ as $n\to\infty$.
  We write that $f = \Omega(g)$ as $n\to\infty$ if $g = O(f)$ in the same limit.
  Write that $f = \Theta(g)$ as $n\to\infty$ if simultaneously $f = O(g)$ and
  $g = O(f)$ as $n\to\infty$.

  \paragraph*{Polylogarithmic factors} In this work, we are mostly concerned
  with functions that have exponential growth and, in some parts of the
  analysis, multiplicative polylogarithmic may be safely ignored (because they
  are of lower order). To this end, we define ``tilde'' versions of the previous
  asymptotic notation that hold up to polylogarithmic factors. We say that $f$
  is in the big order of $g$ up to polylogarithmic factors as $n\to\infty$ and
  write that $f = \tilde{O}(g)$ as $n\to\infty$ if there are $n_0\in\mathbf{N}$,
  and $C> 0$, and $a\in \mathbf{R}$ such that $|f(n)| \leq C |g(n)|(\ln n)^{a}$
  anytime that $n>n_0$. We define $\tilde{\Omega}$ and $\tilde{\Theta}$
  analogously.

  \paragraph*{Gaussian distributions} We write
  $\varphi(x) := (1/\sqrt{2\pi})\exp(-x^2 / 2)$ for the density of a standard
  Gaussian distribution. We write $\Phi(x) := \int_{-\infty}^x\varphi(w)\rmd w$
  for its cumulative distribution function, and we write
  $\overline{\Phi}(x) := 1 - \Phi(x)$ for its complementary cumulative
  distribution. At the risk of a slight abuse of notation, we also use
  $\Phi(a,b)$ to denote $\Phi(a,b) := \int_{a}^b\varphi(w)\rmd w$.

\subsection{Proof of Theorem~\ref{thm:lr_achieves_boundary}}
\label{sec:proof-lr_achieves_boundary}

The proof of Theorem~\ref{thm:lr_achieves_boundary} is a consequence of two
lemmas: Lemma~\ref{lem:oracle_upper_bound}, which is derived by bounding the
logarithmic growth of $E^\star_t$ from above, and
Lemma~\ref{lem:oracle_lower_bound}, which relies on lower bounds for the same
quantity.
\begin{lemma}\label{lem:oracle_upper_bound}
  Under the setup of Theorem~\ref{thm:lr_achieves_boundary}, if
  $t / T^\star < \rho(\beta^\star)$, then $\mathbf{E}^\star[\ln E^\star_t]\to 0$ as
  $K\to \infty$.
\end{lemma}
\begin{lemma}\label{lem:oracle_lower_bound}
  Under the setup of Theorem~\ref{thm:lr_achieves_boundary}, if
  $t / T^\star > \rho(\beta^\star)$, then $\mathbf{E}^\star[\ln E^\star_t]\to \infty$ as
  $K\to \infty$.
\end{lemma}
The two lemmas together yield Theorem~\ref{thm:lr_achieves_boundary}. The proof
of Lemma~\ref{lem:oracle_upper_bound} is technical, and is presented in
Section~\ref{sec:oracle_upper_bound}. Lemma~\ref{lem:oracle_lower_bound} is
proven in Section~\ref{lem:oracle_lower_bound} and employs known relationships
between the Kullback-Leibler divergence and the total variation distance, and
the relation of the latter to the risk of fixed-sample tests.

\subsubsection{Proof of Lemma~\ref{lem:oracle_upper_bound}}
\label{sec:oracle_upper_bound}

By sufficiency, $E^\star_t$ is a function of the stream averages only. Define
the standardized average of the stream $i\in [K]$ as
$Z_{i,t} := (1 / \sqrt{t})\sum_{s\leq t} X_{i, s}$. Furthermore let
$r^\star := t / T^\star$ and
$\Delta^\star := \sqrt{t}\delta^\star = \sqrt{2r^\star\ln K}$, and define
$S_{i,t} := (1 - \varepsilon^\star) +
\varepsilon^\star p_{\Delta^\star}(Z_{i,t}) / p_0(Z_{i,t})$. Notice that
when $X_{i, 1}, \dots, X_{i, t}\iidsim \mathrm{Normal}(\delta^\star, 1)$ we have
$Z_{i,t}\sim \mathrm{Normal}(\Delta^\star, 1)$. With these rewritings
$$\mathbf{E}^\star\sqbrack{ \ln E^{\star}_t } = \sum_{i\in [K]}\mathbf{E}^\star\sqbrack{ \ln S_{i,t} } = K\mathbf{E}^\star\sqbrack{ \ln S_{1,t}}.$$
Recall that
$\mathbf{E}^\star\sqbrack{ \ln S_{1,t} } = \varepsilon^\star
\mathbf{E}_{\delta^\star}\sqbrack{ \ln S_{1,t} } + (1 - \varepsilon^\star)
\mathbf{E}_{0}\sqbrack{ \ln S_{1,t} }$. By Jensen's inequality,
$ \mathbf{E}_{0}\sqbrack{ \ln S_{1,t}} \leq 0$ and consequently
  \begin{equation}\label{eq:basic_bound_KL}
    \mathbf{E}^\star\sqbrack{
      \ln
      E^\star_t
    }
    \leq
    K\varepsilon^\star
    \mathbf{E}_{\delta^\star}\sqbrack{
      \ln
      S_{1,t}
    }.
  \end{equation}
  We now bound $\mathbf{E}_{\delta^\star}\sqbrack{ \ln S_{1,t} }$ separately in
  several regions of the $\beta^\star$-$r^\star$~plane. First, assume that
  $r^\star < \beta^\star - 1/2$. Let $Z$ be a standard normal random variable.
  Using the inequality $\ln x \leq x - 1$, we obtain
  \begin{equation}\label{eq:basic_trick}
    \mathbf{E}^\star\sqbrack{
      \ln
      E^\star_t
    }
    \leq
    K\varepsilon^{\star 2}
    \mathbf{E}_{\delta^\star}
    \sqbrack{
      \frac{p_{\Delta^\star}(Z_{1,t})}{p_0(Z_{1,t})}
      - 1
    }
    =
    K\varepsilon^{\star 2}
    \mathbf{E}_{0}
    \sqbrack{
      \rme^{\Delta^\star(Z + \Delta^\star) - \Delta^{\star 2} / 2}
      - 1
    },
  \end{equation}
  where, in the last step we used the location nature of the normal
  distribution. Recall that
  $\mathbf{E}_{0}[\rme^{\Delta^{\star} Z}] = \rme^{\Delta^{\star 2} / 2}$ and
  use this to obtain that
  $ \mathbf{E}^\star\sqbrack{ \ln E^\star_t } \leq K \varepsilon^{\star 2}
  \paren{ \rme^{\Delta^{\star 2}} - 1 }$. In the $\beta^\star$-$r^\star$
  parametrization this translates to
  $ \mathbf{E}^\star\sqbrack{ \ln E^\star_t } \leq K^{1 - 2\beta^\star} \paren{
    K^{2r^\star} - 1 }$, which tends to zero as $K\to \infty$ because
  $1-2\beta^\star+2r^\star < 0$ when $r^\star < \beta^\star - 1/2$.

  Now, we move to a different region in the $\beta^\star$-$r^\star$~plane. The
  proof uses a truncation argument, reminiscent of the truncated second-moment
  method used to show lower bounds in hypothesis
  testing~\citep{ingster1997,stoepker_anomaly_2022}. Let
  $3/4 < \beta^\star \leq 1$, and assume that
  $(1 - \sqrt{1 - \beta^\star})^2 > r^\star \geq \beta^\star - 1/2$. In
  particular, $r^\star > 1/4$. Consider the event
  $\Psi = \{Z_{1,t} \leq \sqrt{2\ln K}\}$. Using
  (\ref{eq:basic_bound_KL}), we split the right-hand side using $\Psi$; that
  is,
  \begin{equation}\label{eq:target_two_part}
    \mathbf{E}^\star[\ln E^\star_t]
    \leq
    K\varepsilon^\star\mathbf{E}_{\delta^\star}[\indicator{\Psi}\ln S_{1,t}]
    +
    K\varepsilon^\star\mathbf{E}_{\delta^\star}[\indicator{\Psi^c}\ln S_{1,t}].
  \end{equation}
  The rest of the proof is dedicated to bound both terms of
  (\ref{eq:target_two_part}). For the first term, in the event $\Psi$, we use
  again the inequality $\ln x \leq x - 1$ as in (\ref{eq:basic_trick}) to deduce
  that
  $\mathbf{E}_{\delta^\star}[\indicator{\Psi}\ln S_{1,t}] \leq
  \varepsilon^\star\mathbf{E}_{\delta^\star}\sqbrack{ \indicator{\Psi} \paren{
      p_{\Delta^\star}(Z_{1,t}) / p_0(Z_{1,t}) - 1} }$. Use a change of
  measure and the location nature of the normal distribution to deduce that
  \begin{equation*}
    \mathbf{E}_{\delta^\star}\sqbrack{ \indicator{\Psi}
      \frac{p_{\Delta^\star}(Z_{1,t})}{p_0(Z_{1,t})}}
    =
    \rme^{\delta^{\star 2}}
    \mathbf{P}_{0}\bracks{Z + 2\Delta^\star \leq \sqrt{2\ln K}}
    \leq
    K^{2r^\star-(1 - \sqrt{4r^\star})^2}.
  \end{equation*}
  We can conclude that the first term of (\ref{eq:target_two_part}) is bounded
  by $K^{1-2\beta^\star + 2r^\star-(1 - \sqrt{4r^\star})^2}$, which tends to
  zero because both the inequalities
  $1-2\beta^\star + 2r^\star-(1 - \sqrt{4r^\star})^2< 0$ (in the exponent of
  $K$) and $r^\star < (1 - \sqrt{1 - \beta^\star})^2 $ (our assumption) define
  the same region in the $\beta^\star$-$r^\star$~plane. Now, we come to bound
  the second term in the right hand side of (\ref{eq:target_two_part}). Here, we
  use a tangent bound for $x\mapsto \ln x$, the fact that
  $\ln x \leq \ln a + x / a - 1$ for any $a > 0$, at $a = K^{2r^\star}$, that
  is, $\ln x \leq 2r^\star\ln K + x / K^{2r^\star} - 1$. Indeed,
  \begin{multline*}
    \mathbf{E}_{\delta^\star}\sqbrack{
    \indicator{\Psi^c}
    \ln S_{1,t}
    }
    \leq
      \mathbf{E}_{\delta^\star}\sqbrack{
      \indicator{\Psi^c}
      \paren{
      2r^\star\ln K
      +
      \frac{
      1
      +
      \varepsilon^\star
      \paren{
      \frac{p_{\Delta^\star}(Z_{1,t})}{p_0(Z_{1,t})} - 1
      }
      }{
      K^{2r^\star}
      }
      -1
      }
      }\\
    =
      \mathbf{P}_{\delta^\star}\{
      \Psi^c
      \}
      \paren{
      2r^\star\ln K
      +
      \frac{1}{K^{2r^\star}}
      -
      1
      }
      +
      \frac{\varepsilon^\star}{K^{2r^\star}}
      \mathbf{E}_{\delta^\star}\sqbrack{
      \indicator{\Psi^c}
      \paren{
      \frac{p_{\Delta^\star}(Z_{1,t})}{p_0(Z_{1,t})} - 1
      }
      }.
  \end{multline*}
  Use that
  $K\varepsilon^\star\mathbf{P}_{\delta^\star}\{ \Psi^c \} =
  K\varepsilon^\star\mathbf{P}_{0}\{ Z + \Delta^\star \geq \sqrt{2\ln K} \}\leq
  K^{1-\beta^\star}K^{-(1 - \sqrt{r^\star})^2}$ and the not-so-subtle bound
  $\mathbf{E}_{\delta^\star}\sqbrack{ \indicator{\Psi^c} \paren{
     p_{\Delta^\star}(Z_{1,t}) / p_0(Z_{1,t}) - 1 } } \leq
  \mathbf{E}_{\delta^\star}\sqbrack{
    p_{\Delta^\star}(Z_{1,t}) / p_0(Z_{1,t}) } = \rme^{\Delta^{\star^2}} =
  K^{2r^\star}$ to obtain that the second term in (\ref{eq:target_two_part})
  satisfies
  \begin{equation*}
    K\varepsilon^\star
    \mathbf{E}_{\delta^\star}\sqbrack{
      \indicator{\Psi^c}
      \ln S_{1,t}
    }
    \leq
    K^{1-\beta^\star-(1 - \sqrt{r^\star})^2}
    \paren{
      2r^\star\ln K
      +
      \frac{1}{K^{2r^\star}}
      -
      1
    }
    +
    K^{1- 2\beta^\star},
  \end{equation*}
  which implies that it tends to zero as $K\to\infty$ for our current selection
  of $\beta^\star$ and $r^\star$. Conclude that (\ref{eq:target_two_part}) tends
  to zero as $K\to\infty$ anytime that $r^\star < \rho(\beta^\star)$. This is
  what was to be shown.

\subsubsection{Proof of Lemma~\ref{lem:oracle_lower_bound}}
\label{sec:oracle_lower_bound_proof2}

Lemma~\ref{lem:oracle_lower_bound} follows from a number of observations and a
result from the literature~\citep[][Theorem~1.2]{donoho_higher_2004}. To show
how, we first recall a number of standard definitions and known facts about
fixed-sample statistical tests, about the total variation distance, and about
its relationship to the Kullback-Leibler divergence.

Recall that, for two probability distributions
$\mathbf{P}_0, \mathbf{P}_1$ defined on the same probability space, the total
variation distance $d_{\mathrm{TV}}(\mathbf{P}_0, \mathbf{P}_1)$ is defined as
$d_{\mathrm{TV}}(\mathbf{P}_0, \mathbf{P}_1) := \sup_{A}|\mathbf{P}_0\{A\} -
\mathbf{P}_1\{A\}|$, where the supremum is over the class of measurable sets. A
test $\psi$ between $\mathbf{P}_0$ and $\mathbf{P}_1$ is a random variable
taking values in $\{0, 1\}$ and its risk is
$R(\psi) := \mathbf{P}_0\{\psi = 1\} + \mathbf{P}_1\{\psi = 0\}$---the sum of
the type-I and type-II error probabilities. The risk $R(\psi)$ of any test $\psi$  satisfies
\citep[Theorem~2.2]{tsybakov_introduction_2009}
\begin{equation}\label{eq:tv-by-risk}
d_{\mathrm{TV}}(\mathbf{P}_0, \mathbf{P}_1) \geq 1 - 2 R(\psi) \ .
\end{equation}
The Kullback-Leibler divergence is defined as
\begin{equation*}
  \KL(\mathbf{P}_1, \mathbf{P}_0) := \int \ln\left(\frac{\rm d
      \mathbf{P}_1}{\rm d \mathbf{P}_0}\right) \rm d \mathbf{P}_1 ,
\end{equation*}
provided that $\mathbf{P}_1$ is absolutely continuous with respect to
$\mathbf{P}_0$, and $\KL(\mathbf{P}_1, \mathbf{P}_0) = \infty$ otherwise. It has
been shown~\citep[][Equation~2.25]{tsybakov_introduction_2009} that
\begin{equation}\label{eq:kl-by-tv}
  d_{\mathrm{TV}}(\mathbf{P}_0, \mathbf{P}_1)
  \leq
  1 - \frac{1}{2}\rme^{-\KL(\mathbf{P}_1, \mathbf{P}_0)} \ .
\end{equation}
Using \eqref{eq:tv-by-risk} and \eqref{eq:kl-by-tv}, we find that
\begin{equation}\label{eq:kl-by-risk}
  \KL(\mathbf{P}_1, \mathbf{P}_0)
  \geq
  - \ln\paren{2(1 - d_{\mathrm{TV}}(\mathbf{P}_0, \mathbf{P}_1))}
  \geq
  - \ln(4 R(\psi)) .
\end{equation}
Now consider a sequence $(\mathbf{P}_{1, K})_{K\geq 1}$ of probability distributions
(of alternative hypotheses), and consider a sequence $(\psi_{K})_{K\geq 1}$ of
tests, where each $\psi_K$ is a test of $\mathbf{P}_0$ against $\mathbf{P}_{1, K}$.
If the sequence of tests is such that $R(\psi_K)\to 0$ as $K\to \infty$, then,
due to \eqref{eq:kl-by-risk}, the Kullback-Leibler divergence tends to infinity
in the same limit.

This last observation is the crucial observation in proving
Lemma~\ref{lem:oracle_lower_bound}.

\begin{proof}[Proof of Lemma~\ref{lem:oracle_lower_bound}]
  With the discussion above in mind, notice that at each value of $K$, in our
  parametrization, an alternative distribution $\mathbf{P}_{1, K}$ is defined on
  the data. Notice also that
  \begin{equation*}
    \mathbf{E}^\star[\ln E^\star_t]  = \KL(\mathbf{P}_0, \mathbf{P}_{1,K}).
  \end{equation*}
  Additionally, there exists a sequence of tests $(\psi_K)_{K \geq 1}$ based on
  the higher-criticism statistic, where each $\psi_K$ is a test of
  $\mathbf{P}_0$ against $\mathbf{P}_{1, K}$, whose risk $R(\psi_K)$ tends to
  zero as $K\to\infty$~\citep[][Theorem 1.2]{donoho_higher_2004}. These
  observations and the discussion above---by \eqref{eq:kl-by-risk} in
  particular---show the claim.
\end{proof}

\subsection{Proof of Theorem~\ref{thm:stop-on-time}}
\label{sec:proof-stop-on-time}

\begin{proof}[Proof of Theorem~\ref{thm:stop-on-time}]
First, we show that $\mathbf{P}^\star\{\tau^\star \leq t\} \to 0$ whenever
$t < t^\star = T^\star\rho(\beta^\star)$. Consider the sequential test that
rejects null hypothesis $\mathcal{H}_0$ in (\ref{eq:generic-alternative}) with
stopping rule $\min\{\tau^\star, t\}$ and rejects the null hypothesis if
$\tau^\star \leq t$. Notice that $\mathbf{P}^\star\{\tau^\star \leq t\}$ is the
power of said sequential test and that the sample size of that test is smaller
than $t$. In particular, this sequential test is less powerful than the most
powerful test at sample size $t$. Since the most powerful test has power tending
to zero as $K\to\infty$ when the sample size
$t < t^\star$~\citep{ingster_minimax_1994,donoho_higher_2004}---because
$t < t^\star$ is equivalent to $\sqrt{t}\delta^\star < \rho(\beta^\star)$---,
then $\mathbf{P}^\star\{\tau^\star \leq t\} \to 0$ necessarily in the same
limit. This shows the first claim.

To prove the second claim, pick $t$ such that
$t > t^\star = T^\star\rho(\beta^\star)$. We show that
$\mathbf{P}^\star\{\tau^\star > t\} \to 0$ as $K\to \infty$. To this end, it is
enough to show that $\mathbf{P}^\star\{E_t^\star < 1 / \alpha\}\to 0$ in the
same limit because
$\mathbf{P}^\star\{\tau^\star > t\} = \mathbf{P}^\star\{\max_{s\leq t} E_s^\star
< 1 / \alpha\}\leq \mathbf{P}^\star\{E_t^\star < 1 / \alpha\}$. To bound this
probability, we use Cram\'{e}r--Chernoff's method on $\ln E_t^\star$, which becomes
a sum of $K$ independent random variables. With this in mind define, similarly
as in Section~\ref{sec:oracle_upper_bound}, the likelihood ratio for each stream
\begin{equation*}
  S_{i, t} := 1 - \varepsilon^\star + \varepsilon^\star \frac{p_{\delta^\star}(X_{i, 1},
    \dots, X_{i, t})}{p_0(X_{i, 1}, \dots, X_{i, t})},
\end{equation*}
  so that
  $\ln E^\star_t = \sum_{i\in [K]}\ln S_{i, t}$. We write
  \begin{equation*}
    \mathbf{P}^\star\left\{
      E^\star_t
      < 1 / \alpha
    \right\}
    =
    \mathbf{P}^\star\left\{
      -\sum_{i\in [K]} \ln S_{i,t}
      >
      -\ln(1 / \alpha)
    \right\}.
  \end{equation*}
  Let
  \begin{equation*}
    \lambda \mapsto \xi(\lambda) := \ln \mathbf{E}^\star[\rme^{- \lambda (\ln S_{i,t})}]
    = \ln \mathbf{E}^\star[S_{i,t}^{-\lambda}]
  \end{equation*}
  be the cumulant generating function of the random variable $- \ln S_{i,t}$.
  Observe both that $\xi(0) = \xi(1) = 0$ and that
  $\lambda \mapsto \xi(\lambda)$ is strictly convex. These two observations
  imply that $\xi(\lambda) < 0$ for each $\lambda \in (0, 1)$. Fix
  $\lambda\in(0, 1)$ now and for the rest of the proof, and use the
  Cram\'{e}r--Chernoff bound to obtain that
  \begin{align}
    \mathbf{P}^\star\left\{
      -\sum_{i\in [K]} \ln S_{i, t}
      >
      -\ln(1 / \alpha)
    \right\}\nonumber
    &=
    \mathbf{P}^\star\left\{
      \exp\paren{-\lambda\sum_{i\in [K]} \ln S_{i, t}}
      >
      \exp\paren{-\lambda\ln(1 / \alpha)}
    \right\}\\
    &\leq\label{eq:tail-prob}
    \rme^{K \xi(\lambda) + \lambda\ln(1/\alpha)}.
  \end{align}
  This inequality suggests that since $\xi(\lambda)< 0$, then
  $\mathbf{P}^\star\{E_t^\star < 1 / \alpha\}\to 0$ exponentially fast as
  $K\to\infty$. However, this is may not necessarily the case because $\xi$ also
  depends on $K$; we need to rule out the possibility that
  $K\xi(\lambda)\not\to -\infty$ as $K\to \infty$. Consequently, it is enough to
  show, for $K\to\infty$, that $\xi(\lambda)\not\to 0$. To achieve this goal, we
  use that
  \begin{equation}
    \label{eq:xi-tv}
    \xi(\lambda) \leq -(1 - \lambda)(\lambda / 2)d_{\mathrm{TV}}^2(p_{0, t},
    p_{\delta^\star,t}),
  \end{equation}
  where $d_{\mathrm{TV}}(p_{\delta^\star,t}, p_{0, t})$ is the total variation
  distance and $p_{\delta,t} = p_{\delta,t}(x_1, \dots, x_t)$ is the joint
  density of $t$ i.i.d. normal random variables with mean $\delta$ and unit
  variance~\citep[][Theorem~31]{van_erven_renyi_2014}. Equation (\ref{eq:xi-tv})
  is a special case of Theorem~31 of \citep{van_erven_renyi_2014} written in our
  notation, but it is originally stated in terms of the Renyi divergence
  $\mathrm{D}_\lambda(p_{0, t}, p_{\delta^\star,t}) := \frac{1}{\lambda -
    1}\mathbf{E}_{\delta^\star}\left[\left\{p_{0, t}(X_{1,1},\ldots,X_{1,t}) /
    p_{\delta^\star, t}(X_{1,1},\ldots,X_{1,t})\right\}^\lambda\right]$, which, in our notation,
  is
  $\xi(\lambda) = -(1 - \lambda)\mathrm{D}_\lambda(p_{0, t},
  p_{\delta^\star,t})$. In terms of the Renyi divergence, (\ref{eq:xi-tv})
  asserts that
  $\mathrm{D}_\lambda(p_{0, t}, p^\star_t) \geq (\lambda /
  2)d_{\mathrm{TV}}^2(p_{0, t}, p_{\delta^\star,t})$, and it is generalization
  of Pinsker's inequality. After this translation, use (\ref{eq:tail-prob}) and
  (\ref{eq:xi-tv}) to conclude that
  \begin{equation*}
    \mathbf{P}^\star\left\{
      E_t^\star < 1 / \alpha
    \right\}
    \leq
    \exp\paren{- K (1 - \lambda)\frac{\lambda ^2}{2}d_{\mathrm{TV}}^2(p_{0, t},
      p^\star_t) +  \lambda\ln(1 / \alpha)}.
  \end{equation*}
  In the case we are treating, where $t/T^\star > \rho(\beta)$ and we are above
  the fixed-sample ``detection boundary'',
  $d_{\mathrm{TV}}^2(p_{0, t}, p_{\delta^\star,t}) \to 1$ as $K\to\infty$.
  Consequently, by the previous display,
  $\mathbf{P}^\star\left\{ E_t^\star < 1 / \alpha \right\}\to 0$ in the same
  limit. By the discussion at the beginning of the proof, conclude that
  $\mathbf{P}^\star\{\tau^\star \leq t\} \to 1$ as $K\to\infty$, as it was to be
  shown.
\end{proof}

\subsection{Proofs of Proposition~\ref{prop:misspecification-cost} and
  Theorem~\ref{thm:prior-works}}
\label{sec:proofs-prior}

Both Proposition~\ref{prop:misspecification-cost} and Theorem~\ref{thm:prior-works}
follow from the next technical proposition, which characterizes the rate at
which $\mathbf{E}^\star[\ln E_t(\varepsilon, \delta)]$ grows to infinity as
$K\to\infty$ above the ``detection boundary''. The proof of this proposition is
technical and it is contained in Section~\ref{sec:proof-master-rates}.
\begin{proposition}\label{prop:master-rates}
  Assume that $1\geq \beta > 1/2$, let $(\beta_K)_{K\geq 1}\subseteq (1/2, 1]$, let
  $(T_K)_{K\geq 1}$ be a nonnegative sequence, and define both
  $r^\star = t / T^\star$ and $r_K = t / T_K$. Let
  $d_T = d_T(K) := 1 - \sqrt{T_K / T^\star}$ and
  $d_\beta = d_\beta(K) := \beta_K - \beta^\star$ and assume that, as
  $K\to\infty$, both $d_T K^{1 - \beta} = O(1)$ and $|d_\beta| \ln K \to 0$. Then, as $K\to\infty$,
  \begin{equation} \label{eq:master-rates}
  \mathbf{E}^\star[\ln E_t(\varepsilon,\delta)]
  =
  \begin{cases}
    \Omega(K^{1 + 2r^\star - 2\beta^\star})  &\text{ if }
                                               1/2 < \beta^\star < 3/4,
    \\
                                             & \text{ and } \beta^\star - 1/2
                                               < r^\star \leq \beta^\star / 3;
    \\
    \tilde{\Omega}(
    K^{1  - \beta^\star -
    \frac{1}{4}\paren{\frac{r^\star
    - \beta^\star}{\sqrt{r^\star}}}^
    2}
    )
                                             &\text{ if }
                                               1/2 < \beta^\star \leq 1,
    \\
                                             &\text{ and }
                                               \max\{\beta^\star/3,
                                               (1 - \sqrt{1 - \beta^\star})^2\}
                                               < r^\star < \beta^\star;
    \\
    \tilde{\Omega}(K^{1 - \beta^\star})
                                             &\text{ if }
                                               1/2 < \beta^\star \leq 1,
    \\
                                             &\text{ and } r^\star \geq \beta^\star.
  \end{cases}
\end{equation}
\end{proposition}
Proposition~\ref{prop:master-rates} shows, among other things, the rates (in
$K$) at which $\mathbf{E}^\star[\ln E_t^\star]$ diverges above the ``detection
boundary'' are exponential and that the same rates can be obtained for
$\mathbf{E}^\star[\ln E_t(\varepsilon,\delta)]$ provided that
$(\varepsilon, \delta)$ and $(\varepsilon^\star,\delta^\star)$ get close fast
enough in an asymptotic sense. The requirements that $d_TK^{1-\beta} = O(1)$ as
$K\to\infty$ and and $|d_\beta| \ln K \to 0$ are necessary to make sure that the cost of having
$(\varepsilon,\delta)\neq(\varepsilon^\star,\delta^\star)$ is of lower
asymptotic order.
With this result at hand,
Proposition~\ref{prop:misspecification-cost} can be shown.
\begin{proof}[Proof of Proposition~\ref{prop:misspecification-cost}]
  If $r^\star < \rho(\beta^\star)$, the fact that
  $\mathbf{E}^\star[\ln E_t(\varepsilon, \delta)]< \mathbf{E}^\star[\ln
  E_t^\star]$ for each $K$ and Theorem~\ref{thm:lr_achieves_boundary} shows that
  $\mathbf{E}^\star[\ln E_t(\varepsilon, \delta)]\to 0$ as $K\to\infty$. If
  $r^\star > \rho(\beta^\star)$, the result follows
  from Proposition~\ref{prop:master-rates}, as the exponents involved in~(\ref{eq:master-rates}) are positive.
\end{proof}
Furthermore, since Proposition~\ref{prop:master-rates} gives exponential rates
for $\mathbf{E}^\star[\ln E_t(\varepsilon, \delta)]$, it ensures that the
balancing act in~(\ref{eq:prior-balance}) is possible and
Theorem~\ref{thm:prior-works} can be proven.
\begin{proof}[Proof of Theorem~\ref{thm:prior-works}]
  By the construction of the grid $\mathcal{G}_C$, there exists
  $(\beta^\circ, \delta^\circ)\in \mathcal{G}_c$ such that
  $d_\beta := \beta^\circ - \beta^\star$ is nonnegative and such that
  $d_\beta \lesssim 1 / (2\ln^2 K)$. Let
  $\varepsilon^\circ := K^{-\beta^\circ}$. We now claim that there exists
  $T^\circ$ such that both
  $\delta^\circ := \sqrt{2(1 / T^\circ) \ln K}\in \mathcal{G}_C$, and
  $d_T := 1 - \sqrt{T^\circ / T^\star}$ is nonnegative and
  $d_T\sim K^{-(1 - \beta^\circ)}$ as $K\to \infty$. Indeed, there is $i$ and
  $i+1$ such that
  $\sqrt{2(1 / T_i)\ln K} =: \delta_i < \delta^\star \leq \delta_{i+1} :=
  \sqrt{2(1 / T_{i+1} )\ln K}$ with
  $T_i = \exp(i / \lceil K\varepsilon^\circ \rceil)$ and
  $T_{i+1} = \exp((i + 1) / \lceil K\varepsilon^\circ \rceil)$, and
  $\delta_i,\delta_{i+1}\in{\cal G}_C$. Then, observe that
  $\ln(\delta_{i}) -\ln(\delta_{i + 1}) = \ln(\sqrt{T_{i+1}/T_i}) = 1 /(2 \lceil
  K\varepsilon^\circ \rceil)\sim K^{-(1 - \beta^\circ)} /2$ . Taking
  $T^\circ = T_{i}$ satisfies the earlier claim. Hence, both
  $|d_\beta|\ln K \to 0$ as $K\to \infty$, and $|d_T|K^{1 - \beta^\circ} = O(1)$
  in the same limit. We have that
  \begin{equation*}
    \mathbf{E}^\star[\ln E_t(\mathbf{\Pi})]
    \geq
    \mathbf{E}^\star[\ln E_t(\varepsilon^\circ, \delta^\circ)]
    - \ln(|\mathcal{G}_C|).
  \end{equation*}
  The proof of Proposition~\ref{prop:misspecification-cost}
  (Proposition~\ref{prop:master-rates} in particular) shows that
  $\mathbf{E}^\star[\ln E_t(\varepsilon^\circ, \delta^\circ)] \gtrsim K^{f(t /
    T^\star, \beta^\star)}$ as $K\to\infty$ for a fixed positive function
  $f$---the growth is exponential. Since $\ln(|\mathcal{G}_C|) = O(\ln K)$ as
  $K\to\infty$, the previous display implies that
  $\mathbf{E}^\star[\ln E_t(\mathbf{\Pi})]\to\infty$ whenever
  $t/T^\star > \rho(\beta^\star)$. In the alternative case when
  $t/T^\star < \rho(\beta^\star)$, together the fact that
  $\mathbf{E}^\star[\ln E_t(\mathbf{\Pi})] \leq \mathbf{E}^\star[\ln E^\star_t]$
  and Theorem~\ref{thm:lr_achieves_boundary} show that
  $\mathbf{E}^\star[\ln E_t(\mathbf{\Pi})] \to 0$ as $K\to \infty$. This is what
  was to be shown.
\end{proof}

\subsection{Proof of Proposition~\ref{prop:master-rates}}
\label{sec:proof-master-rates}

\begin{proof}[Proof of Proposition~\ref{prop:master-rates}]
Let $\Delta = \sqrt{2r\ln K}$ with $r = t / T$ and let
$\Delta^\star = \sqrt{2r^\star\ln K}$ with $r^\star = t / T^\star$. Let
\begin{equation}\label{eq:loglikelihood}
\ell(x) := \ln(1 - \varepsilon + \varepsilon \rme^{\Delta x - \Delta ^ 2 / 2}).
\end{equation}
The function $\ell$ is the single-stream log-likelihood ratio between the null
hypothesis and the alternative determined by the parameter pair
$(\varepsilon, \Delta)$ evaluated at $x$. A number of properties of $\ell$ are
readily checked: the function $x\mapsto \ell(x)$ is differentiable, $\ell$ is
increasing and convex---provided $\varepsilon>0$---, $\ell$ takes the value $0$
at $x = \Delta / 2$, and it is always larger than $\ln(1 - \varepsilon)$. In
terms of $\ell$, the quantity of interest is written as
\begin{equation*}
  \mathbf{E}^\star[\ln E_t(\varepsilon, \delta)]
  =
  K\mathbf{E}^\star[\ln \ell(Z_{1, t})]
  =
  K\{
  (1 - \varepsilon^\star)\mathbf{E}_{0}[\ell(Z_{1, t})]
  +
  \varepsilon^\star\mathbf{E}_{0}[\ell(Z_{1, t} + \Delta^\star)]
  \},
\end{equation*}
where, recall from Section~\ref{sec:state-things},
$Z_{1,t} = (1 / \sqrt{t})\sum_{s\leq t}X_{1, s}$, and
$\delta^\star = \sqrt{2(1/T^\star)\ln K}$, and
$\Delta^\star = \sqrt{t}\delta^\star$. The quantity of interest consists of two
Gaussian integrals of translations of $\ell$. If $\Delta^\star$ and $\Delta$ are
close to each other, the function $x\mapsto \ell(x)$ has a zero close to
$x = \Delta^\star/2$ and the function $x\mapsto \ell(x + \Delta^\star)$ has a
zero close to $x = - \Delta^\star/2$. In order to align these two (near) zeros
and the main features of each term in the previous display, the authors found it
useful to translate and reflect $\ell$ around $x = \Delta^\star /2$. This tool
will be used repeatedly and it is stated in the following lemma---its proof is
found in Section~\ref{sec:proof-translate-center}.
\begin{lemma}\label{lem:translate-center}
  Let $((X_{i, t})_{t\geq 1}; i\in [K])\sim \mathbf{P}^\star$ and
  $Z_{i, t} = (1 / \sqrt{t})\sum_{s\leq t}X_{i, s}$ for each $i\in[K]$.
  Let $f:\mathbf{R}\to\mathbf{R}$ be such that
  $\mathbf{E}^\star[|f(Z_{1, t})|] < \infty$, and let
  $\Delta^\star = \sqrt{t}\delta^\star$. Then,
  $$\mathbf{E}^\star[f(Z_{1, t})] = \int_0^{\infty} i[f](x) \rmd x,$$
  where $i$ is the functional operator
  \begin{equation}\label{eq:integrand}
  i[f](x):=(1 - \varepsilon^\star)i_0[f](x) + \varepsilon^\star i_1[f](x)
  \end{equation} with
  \begin{equation}\label{eq:i0-def}
    i_0[f](x)
    :=
    \varphi(\Delta^\star / 2 + x)f(\Delta^\star / 2 + x)
    +
    \varphi(\Delta^\star / 2 - x)f(\Delta^\star / 2 - x)
  \end{equation}
  and
  \begin{equation}\label{eq:i1-def}
    i_1[f](x)
    :=
    \varphi(\Delta^\star / 2 - x)f(\Delta^\star / 2 + x)
    +
    \varphi(\Delta^\star / 2 + x)f(\Delta^\star / 2 - x).
  \end{equation}
\end{lemma}

Note that $i$, $i_0$ and $i_1$ are positive linear operators. Using Lemma~\ref{lem:translate-center}, we derive that
\begin{equation}
  \label{eq:representing}
  \mathbf{E}^\star[\ln E^\star_t] = K\int_0^\infty i[\ell](x)\rmd x,
\end{equation}
where $i$ is the functional defined in~(\ref{eq:integrand}) and $\ell$ is the
single-stream log-likelihood.

For the next step we split the integral in the right hand side
of~(\ref{eq:representing}) in two parts depending on whether $x < x_*$ or
$x\geq x_*$ for a carefully chosen value $x_*$. The choice of $x_*$ is driven by the
behavior of
$x\mapsto \ell(x + \Delta^\star/2) = \ln(1 - \varepsilon +
\varepsilon\rme^{\Delta x + (\Delta\Delta^{\star} - \Delta^2)/2})$ because, as
we will see, the terms containing this function constitute the main
contributions to (\ref{eq:representing}). We choose $x_*$ to be the point where
$\varepsilon\rme^{\Delta x_* + (\Delta\Delta^{\star} - \Delta^2)/2} = 1$; that
is,
$x_* = \ln(1 /\varepsilon) / \Delta + (1 / 2)(\Delta - \Delta^{\star})$.
Notice that if it was the case that $\Delta = \Delta^\star$, in the
$(\beta^\star,r^\star)$ parametrization we would have
$x_* = (\beta^{\star} / 2\sqrt{r^\star})\sqrt{2\ln K}$. However, since
$\Delta \neq \Delta^\star$ in general, and, recall from the statement of the
proposition, $d_T = 1 - \sqrt{r^\star / r} = o(1/\ln K)$, we have, as $K\to\infty$,
\begin{equation}\label{eq:def-star}
  x_* := \frac{\beta}{2\sqrt{r}}\sqrt{2\ln K} + \frac12 d_T\sqrt{2 r \ln K} =
  \frac{\beta^\star}{2\sqrt{r^\star}}\sqrt{2\ln K} + o(1/\sqrt{\ln K}).
\end{equation}
This means both that $x_*$ is not too far from the value it would take if
$\Delta = \Delta^\star$, and it means that the difference decays, as we will
see, sufficiently fast for our purposes.

In order to simplify the presentation of the results, the following analysis is
written in terms of $d_r := rd_T$. By the assumption that
$d_TK^{1-\beta} = O(1)$ as $K\to\infty$, we have that $d_r K^{1-\beta} = O(1)$
as $K\to\infty$ as well.

The next lemma, Lemma~\ref{lem:asymp-before-x0}, contains the analysis for
$x < x_*$. The analysis consists of rewriting $x\mapsto \ell(x)$ as as
$\ell(x) = \ln(1 + u(x))$ with
$u(x) := \varepsilon(\rme^{\Delta x - \Delta^2/2} - 1)$ and using the Taylor
expansion $\ln(1 + z) = z + z^2/2 + O(z^3)$ as $z\to 0$. Indeed, we expand
$\ell(x) = \ln(1 + u(x)) = u(x) - \frac12 u^2(x) + O(u^3(x))$, and control the
$O(u^3)$-error term explicitly. This line of development requires analyzing both
$I_1 := K\int_0^{x_*}i[u](x)\rmd x$ and $I_2 := K\int_0^{x_*}i[u^2](x)\rmd x$
asymptotically as $K\to\infty$, which requires the bulk of the technical efforts
in this proof.
The analysis shows that the term $I_1$ contains the main cost of having
$\Delta\neq \Delta^\star$---the cost of misspecification---and the requirement
that $K^{1-\beta}d_r = O(1)$ as $K\to\infty$ is essential to maintain said cost
of lower asymptotic order in the same limit. In the analysis of $I_2$, it is
sufficient that $d_r\ln K, d_\beta\ln K = o(1)$ as $K\to\infty$. The following
lemma contains the result of this analysis and it is proven in
Section~\ref{sec:proof-asymp-before-x0}.
\begin{lemma}\label{lem:asymp-before-x0}
  Adopt the assumptions of Lemma~\ref{prop:master-rates}; in particular, recall,
  $d_rK^{1 - \beta} = O(1)$ and $d_\beta\ln K = o(1)$ as $K\to\infty$. Let
  $\ell$ be as in~(\ref{eq:loglikelihood}), let $i$ be as
  in~(\ref{eq:integrand}), and let $x_*$ be as in~(\ref{eq:def-star}). Then, as $K\to\infty$,
  \begin{equation}\label{eq:ell-lower-rates}
    K\int_0^{x_*}i[\ell](x)\rmd x
    =
    \begin{cases}
      \Omega(K^{1 + 2r^\star - 2\beta^\star})  &\text{ if }
                                                 1/2 < \beta^\star < 3/4,
                                                 \\
                                               & \text{ and } \beta^\star - 1/2
                                                 < r^\star \leq \beta^\star / 3;
      \\
      \tilde{\Omega}(
      K^{1  - \beta^\star -
      \frac{1}{4}\paren{\frac{r^\star
      - \beta^\star}{\sqrt{r^\star}}}^
      2}
      )
                                               &\text{ if }
                                                 1/2 < \beta^\star \leq 1,
      \\
                                               &\text{ and }
                                                 r^\star > \max\{\beta^\star/3,
                                                 (1 - \sqrt{1 - \beta^\star})^2\}.
    \end{cases}
  \end{equation}
\end{lemma}

The next lemma, Lemma~\ref{lem:asymp-after-x0}, contains the analysis
of~(\ref{eq:representing}) for $x \geq x_*$. Its proof uses that
$x\mapsto \ell(x)$ is increasing for $x\geq x_*$ and that consequently
$\ell(x \pm \Delta^\star /2)\geq \ell(x_* \pm \Delta^\star /2)$ as long as
$x\geq x_*$. Analyzing the remaining Gaussian integrals yields the following
lemma, which is proven in Section~\ref{sec:u-bound-after-x0}.
\begin{lemma}\label{lem:asymp-after-x0}
  Adopt the assumptions of Lemma~\ref{prop:master-rates} with a modification:
  instead of assuming that $d_rK^{1 - \beta} = O(1)$ and $d_\beta\ln K = o(1)$
  as $K\to\infty$, assume that $d_r\ln K, d_\beta\ln K \to 0$ in the same limit
  (a weaker assumption). Let $\ell$ be as in~(\ref{eq:loglikelihood}), let $i$
  be as in~(\ref{eq:integrand}), and let $x_*$ be as in~(\ref{eq:def-star}).
  Then, as $K\to\infty$
  \begin{equation}
    \label{eq:asymp-after-x0}
    K\int_{x_*}^{\infty}i[\ell](x) \rmd x
    =
    \begin{cases}
      \tilde{\Omega}(K^{1 - \beta^\star - \frac14\paren{\frac{r^\star - \beta^\star}{\sqrt{r^\star}}}^2})
      &\text{ if } r^\star <  \beta^\star;\\
      \tilde{\Omega}(K^{1 - \beta^\star})
      &\text{ if } r^\star \geq \beta^\star.
    \end{cases}
  \end{equation}
\end{lemma}
Using the Lemma~\ref{lem:asymp-before-x0} and Lemma~\ref{lem:asymp-after-x0}
and~(\ref{eq:representing})---in particular, adding~(\ref{eq:ell-lower-rates})
and~(\ref{eq:asymp-after-x0})---, we obtain (\ref{eq:master-rates}), concluding the proof.

\end{proof}

\subsubsection{Proof of Lemma~\ref{lem:translate-center}}
\label{sec:proof-translate-center}
\begin{proof}[Proof of Lemma~\ref{lem:translate-center}]
  Notice that
  $Z_{1,t}\sim (1 - \varepsilon^\star)\mathrm{Normal}(0, 1) +
  \varepsilon^\star\mathrm{Normal}(\Delta^\star, 1)$ and that
  $\mathbf{E}^\star[f(Z_{1, t})] = (1 - \varepsilon^\star)\mathbf{E}_0[f(Z_{1,
    t})] + \varepsilon^\star\mathbf{E}_{\delta^\star}[f(Z_{1, t})]$. We
  translate the integrands by $\Delta^\star/2$, split the integrals in two, and
  perform a reflection. For the first expectation we have that
  \begin{align*}
    \mathbf{E}_{0}[f(Z_{1, t})]
    &= \int_{-\infty}^{\infty} \varphi(x)f(x)\rmd x\\
    &= \int_{-\infty}^{\infty} \varphi(\Delta^\star/2 + x)f(\Delta^\star/2 + x)\rmd x\\
    &= \int_{-\infty}^{0} \varphi(\Delta^\star/2 + x)f(\Delta^\star/2 + x)\rmd x
      + \int_{0}^{\infty} \varphi(\Delta^\star/2 + x)f(\Delta^\star/2 + x)\rmd x\\
    &= \int_{0}^{\infty} \varphi(\Delta^\star/2 - x)f(\Delta^\star/2 - x)\rmd x
      + \int_{0}^{\infty} \varphi(\Delta^\star/2 + x)f(\Delta^\star/2 + x)\rmd x\\
    &=
      \int_{0}^{\infty} \bracks{
      \varphi(\Delta^\star/2 + x)f(\Delta^\star/2 + x) +
      \varphi(\Delta^\star/2 - x)f(\Delta^\star/2 - x)
      }\rmd x.
  \end{align*}
  The integrand is exactly $i_0[f]$. Starting with
  $\mathbf{E}_{\delta^\star}[f(X)] = \int_{-\infty}^{\infty}\varphi(x -
  \Delta^\star)f(x)\rmd x$ and repeating the same steps, one derives $i_1[f]$
  as claimed.
\end{proof}

\subsubsection{Proof of Lemma~\ref{lem:asymp-before-x0}}
\label{sec:proof-asymp-before-x0}

\begin{proof}[Proof of Lemma~\ref{lem:asymp-before-x0}.]
  Recall that $u(x) = \varepsilon(\rme^{\Delta x - \Delta^2 / 2} - 1)$ and that
  $\ell(x) = \ln(1 + u(x))$. We bound $\int_0^{x_*}i[\ell](x)\rmd x$ from below
  using a lower bound on $z\mapsto\ln(1 + z)$. Indeed, using the Taylor expansion
  of the natural logarithm and controlling the terms larger than third order, the
  elementary inequality
  \begin{equation*}
    \ln(1 + z) \geq z - \frac12 z^2 + \frac{\min\{0,z^3\}}{3(1+\min\{0,z\})}
  \end{equation*}
  holds for any $z>-1$. Note that the last term in the previous bound is
  increasing in $z$, which implies that, since $u(x)\geq - \varepsilon$ and
  $\ell(x) = \ln(1 + u(x))$,
  \begin{equation}\label{eq:ell-lowebound}
    \ell(x) \geq u(x) - \frac12 u(x)^2 - \frac{\varepsilon^3}{3(1-\varepsilon)}.
  \end{equation}
  Next, we apply $i$, the operator from \eqref{eq:integrand}, on both sides
  of~\eqref{eq:ell-lowebound} and integrate in $[0, x_*]$. To this end, recall
  that, from the definition of $i$, for a function $f$ and a constant $C$,
  \begin{equation}\label{eq:i-function-constant}
    i[f+C](x) = i[f](x) + C\left\{\varphi(\Delta^\star/2 + x) + \varphi(\Delta^\star/2-x)\right\} .
  \end{equation}
  With this in mind, define $I_1 := K\int_0^{x_*}i\sqbrack{u}(x)\rmd x$ and
  $I_2 := K\int_0^{x_*}i\sqbrack{u^2}(x)\rmd x$. Using the linearity of the
  operator $i$ along with~\eqref{eq:ell-lowebound} and~\eqref{eq:i-function-constant},
  \begin{align}
    K\int_0^{x_*}i[\ell](x)\rmd x
    &\geq I_1 -   \frac12 I_2
      - \frac{K\varepsilon^3}{3(1-\varepsilon)}\int_0^{x_*}\left\{\varphi(\Delta^\star/2 + x) + \varphi(\Delta^\star/2-x)\right\}\rmd x\nonumber\\
    &\geq I_1-\frac12 I_2-\sqrt{\frac{2}{\pi}}
      \frac{x_* K^{1 - 3\beta}}{3 (1 - K^{-\beta})},\label{eq:ell-lbd}
  \end{align}
  where the second inequality follows from the crude bound $\varphi(x)\leq 1/\sqrt{2\pi}$. Note the last term is $\tilde{\Theta}(K^{1 - 3\beta})$ as $K\to\infty$; it tends to zero as $K\to\infty$.

  We continue with the asymptotic analysis of $I_1$ and $I_2$ in
  Lemma~\ref{lem:linear-lemma} and Lemma~\ref{lem:quadratic-lemma},
  respectively. Lemma~\ref{lem:quadratic-lemma} gives an exact asymptotic
  characterization of $I_2$ while only an asymptotic lower bound for $I_2$ is
  needed to show Lemma~\ref{lem:asymp-before-x0}; that is,
  Lemma~\ref{lem:quadratic-lemma} is stronger than is currently needed. The
  reason for this (and for splitting the analysis in
  Lemma~\ref{lem:linear-lemma} and Lemma~\ref{lem:quadratic-lemma}) is that
  Lemma~\ref{lem:quadratic-lemma} is also necessary in the proof of
  Lemma~\ref{lem:power-crux} (where an asymptotic upper bound for $I_2$ is
  needed). We remark that in order to control the difference $I_1 - \frac12 I_2$
  appearing in~(\ref{eq:ell-lbd}), it is necessary to keep track of all leading
  constants in the asymptotic expressions for both $I_1$ and $I_2$ because these
  two quantities are of the same asymptotic order in certain regimes. To this
  end, the asymptotic ``$\sim$'' notation is used (see
  Section~\ref{sec:notation-proofs}).

  \paragraph*{Analysis of $I_1$} In Lemma~\ref{lem:linear-lemma}, the definite
  integral $I_1$ is analyzed by quantifying how far its value is from the value
  it would have if there was no misspecification. Indeed, if $d_r\ln K = o(1)$
  and $d_\beta\ln K = o(1)$ (weaker conditions than stated in
  Lemma~\ref{lem:asymp-before-x0}, which we are proving) then $I_1$ equals
  approximately $I^\star_1 := K\int_0^{x_*}i\sqbrack{u^\star}(x)\rmd x$ with
  $u^\star(x) := \varepsilon^\star(\rme^{\Delta^\star x - \Delta^{\star 2} / 2}
  - 1)$. Note that $u(x)$ would be equal to its starred version, $u^\star(x)$, if
  $(\varepsilon, \Delta)$ was equal to $(\varepsilon^\star, \Delta^\star)$; in
  other words, if there was no misspecification. The difference between $I_1$
  and $I_1^\star$, understood as an approximation error, is denoted by $\Xi_K$
  in \eqref{eq:I1-and-I1star}; it captures the cost of ``misspecification'' when
  $\Delta\neq \Delta^\star$ and $\varepsilon \neq \varepsilon^\star$.
  Corollary~\ref{cor:linear-corollary} of Lemma~\ref{lem:linear-lemma} shows,
  under the stronger assumption that $K^{1-\beta}d_r = O(1)$ as $K\to\infty$,
  that the misspecification cost $\Xi_K$ is of lower order than $I_1^\star$.
  Additionally, both an asymptotic expressions for $I^\star_1$ and an asymptotic
  upper bound for $|\Xi_K|$ are derived as $K\to\infty$; the proof is contained
  in Section~\ref{sec:proof-linear-lemma}.
  \begin{lemma}\label{lem:linear-lemma}
    Recall that $u(x) = \varepsilon(\rme^{\Delta x - \Delta^2 / 2} - 1)$, let
    $u^\star(x) = \varepsilon^\star(\rme^{\Delta^\star x - \Delta^{\star 2} / 2}
    - 1)$, let $I_1 = K\int_0^{x_*}i[u](x)\rmd x$, and let
    $I_1^\star = K\int_0^{x_*}i[u^\star](x)\rmd x$, where the operator $i$ is
    defined in~(\ref{eq:integrand}) and $x_*$ is defined in~(\ref{eq:def-star}).
    Then, as $K\to\infty$,
    \begin{equation}\label{eq:linear_rates}
      I_1^\star
      \sim
      \begin{cases}
        K^{1 + 2r^\star - 2\beta^\star}  &\text{ if }
                                           1/2 < \beta^\star < 3/4
                                           \text{ and }
                                           \beta^\star - 1/2 \leq r^\star < \beta^\star / 3;\\
        \frac12 K^{1 + 2r^\star - 2\beta^\star} &\text{ if }
                                                  1/2 < \beta^\star \leq 3/4
                                                  \text{ and }
                                                  r^\star = \beta^\star/3; \\
        \frac{1}{\sqrt{2\pi}}
        \frac{2\sqrt{r^\star}}{3r^\star
        - \beta^\star}\frac{K^{1  - \beta^\star -
        \frac{1}{4}\paren{\frac{r^\star
        - \beta^\star}{\sqrt{r^\star}}}^
        2}}{\sqrt{2\ln K}}
                                         &\text{ if }
                                           1/2 < \beta^\star < 3/4
                                           \text{ and }
                                           r^\star > \beta^\star/3,\\
                                         &\text{ or if }
                                           3/4 \leq \beta^\star \leq 1
                                           \text{ and }
                                           r^\star \geq (1 - \sqrt{1 - \beta^\star})^2.
      \end{cases}
    \end{equation}
    Additionally, if both $d_r\ln K = o(1)$ and $d_\beta\ln K = o(1)$ as
    $K\to\infty$, then, as $K\to \infty$,
    \begin{equation}\label{eq:I1-and-I1star}
      I_1
      =
      I_1^\star
      (1 + o(1))
      -
      \Xi_K
    \end{equation}
    for a real-valued sequence $(\Xi_K)_{K\geq 1}$ such that, in the same limit,
    \begin{equation}\label{eq:xi-asymp}
      |\Xi_K| \lesssim
      \begin{cases}
        K^{1 - \beta^\star}
        |d_r|\ln K
        &\text{ if } r^\star <
          \beta^\star;\\
        \frac12 K^{1 - \beta^\star}
        |d_r|\ln K
        &\text{ if } r^\star =
          \beta^\star;\\
        \frac{1}{\sqrt{2\pi}}
        \frac{2\sqrt{r^\star}}{r^\star
        - \beta^\star}
        \frac{K^{1 - \beta^\star - \frac14\paren{\frac{r^\star - \beta^\star}{\sqrt{r^\star}}}^2}}{\sqrt{2\ln K }}
        |d_r|\ln K
        &\text{ if } r^\star >  \beta^\star.\\
      \end{cases}
    \end{equation}
  \end{lemma}
  With the asymptotic expression~(\ref{eq:I1-and-I1star}) for $I_1$ in terms of
  $I_1^\star$ and $\Xi_K$, and the asymptotic upper bound for the absolute
  misspecification cost $|\Xi_K|$, the assumption that $d_rK^{1-\beta} = O(1)$
  is now motivated: under this additional assumption, $\Xi_K$ is of lower order
  than $I^\star_1$ in~(\ref{eq:I1-and-I1star}) and in consequence $I_1$ is
  asymptotic to the same expression as $I_1^\star$ from~(\ref{eq:linear_rates})
  as $K\to\infty$.
  \begin{corollary}\label{cor:linear-corollary}
    Let $I_1$ and $I_1^\star$ be as in Lemma~\ref{lem:linear-lemma}. Under that
    lemma's assumptions, if $r^\star > \rho(\beta^\star)$, and if, as
    $K\to\infty$, both $d_{r}K^{1 - \beta^\star} = O(1)$ and
    $d_\beta \ln K = o(1)$, then $I_1 \sim I_1^{\star}$ as $K\to \infty$.
    Consequently, $I_1$ is asymptotic to the right hand side
    of~(\ref{eq:linear_rates}).
  \end{corollary}
  \begin{proof}[Proof of Corollary~\ref{cor:linear-corollary}]
    Lemma~\ref{lem:linear-lemma} shows that, whenever
    $r^\star > \rho(\beta^\star)$,
    $I^\star_1 = \tilde{\Theta}(K^{f(\beta^\star, r^\star)})$ for a nonnegative
    $f$ given by~(\ref{eq:linear_rates}), which is independent of $K$, and it
    shows that $I^\star_1\to\infty$ in the same case. Notice that our assumption
    on $d_r = d_r(K)$ that $K^{1 - \beta^\star}d_r = O(1)$ as $K\to\infty$
    implies that $d_r\ln K = o(1)$ in the same limit; the assumption of
    Lemma~\ref{lem:linear-lemma} on $d_r$ is satisfied.
    Hence,~(\ref{eq:I1-and-I1star}) holds and, in order to conclude that
    $I_1 \sim I_1^\star$ as $K\to\infty$, it is necessary to analyze $\Xi_K$
    asymptotically as $K\to \infty$. We do this using our assumption on $d_r$
    (that $d_{r}K^{1 - \beta^\star} = O(1)$ as
    $K\to\infty$),~(\ref{eq:I1-and-I1star}), and~(\ref{eq:xi-asymp}). Two cases
    are needed: (1) if $r^\star\leq \beta^\star$, then , as $K\to\infty$,
    $\Xi_K = O(\ln K)$ and consequently $\Xi_K = o(I_1^\star)$ because
    $I_1\sim K^{f(\beta^\star, r^\star)}$ in the same limit; (2) if
    $r^\star > \beta^\star$, the same conclusion holds because, as $K\to\infty$,
    $I^\star_1 = \Theta(K^{f(\beta^\star, r^\star)}(\ln K)^{-1/2})$ and
    $\Xi_K = O(K^{f(\beta^\star, r^\star)}d_r(\ln K)^{1/2})$, and
    $d_r(\ln K)^{1/2} = o((\ln K)^{-1/2})$ as $K\to\infty$ by our assumption.
    This implies the claim.
  \end{proof}

  \paragraph*{Analysis of $I_2$} For this analysis, contained in
  Lemma~\ref{lem:quadratic-lemma}, it is enough to have $d_r\ln K = o(1)$ and
  $d_\beta\ln K = o(1)$ as $K\to\infty$. Contrary to the asymptotic expression
  derived for $I_1$ in \eqref{eq:I1-and-I1star}, misspecification does not add
  terms to $I_2$ asymptotically under these assumptions. The proof of the following
  lemma is contained in Section~\ref{sec:proof-quadratic-lemma}.
  \begin{lemma}\label{lem:quadratic-lemma}
    Let $I_2 = K\int_0^{x_*}i[u^2](x)\rmd x$, where
    $u(x) = \varepsilon(\rme^{\Delta x - \Delta^2/2} - 1)$, the operator $i$ is
    defined in (\ref{eq:integrand}) and $x_*$ is defined in~(\ref{eq:def-star}).
    Adopt the assumptions of Lemma~\ref{lem:asymp-before-x0} with a
    modification: instead of assuming that $d_rK^{1 - \beta} = O(1)$ and
    $d_\beta\ln K = o(1)$ as $K\to\infty$, assume that $d_r\ln K = o(1)$ and
    that $d_\beta\ln K = o(1)$ in the same limit (a weaker assumption). Then, as
    $K\to\infty$,
    \begin{equation}\label{eq:quadratic-rates}
      I_2 \sim
      \begin{cases}
        K^{1 + 2r^\star - 2\beta^\star}
            &\text{if } 1/2 < \beta^\star < 3/4 \text{ and } \\
            &\beta^\star - 1/2 \leq r^\star < \beta^\star / 3;
        \\
        \frac12 K^{1 + 2r^\star - 2\beta^\star}
            &\text{if } 1/2 < \beta^\star \leq 3/4 \text{ and } r^\star = \beta^\star / 3;\\
        \frac{C(r^\star, \beta^\star)}{\sqrt{2\pi}}
        \frac{
        K^{1 - \beta^\star -
        \frac14\paren{\frac{\beta^\star - r^\star}{\sqrt{r^\star}}}^2}
        }{\sqrt{2\ln K}}
            &\text{if } 1/2 < \beta^\star \leq 3/4
              \text{ and } r^\star > \beta^\star / 3, \\
            &\text{or if } 3/4 < \beta^\star \leq 1 \text{ and }\\
            & r^\star \geq (1 - \sqrt{1-\beta^\star})^2
      \end{cases}
    \end{equation}
    with $C(r^\star, \beta^\star) := 2\sqrt{r^\star}
        \left[
        (3r^\star - \beta^\star)^{-1}
        +
        (5r^\star - \beta^\star)^{-1}
        \right]$.
  \end{lemma}

Using this asymptotic analysis of $I_1$ and $I_2$ in~(\ref{eq:ell-lbd}), we can finish the proof of Lemma~\ref{lem:asymp-before-x0}. By Corollary~\ref{cor:linear-corollary}, the rates in~(\ref{eq:linear_rates})
  are the same as those of $I_1 = K\int_0^{x_*}i[u](x)\rmd x$. Furthermore (\ref{eq:quadratic-rates}) give the rates for
  $I_2 =K\int_0^{x_*}i[u^2](x)\rmd x$. Putting these two together,
  \begin{equation}\label{eq:approx-ell-lower-rates}
    K\int_0^{x_*}i\sqbrack{u - \frac12 u^2}(x)\rmd x
    \gtrsim
    \begin{cases}
       \frac12 K^{1 + 2r^\star - 2\beta^\star}      &\text{if } 1/2 < \beta^\star < 3/4 \text{ and }\\
                                                   &\beta^\star - 1/2 < r^\star < \beta^\star/3;\\
       \frac14 K^{1 + 2r^\star - 2\beta^\star}      &\text{if } 1/2 < \beta^\star \leq 3/4\text{ and } \\
                                                   &r^\star = \beta^\star/3;\\
      \frac{C'(r^\star, \beta^\star)}{\sqrt{2\pi}}
      \frac{K^{1  - \beta^\star -
      \frac{1}{4}\paren{\frac{r^\star
      - \beta^\star}{\sqrt{r^\star}}}^
      2}}{\sqrt{2\ln K}}  &\text{if } 1/2 < \beta^\star < 3/4 \text{ and }\\
                                                   & r^\star > \beta^\star/3,\\
                                               &\text{or if }
                                                 3/4 \leq \beta^\star \leq 1 \text{ and }\\
                                                 & r^\star > (1 - \sqrt{1 - \beta^\star})^2,
    \end{cases}
    \end{equation}
  where, in the last of the cases, $C'(r^\star, \beta^\star) :=   \sqrt{r^\star}\left[
      (3r^\star
      - \beta^\star)^{-1}
      -
      (5r^\star
      - \beta^\star)^{-1}
      \right]$ is positive.
  This shows that $K\int_0^{x_*}i\sqbrack{u - \frac12 u^2}(x)\rmd x \to \infty$
  as $K\to\infty$ anytime that $r^\star > \rho(\beta^\star)$. Since the last
  term in~(\ref{eq:ell-lbd}) is of lower order (it is
  $\tilde{\Theta}(K^{1 - 3\beta})$ as $K\to\infty$ and it consequently tends to
  zero in the same limit), the result follows
  as~(\ref{eq:approx-ell-lower-rates}) and~(\ref{eq:ell-lbd}) now
  imply~(\ref{eq:ell-lower-rates}).
\end{proof}

\subsubsection{Proof of Lemma~\ref{lem:asymp-after-x0}}
\label{sec:u-bound-after-x0}
\begin{proof}[Proof of Lemma~\ref{lem:asymp-after-x0}]
  The proof strategy in this lemma is different than the one used for
  Lemma~\ref{lem:asymp-before-x0}. Rather than using the representation from Lemma~\ref{lem:translate-center}, it is more convenient to rewrite the integrand \eqref{eq:integrand} $i$ as
  $i[\ell](x) = i_+(x) + i_-(x)$ with
  \begin{equation*}
    i_+(x) =
    \ell(\Delta^\star / 2 + x)
    [(1 - \varepsilon^\star)\varphi(\Delta^\star / 2 + x)
    +
    \varepsilon^\star\varphi(\Delta^\star / 2 - x)],
  \end{equation*}
  and
  \begin{equation*}
    i_-(x)
    =
    \ell(\Delta^\star / 2 - x)[
    (1-\varepsilon^\star)\varphi(\Delta^\star / 2 - x)
    +
    \varepsilon^\star\varphi(\Delta^\star / 2 + x)].
  \end{equation*}
  On $i_-$, since $x\mapsto \ell(\Delta^\star/2 - x)$ is
  decreasing, we conclude that $\ell(\Delta^\star/2 - x) \geq \ln(1-\varepsilon)$, its value
  at infinity. Hence,
  \begin{equation*}
    i_-(x)
    \geq
    \ln(1 - K^{-\beta})
    \sqbrack{
    (1-\varepsilon^\star)\varphi(\Delta^\star / 2 - x)
    +
    \varepsilon^\star\varphi(\Delta^\star / 2 + x)
    }.
  \end{equation*}
  Integrating, one obtains that
  \begin{equation}\label{eq:iminus-bound}
    \int_{x_*}^{\infty} i_-(x) \rmd x
    \geq
      \ln(1 - K^{-\beta})
      [
      (1-\varepsilon^\star)\overline{\Phi}(x_* - \Delta^\star / 2)
      +
      \varepsilon^\star\overline{\Phi}(x_* + \Delta^\star / 2)
      ].
  \end{equation}
  On $i_+$ we use that $\ell$ is increasing; therefore, for any $x \geq x_*$, we
  have that
  $\ell(x + \Delta^\star / 2)\geq \ell(x_* + \Delta^\star /
  2) \geq \ln(2-\varepsilon)$. Using this bound and integrating shows that
  \begin{equation}    \label{eq:iplus-bound}
    \int_{x_*}^\infty i_+(x)\rmd x
    \geq
    \ln(2 - K^{-\beta})
    [
    (1-\varepsilon^\star)\overline{\Phi}(x_* + \Delta^\star / 2)
    +
    \varepsilon^\star\overline{\Phi}(x_* - \Delta^\star / 2)
    ].
  \end{equation}
  Adding~(\ref{eq:iminus-bound}) and~(\ref{eq:iplus-bound}) and regrouping, we
  obtain that
  \begin{multline}\label{eq:i-bound-after-x0}
    \int_{x_*}^\infty i(x)\rmd x
    \geq
    \overline{\Phi}(x_* - \Delta^\star / 2)
    [
    (1-\varepsilon^\star)\ln(1 -
    K^{-\beta})
    +
    \varepsilon^\star \ln(2 - K^{-\beta})
    ]
    \\
    +
    \overline{\Phi}(x_* + \Delta^\star / 2)
    [
    (1-\varepsilon^\star) \ln(2 - K^{-\beta})
    +
    \varepsilon^\star\ln(1 - K^{-\beta})
    ]
    .
  \end{multline}
  With an eye on the previous display, notice that asymptotically as
  $K\to\infty$, we have that $\ln(2 - K^{-\beta})\to \ln(2)$ and
  $\ln(1 - K^{-\beta})\sim - K^{-\beta}$, which, under our assumptions on
  $d_\beta$, implies that $\ln(1 - K^{-\beta}) \sim -K^{-\beta^\star}$. Also,
  notice that
  \begin{equation*}
    \overline{\Phi}\left(x_* + \Delta^\star / 2\right) =
    \tilde{\Theta}(K^{-\frac14\paren{\frac{\beta^\star + r^\star}{\sqrt{r}}}^2}) =
    \tilde{\Theta}(K^{-\beta^\star - \frac14\paren{\frac{r^\star - \beta^\star
          }{\sqrt{r}}}^2}).
  \end{equation*}
For $\overline{\Phi}(x_* - \Delta^\star / 2)$, there are two cases leading to different conclusions concerning~(\ref{eq:i-bound-after-x0}): (i) if $\beta^\star > r^\star$, then
\begin{equation*}
  \overline{\Phi}(x_* - \Delta^\star / 2) = \tilde{\Theta}(K^{-\frac14\paren{\frac{\beta^\star - r^\star}{\sqrt{r^\star}}}^2})
\end{equation*}
and, therefore, using (\ref{eq:i-bound-after-x0}) we conclude that
  \begin{multline*}
    K\int_{x_*}^{\infty}i(x)
    \rmd x
     \\
    \begin{aligned}
      &\gtrsim K\left\{\overline{\Phi}(x_* + \Delta^\star / 2)
      [
      \ln(2) - K^{-2\beta^\star}
      ]
      +
      \overline{\Phi}(x_* - \Delta^\star / 2)
      [
      K^{-\beta^\star}(\ln(2) -K^{-\beta^\star})
      ]\right\}
      \\
      &\sim K\left\{\ln(2)\overline{\Phi}(x_* + \Delta^\star / 2)
        +
        \ln(2)K^{-\beta^\star}\overline{\Phi}(x_* - \Delta^\star / 2)\right\}
      \\
      &=
        \tilde{\Theta}(K^{1-\beta^\star - \frac14\paren{\frac{\beta^\star - r^\star}{\sqrt{r^\star}}}^2});
    \end{aligned}
  \end{multline*}
  (ii) if $\beta^\star \leq r^\star$, then
  $\overline{\Phi}(x_* - \Delta^\star / 2) = \Theta(1)$, and, because of both the
  previous observations and~(\ref{eq:i-bound-after-x0}),
  \begin{equation*}
    K\int_{x_*}^{\infty}i(x) \rmd x
    = \tilde{\Omega}(K^{1-\beta^\star}).
  \end{equation*}
  The last two displays imply~(\ref{eq:asymp-after-x0}), which is what was to be
  shown.
\end{proof}

  \subsubsection{Proof of Lemma~\ref{lem:linear-lemma}}
  \label{sec:proof-linear-lemma}
\begin{proof}[Proof of Lemma~\ref{lem:linear-lemma}]
  Let $u(x) = \varepsilon(\rme^{\Delta x - \Delta^2 / 2} - 1)$. We compute
  $i[u]$, where $i$ is given in Lemma~\ref{lem:translate-center}. There are
  several ways to perform this calculation, which involves computing all
  possible products $u(x \pm \Delta^\star / 2)\varphi(x \pm \Delta^\star/2)$. To
  this end, it is necessary to have both the parameters $\varepsilon$ and
  $\Delta$, and the variable $x$ in a common parametrization. Our choice is to
  write everything in terms of $K$ through the parametrization
  $x = x(\xi) = \sqrt{2 \xi \ln K}$ with $\xi \geq 0$. Notice that this relation
  is invertible for $x\geq 0$ and that we can write $\xi = \xi(x) = x^2 / (2 \ln K)$.
  Define $m(\xi) := 2 d_r\sqrt{\xi}/r$, which implies that
  $2\sqrt{\xi r} = 2\sqrt{\xi r^\star} + m(\xi)$. Notice that $m(\xi)$ is
  proportional to $d_r$. We first compute $i_0[u]$ and $i_1[u]$ from
  \eqref{eq:i0-def} and \eqref{eq:i1-def}. After the dust has settled,
  \begin{multline}\label{eq:linear-lemma-i0}
    \sqrt{2\pi}i_0[u](x)
    =
    K^{-(\sqrt{\xi} - \sqrt{r^\star} / 2)^2 + m(\xi)- (\beta + d_r)}
    -
    K^{-(\sqrt{r^\star}/2 + \sqrt{\xi})^2 - \beta}
    \\
    +
    K^{-(\sqrt{\xi} + \sqrt{r^\star} / 2)^2 - m(\xi)- (\beta + d_r)}
    -
    K^{-(\sqrt{r^\star}/2 - \sqrt{\xi})^2 - \beta},
  \end{multline}
and
\begin{multline}\label{eq:linear-lemma-i1}
  \sqrt{2\pi}i_1[u](x)
  =
  K^{-(\sqrt{\xi} - 3 \sqrt{r^\star} / 2)^2 + 2 r^\star + m(\xi)- (\beta + d_r)}
  -
  K^{-(\sqrt{r^\star}/2 - \sqrt{\xi})^2 - \beta}
  \\
  +
  K^{-(\sqrt{\xi} + 3\sqrt{r^\star} / 2)^2 + 2r^\star - m(\xi)- (\beta + d_r)}
  -
  K^{-(\sqrt{r^\star}/2 + \sqrt{\xi})^2 - \beta}.
\end{multline}
Recalling the definition of the functional $i$ in \eqref{eq:integrand}, using
\eqref{eq:linear-lemma-i0} and \eqref{eq:linear-lemma-i1} and reverting back to
the original $x$ parametrization, we find that
\begin{multline}
  i[u](x)  =
    \varphi(x - 3\Delta^\star/2)
  K^{2 r^\star - 2\beta^\star + m(\xi)- (d_\beta + d_r)}
  -
  \varphi(x - \Delta^\star/2)
  K^{- 2\beta^\star - d_\beta}
  \\
   +
  \varphi(x + 3\Delta^\star/2)
  K^{2r^\star - 2\beta^\star - m(\xi)- (d_\beta + d_r)}
  -
  \varphi(x + \Delta^\star/2)
  K^{- 2\beta^\star - d_\beta}
    - \eta_K(x) ,
 \label{eq:iu}
\end{multline}
where, recall, $\xi = \xi(x)$, and we defined the last term $\eta_K(x)$ as
\begin{multline}\label{eq:xi-precursor}
	  \eta_K(x) :=
	\varphi(x - \Delta^\star/ 2)K^{- \beta}
	(1 - \varepsilon^\star)
	(
	1
	-
	K^{m(\xi) - d_r}
	)
	\\
	+
	\varphi(\Delta^\star/2 + x)
	K^{- \beta}
	(1 - \varepsilon^\star)
	(
	1
	-
	K^{- m(\xi) - d_r}
	).
\end{multline}
The terms in \eqref{eq:iu} are arranged to make it easier to compare $i[u](x)$
and $i[u^\star](x)$. Indeed, $i[u^\star](x)$ is obtained from \eqref{eq:iu} by
setting $d_r = d_\beta = 0$, which implies that $m(\xi) = 0$ and that
consequently $\eta_K(x)=0$. Hence,
\begin{multline}
  \label{eq:iustar}
  i[u^\star](x)  =
  \varphi(x - 3\Delta^\star/2)
  K^{2 r^\star - 2\beta^\star}
  -
  \varphi(x - \Delta^\star/2)
  K^{- 2\beta^\star}
  \\
  +
  \varphi(x + 3\Delta^\star/2)
  K^{2r^\star - 2\beta^\star}
  -
  \varphi(x + \Delta^\star/2)
  K^{- 2\beta^\star}.
\end{multline}
Since $\eta_K\neq 0$ when $d_r, d_\beta\neq 0$, we may interpret $\eta_K$ as
being caused by misspecification; in fact, $\eta_K$ is the precursor of the term
$\Xi_K$ from the claim.

Looking to prove~(\ref{eq:I1-and-I1star}), it might be tempting to jump to a
premature conclusion since the first four terms~(\ref{eq:iu})
and~(\ref{eq:iustar}) coincide, for fixed $x$, up to $1 + o(1)$ factors as
$K\to\infty$. Indeed, since $d_\beta, d_r = o(1 / \ln K)$ as $K\to\infty$, then
$m(\xi) = o(1 / \ln K)$ in the same limit, and, consequently each of the first
four terms of (\ref{eq:iu}) is asymptotic its respective ``well-specified
version'' in~(\ref{eq:iustar}). This is because, if a sequence $c = c(K)$ is
such that $c = o(1 / \ln K)$, then $K^c = \exp(c \ln K) \to 1$ as $K\to\infty$.
However, two obstacles still need to be overcome in this argument: (1) the
``$1 + o(1)$'' in involved in the asymptotic equivalence statement for each term
is, at this stage, dependent on $x$, and (2) asymptotic equivalence is, in
general, not compatible with addition (if $f\sim g$ and $h\sim k$ it is not true
in general that $f + h \sim g + k$ in the same limits). In order to overcome
these two obstacles, a detour is needed.

Before relating $I_1$ to $I_1^\star$, we analyze $I^\star_1$ asymptotically.
Multiply~\eqref{eq:iustar} by $K$ and integrate between $0$ and $x_*$ to obtain
$I^\star_1$. Explicitly,
\begin{multline}\label{eq:B-lbd}
  I_1^\star =
  \Phi(- 3\Delta^\star/2, x_* - 3\Delta^\star/2)
  K^{1 + 2 r^\star - 2\beta^\star}
  -
  \Phi(- \Delta^\star/2, x_* - \Delta^\star/2)
  K^{1 - 2\beta^\star}
  \\
  +
  \Phi(3\Delta^\star/2, x_* + 3\Delta^\star/2)
  K^{1 + 2r^\star - 2\beta^\star}
  -
  \Phi(\Delta^\star/2, x_* + \Delta^\star/2)
  K^{1 - 2\beta^\star}.
\end{multline}
Item~\ref{item:asymp-ccdf-pert} of Lemma~\ref{lem:gaussian-things}
and~(\ref{eq:def-star}) allows us to perform the following asymptotic analysis
of $\Phi$ as if $x_*$ was equal to
$(\beta^\star / 2\sqrt{r^\star})\sqrt{2\ln K}$ because $x_*$ gets close to this
value fast enough (see the discussion around~(\ref{eq:def-star})).
Asymptotically, for $r^\star > \rho(\beta^\star)$, the first term on the right
hand side (\ref{eq:B-lbd}) is of larger order than the rest: (1) the first term
tends to infinity as $K\to\infty$ because it is either
$\Theta(K^{1 + 2r^\star - 2\beta^\star})$ if $r^\star\leq \beta^\star/3$ or
$\tilde{\Theta}(K^{1 - \beta^\star - \frac14\paren{\frac{r^\star -
      \beta^\star}{\sqrt{r^\star}}}^2})$ if $r^\star > \beta^\star/3$; (2) the
second term tends to zero because, as $K\to\infty$, it is either
$\tilde{\Theta}(K^{1 - 2\beta^\star})$ if $\beta^\star \geq r^\star$, or
$o(K^{1 - 2\beta^\star})$ if $\beta^\star < r^\star$; (3) the last last two
terms in the right hand side of (\ref{eq:B-lbd}) tend to zero because they are
both $\tilde{\Theta}( K^{1 - r^\star/4 - 2\beta^\star})$ in the same limit.
Taking care of the case that $r^\star = \rho(\beta^\star)$, these observations
and Lemma~\ref{lem:gaussian-things} (items~\ref{item:asymp-ccdf}--\ref{item:asymp-ccdf-transl}) imply that, as $K\to\infty$,
\begin{equation}\label{eq:i1star-asymp}
  I_1^\star \sim
  \begin{cases}
    K^{1 + 2r^\star - 2\beta^\star}  &\text{ if }
                                       1/2 < \beta^\star < 3/4
                                       \text{ and }
                                       \beta^\star - 1/2 \leq r^\star < \beta^\star / 3;\\
    \frac12 K^{1 + 2r^\star - 2\beta^\star} &\text{ if }
                                       1/2 < \beta^\star \leq 3/4
                                              \text{ and }
                                              r^\star = \beta^\star/3; \\
    \frac{1}{\sqrt{2\pi}}
    \frac{2\sqrt{r^\star}}{3r^\star
    - \beta^\star}\frac{K^{1  - \beta^\star -
    \frac{1}{4}\paren{\frac{r^\star
    - \beta^\star}{\sqrt{r^\star}}}^
    2}}{\sqrt{2\ln K}}
                                     &\text{ if }
                                       1/2 < \beta^\star < 3/4
                                       \text{ and }
                                       r^\star > \beta^\star/3,\\
                                     &\text{ or if }
                                       3/4 \leq \beta^\star \leq 1
                                       \text{ and }
                                       r^\star \geq (1 - \sqrt{1 - \beta^\star})^2.
  \end{cases}
\end{equation}
This proves claim \eqref{eq:linear_rates} about the asymptotic expression for
$I_1^\star$.

In order to relate $I_1$ to $I_1^\star$, we begin by relating~\eqref{eq:iu} to
its well-specified analogue~\eqref{eq:iustar}. To that end, it is necessary to
analyze $m(\xi)$. Recall that $\xi$ depends on $x$ and that, for fixed $\xi$,
$m(\xi) = o(1/\ln K)$ as $K\to \infty$. Since $I_1$ is obtained
from~(\ref{eq:iu}) by integration with respect to $x$ (and multiplication by
$K$), for the purpose of relating $I_1$ to $I^\star_1$, it is useful to show
that the asymptotic statement that $m(\xi) = o(1/\ln K)$ as $K\to \infty$ also
holds uniformly over $0 \leq x \leq x_*$. To this end, define
$\xi_* := \xi(x_*)$. Since $x \mapsto \xi(x)$ is increasing and $x \leq x^*$
(where $x_*$ is as in \eqref{eq:def-star}), then
\begin{equation}\label{eq:bound_mxi}
  m(\xi)
  \leq
  m(\xi_*)
  =
  \frac{2d_r}{r}\frac{x_*}{\sqrt{2\ln K}}
  =
  \frac{2 d_r}{r}\left(\frac{\beta}{2\sqrt{r}} + \frac{d_r}{\sqrt{2r}}\right) .
\end{equation}
This upper bound is independent of $x$ and, consequently,
$\sup_{0\leq \xi \leq \xi_*}m(\xi) = o(1 / \ln K)$ as $K\to \infty$.
Additionally, since $m(\xi) \geq 0$ (because $x\geq 0$),
$\inf_{0\leq \xi \leq \xi_*}m(\xi) = 0$. Hence, the quantities $d_r, d_\beta$,
and $\sup_{0\leq \xi \leq \xi_*}\{\pm m(\xi)\} $ are all $o(1 / \ln K)$ as
$K\to \infty$.

Now, we relate each of the four terms in the right hand side of~\eqref{eq:iu} to
their respective well-specified analogue in~\eqref{eq:iustar}. For the argument,
we focus on the exponential factor in $K$ that accompanies each of the first
four terms in~\eqref{eq:iu}. Factors of the form $K^{-d_\beta}$,
$K^{-(d_\beta + d_r)}$, and $K^{\pm m(\xi)}$ appear. For the first kind, we can
use the assumptions on $d_\beta$ to conclude that
$K^{-d_\beta} = \exp(-d_\beta \ln K) = 1 + o(1)$ as $K\to\infty$; this is
independent of $x$. Similarly, independently of $x$,
$K^{-(d_\beta + d_r)} = 1 + o(1)$ as $K\to\infty$. Now, consider $K^{m(\xi)}$.
Since $x\mapsto m(\xi(x))$ is increasing, then $K^{m(\xi)}\in [1, K^{m(\xi_*)}]$
any time that $x\in [0, x_*]$. Consequently, since also
$m(\xi_*) = o(1 / \ln K)$ as $K\to \infty$ (by the previous paragraph), then
$\sup_{0\leq\xi\leq\xi_*}|K^{m(\xi)} - 1| = o(1)$ in the same limit. Similarly,
$\sup_{0\leq\xi\leq\xi_*}|K^{-m(\xi)} - 1| = o(1)$ as $K\to\infty$. Using these
arguments and the fact that if $f\sim g$ and $h\sim k$ in a certain limit then
$fh\sim gk$ in the same limit, we obtain that, as $K\to\infty$, each of the
first four terms in~\eqref{eq:iu} is asymptotic to its well-specified analogue
in~\eqref{eq:iustar} and that the $1 + o(1)$ factor that relates them can be
upper and lower bounded by two $1 + o(1)$ sequences that are independent of $x$
any time that $0\leq x \leq x_*$. This uniform (in $0\leq x\leq x_*$) asymptotic
expression is useful to relate $I_1$ to $I_1^\star$.

With the arguments so far, we can relate $I_1$ to $I^\star_1$. As a consequence
of the argument in the previous paragraph, multiplying both sides
of~(\ref{eq:iu}) by $K$ and integrating between $0$ and $x_*$ with respect to
$x$, we obtain that $I_1$ equals, up to $1 + o(1)$ factors (as $K\to\infty$),
the sum of each of the terms in $I_1^\star$ as in~(\ref{eq:B-lbd}) plus a
misspecification term; that is, as $K\to\infty$,
 \begin{multline}\label{eq:i1-up-to-o1}
   I_1 =\\
     \Phi(- 3\Delta^\star/2, x_* - 3\Delta^\star/2)
     K^{1 + 2 r^\star - 2\beta^\star}(1 + o(1))
     - \Phi(- \Delta^\star/2, x_* - \Delta^\star/2)
     K^{1 - 2\beta^\star}(1 + o(1))\\
     + \Phi(3\Delta^\star/2, x_* + 3\Delta^\star/2)
     K^{1 + 2r^\star - 2\beta^\star}(1 + o(1))
     - \Phi(\Delta^\star/2, x_* + \Delta^\star/2)
     K^{1 - 2\beta^\star}(1 + o(1))\\
     - \Xi_K,
 \end{multline}
 where
 \begin{equation}\label{eq:xi-def}
 	\Xi_K
 	:=
 	K\int_0^{x_*}
 	\eta_K
 	\rmd x.
 \end{equation}
 The same asymptotic analysis that was used on~(\ref{eq:B-lbd}) to
 obtain~(\ref{eq:i1star-asymp}) can be used on the first four terms
 of~(\ref{eq:i1-up-to-o1}); those four terms are asymptotic to the right hand
 of~(\ref{eq:i1star-asymp}) as $K\to\infty$ and, consequently, to $I_1^\star$
 as well in the same limit. Therefore, as $K\to\infty$,
 \begin{equation*}
   I_1 = I_1^\star(1 + o(1)) - \Xi_K.
 \end{equation*}
This shows~(\ref{eq:I1-and-I1star}).

Now, we turn to showing the last claim of the lemma that is
left,~(\ref{eq:xi-asymp}), which requires analyzing the misspecification
cost~\eqref{eq:xi-def}. Computing~(\ref{eq:xi-def}) in terms
of~(\ref{eq:xi-precursor}) and using both the triangle inequality and that
$0 \leq m(\xi) \leq m(\xi_*)$ from the bound in~\eqref{eq:bound_mxi}, an upper
bound is
\begin{multline*}
  |\Xi_K|
  \leq
  \Phi(- \Delta^\star/ 2, x_* - \Delta^\star/ 2)K^{1 - \beta}
  (1 - \varepsilon^\star)
  |
  1
  -
  K^{- d_r}
  |\\
  +
  \Phi(\Delta^\star/2, \Delta^\star/2 + x_*)
  K^{1 - \beta}
  (1 - \varepsilon^\star)
  |
  1
  -
  K^{- m(\xi_*) - d_r}
  |,
\end{multline*}
Recall $m(\xi_*)$ from \eqref{eq:bound_mxi} and that, for $c\in\mathbf{R}$,
$K^{cd_r} = \rme^{cd_r \ln K} \sim 1 + cd_r\ln K$ as $K\to\infty$ by our
assumption on $d_r$. Then, the previous bound for $|\Xi_K|$ implies that
\begin{multline}
  |\Xi_K|\label{eq:xi-lesssim}
  \lesssim
  \underbrace{\Phi(- \Delta^\star/ 2, x_* -
  \Delta^\star/ 2)K^{1 - \beta^\star}
  |d_r|\ln K}_{=:A_K}
  \\ +
  \underbrace{\Phi(\Delta^\star/2, \Delta^\star/2 + x_*)
  K^{1 - \beta^\star}
  |d_r + m(\xi_*)|\ln K}_{=:B_K},
\end{multline}
where $A_K$ and $B_K$ are defined, respectively, as the first and second term of
the right hand side. We now analyze the right hand side of the last display
asymptotically as $K\to \infty$. We do this in two steps: firstly, we show that
the right hand side of the previous display is asymptotic to $A_K$, the first
term on the right hand side; secondly, we derive an asymptotic upper bound for
$A_K$, which, given the first step, also results in an asymptotic upper bound
for $|\Xi_K|$.

On the one hand, $B_K$, the second term in the right hand side of
\eqref{eq:xi-lesssim}, is $\tilde{\Theta}(K^{1 - r^\star/4 - \beta^\star}d_r)$
as $K\to\infty$. On the other hand, $A_K$ is analyzed in two cases: (i) if
$r^\star \leq \beta^\star$, then $x^\star \geq \Delta^\star / 2$ and the first
term on the right hand side is $\tilde{\Theta}(K^{1 - \beta^\star}d_r)$ and
hence $B_K = o(A_K)$; (ii) if $r^\star > \beta^\star$, then
$x^\star - \Delta^\star/2 < 0$ and the first term of the right hand side is
$\tilde{\Theta}\Big(K^{1 - \beta^\star - \frac14\paren{\frac{r^\star -
      \beta^\star}{\sqrt{r^\star}}}^2}d_r\Big)$. Since
$(r^\star - \beta^\star)^2< r^{\star 2}$ for $r^\star > \beta^\star$, then
$-(1 / 4) r^\star < -(1 / 4)\sqbrack{(r^\star - \beta^\star) /
  \sqrt{r^\star}}^2$ and consequently $B_K = o(A_K)$ as well. Therefore, we
conclude that $B_K=o(A_K)$ always. Finally, using
Lemma~\ref{lem:gaussian-things} (items~\ref{item:asymp-ccdf}
and~\ref{item:asymp-ccdf-transl}) and (\ref{eq:xi-lesssim}), we derive an
asymptotic upper bound for $A_K$, which is also an asymptotic upper bound for
$\Xi_K$ as $K\to\infty$. Indeed, in that limit
\begin{align*}
  |\Xi_K|
  &\lesssim
  \begin{cases}
  K^{1 - \beta^\star}
  |d_r|\ln K
        &\text{ if } r^\star <
          \beta^\star;\\
  \frac12 K^{1 - \beta^\star}
  |d_r|\ln K
    &\text{ if } r^\star =
      \beta^\star;\\
    \frac{1}{\sqrt{2\pi}}
    \frac{2\sqrt{r^\star}}{r^\star
    - \beta^\star}
    \frac{K^{1 - \beta^\star - \frac14\paren{\frac{r^\star - \beta^\star}{\sqrt{r^\star}}}^2}}{\sqrt{2\ln K }}
    |d_r|\ln K
        &\text{ if } r^\star >  \beta^\star.\\
  \end{cases}
\end{align*}
This is the asymptotic upper bound for $|\Xi_K|$ from the claim, concluding the
proof of the lemma.
\end{proof}

\subsubsection{Proof of Lemma~\ref{lem:quadratic-lemma}}
\label{sec:proof-quadratic-lemma}

\begin{proof}[Proof of Lemma \ref{lem:quadratic-lemma}]
  The following arguments are similar in spirit---but different in the
  details---to those used in the proof of Lemma~\ref{lem:linear-lemma}, in
  Section~\ref{sec:proof-linear-lemma}. Let $x = \sqrt{2\xi\ln K}$. In this
  $\xi$ parametrization,
  $u(\Delta^\star \pm x) = K^{\pm 2\sqrt{\xi r^\star} \pm m(\xi) - (\beta +
    d_r)} - K^{-\beta}$ with $m(\xi) = 2\sqrt{\xi}d_r/r$. We first compute
  $i_0[u^2]$ and $i_1[u^2]$, where $i_0$ and $i_1$ are defined
  in~\eqref{eq:i0-def} and~\eqref{eq:i1-def} respectively. We perform this
  computation in the $\xi$ parametrization and return to the $x$ parametrization
  later. Indeed, completing the squares, we obtain that
  \begin{align*}
    \sqrt{2\pi}i_0[u^2]&=K^{-(\sqrt{\xi} - 3\sqrt{r^\star} / 2)^2 + 2r^\star + 2m(\xi) - 2(\beta + d_r)}
    +
    K^{-(\sqrt{r^\star}/2 +  \sqrt{\xi})^2 - 2\beta}\\
    &\quad-2 K^{-(\sqrt{\xi} - \sqrt{r^\star}
      / 2)^2 + m(\xi)- (2\beta +
      d_r)}
    +
        K^{ -(\sqrt{\xi} + 3\sqrt{r^\star} / 2)^2 + 2r^\star - 2m(\xi) - 2(\beta + d_r)}\\
    &\quad+
    K^{-(\sqrt{r^\star}/2 -  \sqrt{\xi})^2 - 2\beta}
    -
    2
    K^{-(\sqrt{\xi} + \sqrt{r^\star}
      / 2)^2 - m(\xi)- (2\beta +
      d_r)},
  \end{align*}
  and that
  \begin{align*}
    \sqrt{2\pi}i_1[u^2]
    &=
    K^{-(\sqrt{\xi} - 5 \sqrt{r^\star} / 2)^2 + 6 r^\star + 2m(\xi) - 2(\beta + d_r)}
    +
    K^{-(\sqrt{r^\star}/2 -  \sqrt{\xi})^2 - 2\beta}\\
    &\quad-2 K^{ -(\sqrt{\xi} - 3 \sqrt{r^\star} / 2)^2 + 2r^\star + m(\xi)- (2\beta + d_r)}
    +
        K^{-(\sqrt{\xi} + 5\sqrt{r^\star} / 2)^2 + 6r^\star - 2m(\xi) - 2(\beta + d_r)}\\
        &\quad
    +
    K^{-(\sqrt{r^\star}/2 +  \sqrt{\xi})^2 -\beta}
    -
    2
    K^{-(\sqrt{\xi} + 3\sqrt{r^\star} / 2)^2 + 2r^\star - m(\xi)- (2\beta + d_r)}.
  \end{align*}
Now we compute $i[u^2]$, which, recall from Lemma~\ref{lem:translate-center},
  equals the convex combination
  $\varepsilon^\star i_1[u^2] + (1 - \varepsilon^\star)i_0[u^2]$. Using the
  previous expressions for $i_1[u^2]$ and $i_0[u^2]$, the target $i[u^2]$ is
  computed in terms of the Gaussian density $\varphi$ and the original $x$
  parametrization. For the subsequent analysis, we rearrange all terms as follows:
  \begin{multline*}
    i[u^2](x) =\\
    \begin{aligned}
      &\varphi(x - 5\Delta^\star/2)
      \bracks{
      K^{6 r^\star - 3\beta^\star + 2m(\xi) - 2(d_\beta + d_r)}
      }
      \\
      +
      &\varphi(x - 3\Delta^\star/2)
      \bracks{
      (1 - \varepsilon^\star)
      K^{2r^\star - 2\beta^\star + 2m(\xi) - 2(d_\beta + d_r)}
      -
      2
      K^{2 r^\star - 3\beta^\star + m(\xi)- (2d_\beta + d_r)}
      }
      \\
      +
      &\varphi(x - \Delta^\star/2)
      \bracks{
      (1 - \varepsilon^\star)
      K^{- 2\beta^\star - 2d_\beta}
      -
      2
      (1 - \varepsilon^\star)
      K^{- 2\beta^\star + m(\xi)- (2d_\beta + d_r)}
      +
      K^{-3\beta^\star - 2d_\beta}
      }
      \\
      +
      &\varphi(x + \Delta^\star/2)
      \bracks{
      K^{-2\beta^\star - d_\beta}
      +
      (1 - \varepsilon^\star)
      K^{- 2\beta^\star - 2d_\beta}
      -
      2
      (1 - \varepsilon^\star)
      K^{- 2\beta^\star - m(\xi)- (2d_\beta + d_r)}
      }
      \\
      +
      &\varphi(x + 3\Delta^\star/2)
      \bracks{
      (1 - \varepsilon^\star)
      K^{2r^\star - 2\beta^\star - 2m(\xi) - 2(d_\beta + d_r)}
      -
      2
      K^{2r^\star - 3\beta^\star - m(\xi)- (2d_\beta + d_r)}
      }
      \\
      +
      &\varphi(x + 5\Delta^\star/2)
      \bracks{
      K^{6r^\star - 3\beta^\star - 2m(\xi) - 2(d_\beta + d_r)}
      }.
    \end{aligned}
  \end{multline*}
  For the integral defining $I_2$, we need to evaluate the expressions above at
  points $x\leq x_*$; consequently, within the $\xi$ reparametrization at
  $0\leq m(\xi) \leq m(\xi_*)$. Using this, note that $K$ times the integral of
  the fourth term in the right hand side of the previous display is
  $o(K^{1 - 2\beta^\star})$ as $K \to \infty$. It will become apparent that this
  is of lower order than the other terms.

  Integrate and multiply the previous display by $K$, use that
  $(1 - \varepsilon^\star)\to 1$, use the observation on its fourth term, and
  use the same argument as the one employed to derive~(\ref{eq:i1-up-to-o1})
  from~(\ref{eq:iu}) to obtain that, as
  $K\to\infty$,
  \begin{equation}\label{eq:zsq_eachterm}
    \begin{aligned}
      I_2^\star
      = & \
             \Phi( - 5\Delta^\star/2, x_* - 5\Delta^\star/2)
          K^{1 + 6 r^\star - 3\beta^\star}
          (1 + o(1))
      \\
           &+
             \Phi( - 3\Delta^\star/2, x_* - 3\Delta^\star/2)
             K^{1 + 2r^\star - 2\beta^\star}
             (1 + o(1))
      \\
           &-
             \Phi(- \Delta^\star/2, x_* - \Delta^\star/2)
             K^{1 - 2\beta^\star}
             (1 + o(1))
                   \\
           &+ o(K^{1 - 2\beta^\star})
      \\
           &+
             \Phi(3\Delta^\star/2, x_* + 3\Delta^\star/2)
             K^{1 + 2r^\star - 2\beta^\star}
             (1 + o(1))
      \\
           &+
             \Phi(5\Delta^\star/2, x_* + 5\Delta^\star/2)
             K^{1 + 6r^\star - 3\beta^\star}          (1 + o(1)).
    \end{aligned}
  \end{equation}

We now proceed to analyze the remaining five terms in (\ref{eq:zsq_eachterm}).
We first claim that the third, fifth and sixth terms are all
$O(K^{1 - 2\beta^\star})$ as $K\to\infty$---and hence $o(1)$ in the same limit.
Indeed, the third term is $O(K^{1 - 2\beta^\star})$ because
$\Phi(- \Delta^\star/2, x_* - \Delta^\star/2) \leq 1$. For the fifth term, we
have that
  \begin{align*}
  \Phi(3\Delta^\star/2, x_* + 3\Delta^\star/2)K^{1+2r^\star-2\beta^\star} &\leq \overline\Phi(3\Delta^\star/2)K^{1+2r^\star-2\beta^\star} \\
  &= O(K^{-\frac{9r}{4}})K^{1+2r^\star-2\beta^\star} ,
  \end{align*}
  where the asymptotic order of $\overline\Phi$ follows from the first statement
  of Lemma~\ref{lem:gaussian-things}. Working out the exponent in the last
  equality implies that the fifth term in~\eqref{eq:zsq_eachterm} is
  $o(K^{1-2\beta^\star})$ as $K\to\infty$. For the sixth term, an analogous
  argument shows that it is $o(K^{1-2\beta^\star})$ as $K\to\infty$. Hence, as
  was claimed, the last four terms of \eqref{eq:zsq_eachterm} are
  $O(K^{1 - 2\beta^\star})$ as $K\to\infty$, and they are consequently $o(1)$ as
  well in the same limit.

  The remainder of the proof is dedicated to characterize the first two terms
  of~\eqref{eq:zsq_eachterm}; we will see that the right hand side
  of~\eqref{eq:zsq_eachterm} is asymptotic to those two terms. Using
  items~\ref{item:asymp-ccdf}--\ref{item:asymp-ccdf-transl} from
  Lemma~\ref{lem:gaussian-things} and recalling the definition of $x_*$ from
  \eqref{eq:def-star}, we find that
  \begin{multline}\label{eq:cincolukas}
      \Phi( - 5\Delta^\star/2, x_* - 5\Delta^\star/2) K^{1 + 6 r^\star -
        3\beta^\star} \\
      \sim K^{1 + 6 r^\star - 3\beta^\star}
      \begin{cases}
        1/2 &\text{if } \beta^\star = 5r^\star;\\
        1   &\text{if } r^\star < \beta^\star/5; \\
        \frac{1}{\sqrt{2\pi}}
        \frac{2\sqrt{r^\star}}{5r^\star
        - \beta^\star}
        \frac{K^{-\frac{1}{4}\paren{\frac{5r^\star
        -
        \beta^\star}{\sqrt{r^\star}}}^2}}{\sqrt{2\ln
        K}}    &\text{if }
                 r^\star > \beta^\star/5;
      \end{cases}
    \end{multline}
    and that
  \begin{multline} \label{eq:treslukas}
    \Phi( - 3\Delta^\star/2, x_* - 3\Delta^\star/2)
    K^{1 + 2r^\star - 2\beta^\star}
     \\
      \sim K^{1 + 2r^\star - 2\beta^\star}
      \begin{cases}
        1/2 &\text{if } \beta^\star = 3r^\star;\\
        1   &\text{if } r^\star < \beta^\star/3;\\
        \frac{1}{\sqrt{2\pi}}
        \frac{2\sqrt{r^\star}}{3r^\star
        - \beta^\star}
        \frac{K^{-\frac{1}{4}\paren{\frac{3r^\star
        -
        \beta^\star}{\sqrt{r^\star}}}^2}}{\sqrt{2\ln
        K}}    &\text{if }
                 r^\star > \beta^\star/3.
      \end{cases}
  \end{multline}
  Notice that that the exponents involved in the respective last cases of
  (\ref{eq:cincolukas}) and (\ref{eq:treslukas}) are actually equal; they
  satisfy
  \begin{align*}
    1 + 6r^\star - 3\beta^\star -
    \frac{1}{4}\paren{\frac{5r^\star -
    \beta^\star}{\sqrt{r^\star}}}^2
    &=
    1 + 2r^\star - 2\beta^\star -
    \frac{1}{4}\paren{\frac{3r^\star -
    \beta^\star}{\sqrt{r^\star}}}^2 \\
    &=
    1 - \beta^\star -
    \frac14\paren{\frac{r^\star -
    \beta^\star}{\sqrt{r^\star}}}^2.
  \end{align*}

  To finish the asymptotic characterization of $I_2$, we need to consider six
  cases to cover the full area of interest in the $\beta^\star$-$r^\star$~plane.
  These cases correspond to all possible combinations of the asymptotic regimes
  from~(\ref{eq:cincolukas}) and~(\ref{eq:treslukas}).
  \begin{description}
  \item[Case 1.] If
    $1/2 < \beta^\star \leq 5/8 $
    and
    $\beta - 1/2 \leq r^\star <
    \beta^\star/5$, then
    \begin{equation*}
      K\int_0^{x_*}i[u^2]\rmd x
      \sim
      K^{1 + 6r^\star - 3\beta^\star}
      +
      K^{1 + 2r^\star - 2\beta^\star}
      \sim
      K^{1 + 2r^\star - 2\beta^\star}
    \end{equation*}
    because
    $r^\star < \beta^\star / 5 <
    \beta^\star / 4$, and
    $6r^\star - 3\beta^\star <
    2r^\star - 2\beta^\star$ if an
    only if
    $r^\star < \beta^\star / 4$. This same argument is used in the next case.
  \item[Case 2.] If
    $1/2 < \beta^\star \leq 5/8 $
    and $r^\star = \beta^\star/5$,
    then
    \begin{equation*}
      K\int_0^{x_*}i[u^2]\rmd x
      \sim
      \frac{1}{2}K^{1 + 6r^\star - 3r^\star}
      +
      K^{1 + 2r^\star - 2\beta^\star}
      \sim
      K^{1 + 2r^\star - 2\beta^\star},
    \end{equation*}
    where, since the order of each term is the same as in the previous case,
    the same reasoning can be used.
  \item[Case 3.] If
    $5/8 < \beta^\star < 3/4$ and
    $r^\star = \beta^\star - 1/2 $,
    then
    \begin{align*}
      K\int_0^{x_*}i[u^2]\rmd x
      &\sim
         \frac{1}{\sqrt{2\pi}}
         \frac{2\sqrt{r^\star}}{5r^\star
         - \beta^\star}
         \frac{K^{1 + 6r^\star -
         3\beta^\star - \frac{1}{4}\paren{\frac{5r^\star
         -
         \beta^\star}{\sqrt{r^\star}}}^2}}{\sqrt{2\ln
         K}}
         +
         K^{1 + 2r^\star - 2\beta^\star}
      \\
       &\sim
         K^{1 + 2r^\star -
         2\beta^\star}
    \end{align*}
    because taking the difference between the exponents of the
    two terms above we find that
    $$\sqbrack{6r^\star - 3\beta^\star - \frac{1}{4}\paren{\frac{5r^\star -
          \beta^\star}{\sqrt{r^\star}}}^2} - [2r^\star - 2\beta^\star] =
    -\frac14\paren{\frac{3r^\star - \beta^\star}{\sqrt{r^\star}}}^2 < 0. $$
    Similarly as happened in Case 1 and Case 2, this reasoning can be used in
    the next two cases.
  \item[Case 4.] When $1/2 < \beta^\star \leq 3/4$ and
    $\max\{\beta^\star/5, \beta^\star - 1/2\} < r < \beta^\star/3$, then
    \begin{align*}
      K\int_0^{x_*}i[u^2]\rmd x
      &\sim
        \frac{1}{\sqrt{2\pi}}
        \frac{2\sqrt{r^\star}}{5r^\star
        - \beta^\star}
        \frac{K^{1 + 6r^\star -3\beta^\star -\frac{1}{4}\paren{\frac{5r^\star
        -
        \beta^\star}{\sqrt{r^\star}}}^2}}{\sqrt{2\ln
        K}}
        +
        K^{1 + 2r^\star - 2\beta^\star}
      \\
      &\sim
        K^{1 + 2r^\star - 2\beta^\star},
    \end{align*}
    where we repeated the same argument as in the previous case.
  \item[Case 5.] If
    $1/2 < \beta^\star \leq 3/4$ and
    $r^\star = \beta^\star / 3$, then
    \begin{align*}
      K\int_0^{x_*}i[u^2]\rmd x
      &\sim
        \frac{1}{\sqrt{2\pi}}
        \frac{2\sqrt{r^\star}}{5r^\star
        - \beta^\star}
        \frac{K^{1 + 6r^\star -
        3\beta^\star -\frac{1}{4}\paren{\frac{5r^\star
        -
        \beta^\star}{\sqrt{r^\star}}}^2}}{\sqrt{2\ln
        K}}
        +
        \frac12 K^{1 + 2r^\star - 2\beta^\star}
      \\
      &\sim
        \frac12 K^{1 + 2r^\star - 2\beta^\star}
    \end{align*}
    where we used again the argument from Case~3.
  \item[Case 6.] If
    $1/2 < \beta^\star \leq 3/4$ and
    $r^\star > \beta^\star / 3$, or
    if $3/4 < \beta^\star \leq 1$
    and
    $r^\star \geq (1 - \sqrt{1 -
      \beta^\star})^2$, then
    \begin{align*}
      K\int_0^{x_*}i[u^2]\rmd x
      &\sim
        \frac{2\sqrt{r^\star}}{\sqrt{2\pi}}\paren{
        \frac{1}{3r^\star - \beta^\star}
        +
        \frac{1}{5r^\star - \beta^\star}
        }
        \frac{
        K^{1 - \beta^\star -
        \frac14\paren{\frac{\beta^\star
        - r^\star}{\sqrt{r^\star}}}^2}
        }{\sqrt{2\ln K}}.
    \end{align*}
  \end{description}
  Using the result of the six cases, (\ref{eq:quadratic-rates}) follows; it is a
  concise summary of the previous cases. This is what was to be shown.
\end{proof}

\subsection{Gaussian integrals}
\label{sec:technical-lemmas}

Recall from Section~\ref{sec:notation-proofs} that
$\varphi(x) = (1 / \sqrt{2\pi})\exp(-x^2 / 2)$ is the density of a standard
Gaussian distribution, that $\Phi(x) = \int_{-\infty}^x\varphi(w)\rmd w$ is its
cumulative distribution function, that $\overline{\Phi}(x) = 1 - \Phi(x)$ is its
complementary cumulative distribution function, and that
$\Phi(a, b) = \Phi(b) - \Phi(a)$ is the Gaussian measure of $[a,b]$ for
$a\leq b$. The next lemma gives asymptotic expressions for Gaussian integrals
that are used in the preceding sections.
\begin{lemma}\label{lem:gaussian-things}
  Let $\xi > 0$ and let $x = \sqrt{2\xi\ln K}$.
  \begin{enumerate}
  \item \label{item:asymp-ccdf} As $K\to\infty$,
    $$\overline{\Phi}(x) \sim \frac{1}{\sqrt{2\pi}} \frac{K^{-\xi}}{\sqrt{2\xi
        \ln K}}.$$
  \item \label{item:asymp-ccdf-pert} If $\zeta_K = o((\ln K)^{-1/2})$ as
    $K\to\infty$, then
    $\overline{\Phi}(x + \zeta_K) \sim \overline{\Phi}(x)$ in the same
    limit.
  \item \label{item:asymp-ccdf-transl} Let $\Delta = \sqrt{2r\ln K}$. Then, as
    $K\to \infty$, both $\Phi(\Delta, x + \Delta)\sim \overline{\Phi}(\Delta)$ and
    \begin{equation*}
      \Phi(-\Delta, x - \Delta)
      \sim
      \begin{cases}
        \frac{1}{\sqrt{2\pi}}\frac{K^{-(r^{1/2} - \xi^{1/2})^2}}{\sqrt{2(r^{1/2} - \xi^{1/2})^2\ln K}}
        &\text{ if } \xi < r;\\
        1 / 2 &\text{ if } \xi = r;\\
        1 &\text{ if } \xi > r.
      \end{cases}
    \end{equation*}
  \item \label{item:asymp-trunc-m1}For each $a\in\mathbf{R}$, we have that
    $\int_0^\infty w\varphi(w - a)\rmd w = \varphi(a) +
    a\Phi(a)$ and consequently, if $a\to\infty$,
    \begin{equation*}
      \int_0^\infty w\varphi(w - a)\rmd w
      \sim
      a
      \text{ and }
      \int_0^\infty w\varphi(w + a)\rmd w
      = o(1).
    \end{equation*}
  \item \label{item:trunc-m2} For each $a\in \mathbf{R}$,
    $\int_0^\infty w^2\varphi(w-a)\rmd w = (a^2 + 1)\Phi(a) +
    \frac{a}{\sqrt{2\pi}}\rme^{-a^2 / 2}$. Consequently, if $a\to\infty$
    \begin{align*}
      \int_0^\infty w^2\varphi(w - a)\rmd w
      \sim
      a^2
      \text{ and }
      \int_0^\infty w^2\varphi(w + a)\rmd w
      \sim
      \overline{\Phi}(a).
    \end{align*}
  \end{enumerate}
\end{lemma}
\begin{proof}[Proof of Lemma~\ref{lem:gaussian-things}]
  Item~\ref{item:asymp-ccdf} is known
  \citep[][Formula~7.1.13]{abramowitz_handbook_1965}. In our notation,
  Formula~7.1.13 \citep{abramowitz_handbook_1965} is equivalent to
  \begin{equation}\label{eqn:abramowitz}
    \sqrt{\frac{2}{\pi}}\frac{\rme^{ - x ^ 2 / 2}}{x + \sqrt{x ^2 + 4}}
    <
    \overline{\Phi}(x)
    <
    \sqrt{\frac{2}{\pi}}\frac{\rme^{ - x ^ 2 / 2}}{x + \sqrt{x ^2 + 8/\pi}}.
  \end{equation}
  Replacing $x = \sqrt{2\xi\ln K}$ and taking $K\to\infty$ shows the first
  claim.
  To prove Item~\ref{item:asymp-ccdf-pert} use again the relation in \eqref{eqn:abramowitz} to see that
  \begin{equation*}
    \sqrt{\frac{2}{\pi}}\frac{K^{-\xi}\rme^{-\zeta_K^2 / 2
        -
        \sqrt{2\xi}\zeta_K
      }}{x + \zeta_K + \sqrt{(x + \zeta_K)^2 + 4}}
    <
    \overline{\Phi}(x+\zeta_K)
    <
    \sqrt{\frac{2}{\pi}}\frac{K^{-\xi}\rme^{-\zeta_K^2 / 2
        -
        \sqrt{2\xi}\zeta_K
      }}{x + \zeta_K + \sqrt{(x + \zeta_K)^2 + 8/\pi}}.
  \end{equation*}
  Note that $\rme^{-\zeta_K^2 / 2 + \sqrt{2\xi\ln K}\zeta_K }\to 1$ as
  $K\to \infty$; here the assumption that $\zeta_K = o(( \ln K)^{-1/2})$ in the
  same limit is crucial. Additionally, since $x + \zeta_K \to \infty$ as
  $K\to\infty$ the denominators in the upper and lower bounds in the previous
  equation are both asymptotic to $2x$ as $K\to\infty$, proving the desired
  result.
  The first claim of Item~\ref{item:asymp-ccdf-transl} follows because
  $\Phi(\Delta, x + \Delta) = \overline{\Phi}(\Delta) - \overline{\Phi}(x +
  \Delta)$ and, by Item~\ref{item:asymp-ccdf},
  $\overline{\Phi}(x + \Delta) = o(\overline{\Phi}(\Delta))$ as $K\to\infty$.
  The second claim of Item~\ref{item:asymp-ccdf-transl} follows from analyzing
  the difference
  \begin{equation*}
    \Phi(-\Delta, x - \Delta)
    =
    \Phi(x - \Delta)
    -
    \Phi(-\Delta)
  \end{equation*}
  By Item~\ref{item:asymp-ccdf}, as $K\to\infty$,
  $\Phi(-\Delta) \sim K^{-r} / (\sqrt{2\pi}\sqrt{2r\ln K})$. To analyze
  $\Phi(x - \Delta)$, three cases are needed; in all cases,
  $\Phi(-\Delta, x - \Delta)$ is asymptotic to $\Phi(x - \Delta)$ as
  $K\to\infty$. Indeed, if $\xi < r$, then, by Item~\ref{item:asymp-ccdf},
  $\Phi(-\Delta, x - \Delta) \sim (1 / \sqrt{2\pi})[K^{-(r^{1/2} -
    \xi^{1/2})^2}/ \sqrt{2(r^{1/2} - \xi^{1/2})^2\ln K}]$ as $K\to\infty$; if
  $r = \xi$, then $\Phi(-\Delta, x - \Delta) \sim 1/2$; and if $\xi > r$, then
  $\Phi(-\Delta, x - \Delta) \sim 1$.
  The first claim in Item~\ref{item:asymp-trunc-m1} follows from the fact the
  integral can be computed analytically yielding
  $\int_0^\infty w\varphi(w - a)\rmd w = \varphi(a) + a\Phi(a)$. The second
  claim in Item~\ref{item:asymp-trunc-m1} follows from observing that
  $\varphi(a) + a\Phi(a)\sim a$ as $a\to\infty$; the third, from observing that
  $\varphi(a) - a\Phi(-a) = \varphi(a) - a\overline{\Phi}(a) \to 0$ as
  $a\to\infty$.
  To prove the first claim of Item~\ref{item:trunc-m2}, use integration by
  parts. The first consequence in Item~\ref{item:trunc-m2} follows from the fact
  that $(a^2 + 1)\Phi(a) + (a / \sqrt{2\pi})\exp(-a^2 / 2) \sim a^2$ as
  $a\to\infty$; the second consequence follows from using that
  $\Phi(-a) = \overline{\Phi}(a) \sim \varphi(a) / a$ as $a\to \infty$ and that
  consequently
  $(a^2 + 1)\Phi(-a) - (a / \sqrt{2\pi})\exp(-a^2 / 2) \sim \overline{\Phi}(a)$
  in the same limit.
\end{proof}

\subsection{Proof of Theorem~\ref{thm:prior-stops-on-time}}
\label{sec:proof-prior-stops}

\begin{proof}[Proof of Theorem~\ref{thm:prior-stops-on-time}]
For the case where $t < t^\star$, an analogous reasoning to that at the start of the proof of Theorem~\ref{thm:stop-on-time} applies; namely, comparing the power of the sequential test with that of the most powerful test at sample size $t$ implies the result.

The rest of the proof pertains exclusively the case $t > t^\star$. Note that
\begin{equation*}
  \mathbf{P}^\star\{\tau^\star \leq t\}
  =
  \mathbf{P}^\star\{\max_{s\leq t}E_s(\mathbf{\Pi}) \geq 1/\alpha\}
  \geq
  \mathbf{P}^\star\{E_t(\mathbf{\Pi}) \geq 1/\alpha\}.
\end{equation*}
Consequently, for any $(\varepsilon,\delta)\in {\cal G}_C$,
\begin{equation}\label{eq:surv-ineq-1}
  \mathbf{P}^\star\{\tau^\star \leq t\}
  \geq
  \mathbf{P}^\star\{E_t(\varepsilon,\delta) \geq  1/(\alpha \pi_K)\}
  =
  \mathbf{P}^\star\{\ln E_t(\varepsilon,\delta) \geq  \ln(1/(\alpha \pi_K))\}
\end{equation}
with $\pi_K = 1 / |{\cal G}_C|$ because
$E_t(\mathbf{\Pi}) \geq E_t(\varepsilon,\delta) / |{\cal G}_C|$. Now, in the
same way as in the proof of Theorem~\ref{thm:prior-works} from
Section~\ref{sec:proofs-prior}, we can find a pair
$(\varepsilon^\circ, \delta^\circ)\in {\cal G}_C$ that is close enough to
$(\varepsilon^\star, \delta^\star)$ for which the rates~(\ref{eq:master-rates}) of Proposition~\ref{prop:master-rates} apply.
Indeed, by the construction of ${\cal G}_C$ , we have that
$\varepsilon^\circ = K^{-\beta^\circ_K}$ and that
$\delta^\circ = \sqrt{2(1 / T^\circ_K)\ln K}$ for a pair
$(\beta^\circ_K, T^\circ_K )$, such that, as $K\to \infty$,
$d_{\beta, K} := \beta^\circ_K - \beta^\star = o((\ln K)^{-1})$ and
$d_{T, K} := 1 - \sqrt{T^\circ_K / T^\star} = O(K^{-(1 - \beta^\star)})$ as
$K\to\infty$. Next, we show that $\ln E_t(\varepsilon^\circ, \delta^\circ)$ is
sufficiently concentrated around its expectation; that is,
$\ln E_t(\varepsilon^\circ, \delta^\circ)$ is concentrated in a band
$\mathbf{E}^\star[\ln E_t(\varepsilon^\circ, \delta^\circ)]\pm c_K$ with two
properties: (1) $K\mapsto c_K$ is a positive sequence such that the lower bound
$b_K = \mathbf{E}^\star[\ln E_t(\varepsilon^\circ, \delta^\circ)] - c_K$ is also
positive; and (2) $K\mapsto b_K$ grows faster than $\ln(1/(\alpha\pi_K))$. This
is made precise in the following lemma, where the pair
$(T_K^\circ, \beta_K^\circ)$ will take the place of $(T_K, \beta_K)$. Its proof
is found in Section~\ref{sec:proof-power-crux}.
\begin{lemma}\label{lem:power-crux}
  Let $\varepsilon = \varepsilon(K) = K^{-\beta_K}$ and
  $\delta = \delta(K) = \sqrt{2(1 / T_K)\ln K}$, and let both
  $d_{\beta} = d_{\beta}(K) = \beta_K - \beta^\star$ and
  $d_T = d_T(K) = 1 - \sqrt{T_K/T^\star}$. Assume that, as $K\to\infty$, both
  $d_{\beta} = o((\ln K)^{-1})$ and $d_T = O(K^{-(1 - \beta^\star)})$. Assume
  that $\ln (1 / (\alpha \pi_K))\to\infty$ and
  $\ln(1 / (\alpha \pi_K)) = O(\ln K)$.
  Then, provided $t > t^\star$, there exists a sequence $K\mapsto c_K$ of
  positive real numbers such that, for
  $b_K := \mathbf{E}^\star[\ln E_t(\varepsilon, \delta)] - c_K$, we have that
  \begin{enumerate}
  \item \label{item:ck-is-fast}the sequence $K \mapsto b_K$ grows faster than
    $K\mapsto \ln(1 / (\alpha \pi_K))$ in the sense that, as $K\to \infty$,
    $\frac{b_K}{\ln(1 / (\alpha \pi_K))}\to \infty$; and
  \item \label{ck-is-good}
    $\mathbf{P}^\star\{\ln E_t(\varepsilon,\delta) \geq b_K\} \to 1$, also as
    $K\to\infty$.
  \end{enumerate}
\end{lemma}

With this lemma at hand and the discussion before it, we can complete the proof
of Theorem~\ref{thm:prior-stops-on-time} for the case that $t>t^\star$. Indeed,
let $b_K = \mathbf{E}^\star[\ln E_t(\varepsilon^\circ, \delta^\circ)] - c_K$
with the sequence $c_K$ guaranteed by Lemma~\ref{lem:power-crux}.
Then,~\eqref{eq:surv-ineq-1} and the fact that
$\ln(1/(\alpha \pi_K)) \leq \max\{\ln(1/(\alpha \pi_K)), \ b_K\}$ together imply
that
  \begin{equation*}
    \mathbf{P}^\star\{\tau^\star \leq t\}
    \geq
    \mathbf{P}^\star\left\{
    \ln E_t(\varepsilon^\circ,\delta^\circ)
    \geq
    \max\{\ln(1/(\alpha \pi_K)), \ b_K\}\right\}.
  \end{equation*}
  Recall from Section~\ref{sec:prior} that $\ln (1 / (\alpha \pi_K))\to\infty$
  and $\ln(1 / (\alpha \pi_K)) = O(\ln K)$. Thus, Item~\ref{item:ck-is-fast} in
  Lemma~\ref{lem:power-crux} shows that $\max\{1/(\alpha \pi_K), \ b_K\} = b_K$
  for $K$ large enough, and Item~\ref{ck-is-good} shows that consequently, as
  $K\to\infty$, the right hand side of the previous display tends to one. Hence,
  in the same limit
  \begin{equation*}
    \mathbf{P}^\star\{\tau^\star \leq t\} \to 1,
  \end{equation*}
  as it was to be shown.
\end{proof}

\subsubsection{Proof of Lemma~\ref{lem:power-crux}}
\label{sec:proof-power-crux}

\begin{proof}[Proof of Lemma~\ref{lem:power-crux}]
To prove Lemma~\ref{lem:power-crux}, Chebyshev's inequality is used. Indeed, if
$\sigma_K^2 := \var^\star(\ln E_t(\varepsilon, \delta))$, then Chebyshev's
inequality implies that
\begin{equation}\label{eq:power-crux-cheby}
  \mathbf{P}^\star\{\ln E_t(\varepsilon, \delta) < b_K\}
  =
  \mathbf{P}^\star\{\ln E_t(\varepsilon, \delta) - \mathbf{E}^\star[\ln
  E_t(\varepsilon, \delta)] < - c_K\}
  \leq
  \frac{\sigma_K^2}{c_K^2} .
\end{equation}
For this reason, the following arguments concern the standard deviation
$\sigma_K$. Recall from~(\ref{eq:loglikelihood}) that
$\ell(x) = \ln(\varepsilon + (1-\varepsilon)\rme^{\Delta x - \Delta^2 / 2})$,
with $\Delta = \sqrt{2(t / T)\ln K}$, is the single-outcome log-likelihood
ratio. Let $Z_{i, t} = (1 / \sqrt{t})\sum_{s\leq t}X_{i, t}$ be the standardized
sum for each stream $i\in [K]$. Since time is fixed and the distribution of
$Z_{i, t}$ does not depend on the stream $i$, we drop the indexes (we write $Z$
instead of $Z_{i,t}$) to simplify the notation. Using independence and that the
variance is bounded by the second moment,
$\sigma_K^2 = K\var^\star(\ell(Z))\leq K\mathbf{E}^\star[\ell^2(Z)]$. Continuing
from the bound~\eqref{eq:power-crux-cheby}, we find that
\begin{equation*}
  \mathbf{P}^\star\{\ln E_t(\varepsilon^\circ, \delta^\circ) < b_K\}
  \leq
  \frac{K  m_{2,K}}{c_K^2},
\end{equation*}
where $m_{2, K} := \mathbf{E}^\star[\ell^2(Z)]$ is the second moment of the
single-outcome likelihood ratio. It is enough to obtain asymptotic upper bounds
for $m_{2, K}$ to control the rate at which $m_{2,K}$ grows. The next lemma
provides such bounds (its proof is contained in Section~\ref{sec:proof-lem-m2k-rates}).
\begin{lemma}\label{lem:m2k-rates}
  For a fixed $t\geq 1$, let $m_{2,K} = \mathbf{E}^\star[\ell^2(Z)]$ be the
  second moment of the single-outcome likelihood ratio. Let
  $r^\star = t / T^\star$. Then, under the assumptions of
  Lemma~\ref{lem:power-crux},
  \begin{equation}
    \label{eq:m2k-rates}
    K m_{2,K}
    =
    \begin{cases}
      O(K^{1 + 2r^\star - 2\beta^\star})  &\text{ if }
                                            1/2 < \beta^\star < 3/4,
      \\
                                          & \text{ and } \beta^\star - 1/2
                                            < r^\star \leq \beta^\star / 3;
      \\
      \tilde{O}(
      K^{1  - \beta^\star -
      \frac{1}{4}\paren{\frac{r^\star
      - \beta^\star}{\sqrt{r^\star}}}^
      2}
      )
                                          &\text{ if }
                                            1/2 < \beta^\star \leq 1,
      \\
                                          &\text{ and }
                                            \max\{\beta^\star/3,
                                            (1 - \sqrt{1 - \beta^\star})^2\}
                                            < r^\star < \beta^\star;
      \\
      \tilde{O}(K^{1 - \beta^\star})
                                          &\text{ if }
                                            1/2 < \beta^\star \leq 1,
      \\
                                          &\text{ and } r^\star \geq \beta^\star.
    \end{cases}
  \end{equation}
\end{lemma}
Comparing the results of the previous lemma, Lemma~\ref{lem:m2k-rates}, with
Lemma~\ref{prop:master-rates}, we observe that the exponents in
\eqref{eq:m2k-rates} coincide with those in~\eqref{eq:master-rates}.
Specifically, these lemmas give a function $f = f(r^\star, \beta^\star)$ such that
that $f$ is independent of $K$, it is positive if $r^\star > \rho(\beta^\star)$,
and $f$ quantifies the growth exponent of our quantities of interest; indeed,
\begin{equation*}
  \mathbf{E}^\star[\ln E_t(\mathbf{\Pi})] = \tilde{\Omega}(K^{f(r^\star,
    \beta^\star)})
  \quad\text{ and }\quad
  K m_{2,K} = \tilde{O}(K^{f(r^\star, \beta^\star)}) .
\end{equation*}
With this in hand, there is some room to choose $c_K$. For instance, any choice of
$c_K$ such that $c_K = \tilde{\Theta}(K^{\frac34 f(r^\star, \beta^\star) })$
satisfies both
$b_K = \mathbf{E}^\star[\ln E_t(\mathbf{\Pi})] - c_K =
\tilde{\Omega}(K^{f(r^\star, \beta^\star)})$ and $Km_{2,K}/c_K^2 \to 0$ as
$K\to\infty$. This implies both claims of the lemma.
\end{proof}

\subsubsection{Proof of Lemma~\ref{lem:m2k-rates}}\label{sec:proof-lem-m2k-rates}
\begin{proof}[Proof of Lemma~\ref{lem:m2k-rates}]
  Recall the definition of the functional $i$ in \eqref{eq:integrand}. In the
  same vein as the proof of Proposition~\ref{prop:master-rates}, we first use
  Lemma~\ref{lem:translate-center} to derive that
  \begin{align*}
    m_{2,K} &= \int_0^\infty i[\ell^2](x)\rmd x \\
    &= \int_0^{x_*} i[\ell^2](x)\rmd x +
      \int_{x_*}^\infty i[\ell^2](x)\rmd x,
  \end{align*}
  where $x_*$ is given by~(\ref{eq:def-star}). The remainder of this proof is
  devoted to bounding the last two integrals. In proving
  Proposition~\ref{prop:master-rates}, similar quantities were analyzed, but these involved $\ell$ rather than $\ell^2$. Despite this difference, some of the
  technical lemmas and results obtained there are also useful for this proof.

  \paragraph*{Bound for $\int_0^{x_*} i[\ell^2](x)\rmd x $} Similarly to the
  proof of Lemma~\ref{lem:asymp-before-x0}, where a lower bound for
  $z\mapsto\ln(1 + z)$ is used, we bound the integral using a lower bound for
  $z\mapsto \ln^2(1 + z)$. We use the following bound, obtained through a Taylor
  expansion and bounding higher-order terms, which holds for $z > -1$:
  \[
  \ln^2(1 + z) \leq z^2 - \frac{\min\{0, z^3\}}{1 + \min\{0, z\}} .
  \]
  Recall that $u(x) = \varepsilon(\rme^{\Delta x - \Delta^2/2} - 1)$ and that
  consequently $\ell(x) = \ln(1 + u(x))$. Notice that the last term in the previous
  bound is increasing in $z$. Consequently, since $u(x) > -\varepsilon$, we
  have that
  \[
  \ell^2(x) = \ln^2(1 + z) \leq u^2(x) + \frac{\varepsilon^3}{1-\varepsilon} .
  \]
  Recalling~\eqref{eq:i-function-constant}, which describes the result of
  applying the operator $i$ from \eqref{eq:integrand} to a function $f$ and a
  constant $C$, we find that
  \begin{align*}
    K\int_0^{x_*}i[\ell^2] \rmd x
    &\leq
    K\int_0^{x_*}i[u^2](x)\rmd x
    + \frac{K\varepsilon^3}{1-\varepsilon}\int_0^{x_*}\left\{\varphi(\Delta^\star/2 + x) + \varphi(\Delta^\star/2-x)\right\}\rmd x \\
    &\leq K\int_0^{x_*}i[u^2](x)\rmd x + \frac{ K^{1-3\beta}}{1-K^{-\beta}}
      x_* \sqrt{\frac{2}{\pi}}.
  \end{align*}
  By Lemma~\ref{lem:quadratic-lemma}, which gives an asymptotic expression for
  $K\int_0^{x_*}i[u^2](x)\rmd x$, the right hand side of the previous display is
  asymptotic to the right hand side of~(\ref{eq:quadratic-rates}). This implies
  that
  \begin{equation}\label{eq:l2-before-x0}
    K\int_0^{x_*}i[\ell^2] \rmd x
    =
    \begin{cases}
      O(K^{1 + 2r^\star - 2\beta^\star})
      &\text{ if }
      1/2 < \beta^\star \leq 3/4
      \\
      &\text{ and }
      \beta^\star - 1/2 \leq r^\star \leq \beta^\star /3;
      \\
      \tilde{O}(K^{1 - \beta^\star -
      \frac14\paren{\frac{\beta^\star
      - r^\star}{\sqrt{r^\star}}}^2})
      &\text{ if }
        1/2 < \beta^\star \leq 3/4 \text{ and }
        r^\star > \beta^\star / 3,\\
      &\text{ or if }
        3/4 < \beta^\star \leq 1 \\
      & \text{ and }
        r^\star \geq (1 - \sqrt{1 -
        \beta^\star})^2.
    \end{cases}
  \end{equation}
  \paragraph*{Bound for $\int_{x_*}^\infty i[\ell^2](x)\rmd x$} We follow a
  similar path as the proof of~\ref{lem:asymp-after-x0}. Write
  $i[\ell](x) = j_+(x) + j_-(x)$ with
  \begin{equation}\label{eq:jplus}
    j_+(x)
    =
    \ell^2(\Delta^\star / 2 + x)
    [(1 - \varepsilon^\star)\varphi(\Delta^\star / 2 + x)
    +
    \varepsilon^\star\varphi(\Delta^\star / 2 - x)],
  \end{equation}
  and
  \begin{equation}\label{eq:jminus}
    j_-(x)
    =
    \ell^2(\Delta^\star / 2 - x)[
    (1-\varepsilon^\star)\varphi(\Delta^\star / 2 - x)
    +
    \varepsilon^\star\varphi(\Delta^\star / 2 + x)].
  \end{equation}
  We bound $j_+$ and $j_-$ from above. The following claim, whose proof is at
  the end of this section, is useful.
  \begin{claim}\label{clm:log-ubd}
    If $x > 0$, then $\ln^2(1 + x) \leq \ln^2(x) + 2/\rme$.
  \end{claim}
  With the claim  and the fact that
  $u(\Delta^\star / 2 + x) = \rme^{\Delta(x - x_*)} - \varepsilon$ at hand, we obtain
  that
  \begin{align*}
    \ell^2(\Delta^\star/2 + x) &= \ln^2(1 + u(\Delta^\star/2 + x)) \\
                               & \leq \ln^2(1 + \rme^{\Delta(x - x_*)})\\
                               & \leq \Delta^2(x - x_*)^2  + 2/\rme ,
  \end{align*}
  where the first inequality is due the fact that $z\mapsto\ln^2(1 + z)$ is
  increasing for $z\geq 0$, and the second inequality uses Claim
  \ref{clm:log-ubd}. We use the previous inequality, together with a change of
  variables, to bound (\ref{eq:jplus}) and the first part of our integral of
  interest as follows:
  \begin{multline}\label{eq:jp-asymp}
    K\int_{x_*}^\infty j_+(x)\rmd x
    \leq
    K
    \Delta^2
    \int_{0}^\infty
    x^2
    [(1 - \varepsilon^\star)\varphi(x + \Delta^\star / 2 + x_*)
    +
    \varepsilon^\star\varphi(x - (\Delta^\star / 2 - x_*))]
    \rmd x\\
    +
    K\frac{2}{\rme}
    [(1 - \varepsilon^\star)\overline{\Phi}(\Delta^\star / 2 + x_*)
    +
    \varepsilon^\star\overline{\Phi}(x_* - \Delta^\star / 2 )].
  \end{multline}
  Since $\Delta = \Theta(\sqrt{\ln K})$, together Item~\ref{item:asymp-ccdf},
  Item~\ref{item:asymp-ccdf-pert} and Item~\ref{item:trunc-m2} from
  Lemma~\ref{lem:gaussian-things} imply that the right hand side of the previous
  display is asymptotic to the first term as $K\to\infty$, by considering two
  cases: (1) if $r^\star < \beta^\star$, the right hand side
  of~(\ref{eq:jp-asymp}) is
  $\tilde{\Theta}(K^{1 - \beta^\star - \frac14\paren{\frac{r^\star -
        \beta^\star}{\sqrt{r^\star}}}^2})$; (2) if $r^\star \geq \beta^\star$,
  then the right hand side of~(\ref{eq:jp-asymp}) is
  $\tilde{\Theta}(K^{1 - \beta^\star})$.
  For $j_-$ in~(\ref{eq:jminus}), we use that
  $\ell^2(\Delta^\star / 2 - x) = \ln^2(1 - \varepsilon + \rme^{-\Delta(x +
    x_*)}) \leq \ln^2(1 - \varepsilon)$ because $z\mapsto \ln^2(z)$ is
  decreasing for $z\in (0, 1]$. Consequently,
  \begin{equation*}
    K\int_{x_*}^\infty j_-(x) \rmd x
    \leq
    K\ln^2(1 - \varepsilon)
    [
    (1-\varepsilon^\star)\overline{\Phi}(\Delta^\star / 2 - x_*)
    +
    \varepsilon^\star\overline{\Phi}(\Delta^\star / 2 + x_*)].
  \end{equation*}
  The previous expression implies that, since $\overline{\Phi} \leq 1$ and
  $\ln(1 - \varepsilon)\sim K^{-\beta}$ as $K\to\infty$, then
  $ K\int_{x_*}^\infty j_-(x) \rmd x = O(K^{1 - 2\beta})$ and consequently
  $K\int_{x_*}^\infty j_-(x) \rmd x\to 0$ as $K\to\infty$ because $\beta>1/2$.
  Use the asymptotic upper bounds of the integrals of $j_+$ and $j_-$ and the fact that
  $i[\ell^2] = j_+ + j_-$ to conclude that
  \begin{equation}\label{eq:l2-after-x0}
    K\int_{x_*}^{\infty}i[\ell^2] \rmd x
    =
    \begin{cases}
      \tilde{O}(K^{1 - \beta^\star})
      &\text{ if }
        r^\star \geq \beta^\star
      \\
      \tilde{O}(K^{1 - \beta^\star -
      \frac14\paren{\frac{\beta^\star
      - r^\star}{\sqrt{r^\star}}}^2})
      &\text{ if }
        \rho(\beta^\star) < r^\star <  \beta^\star.
    \end{cases}
  \end{equation}

  \paragraph*{Conclusion} Putting together~(\ref{eq:l2-before-x0})
  and~(\ref{eq:l2-after-x0}), we obtain~(\ref{eq:m2k-rates}), which is what was
  to be shown.
\end{proof}

We now prove Claim~\ref{clm:log-ubd}.
\begin{proof}[Proof of Claim~\ref{clm:log-ubd}]
  Let $h(x) = \ln^2(x)$. Then $h'(x) = 2\ln(x)/x$ and
  $h''(x) = 2(1 - \ln(x))/x^2$. Notice that $h'(x)\to 0$ as $x\to\infty$ and
  $h'(x)\to-\infty$ as $x\to 0$. Then, $x\mapsto h'(x)$ is maximized at the
  point where $h''(x) = 0$, which is at $x = \rme$, because $h'(\rme) = 2/\rme > 0$. The
  mean value theorem implies that there is a value $c\in [0, 1]$ such that
  $h(1 + x) - h(x) = h'(c)$. Since $h'(c)\leq h'(\rme) = 2/\rme$, the claim
  follows.
\end{proof}

\section{Comparison to a Bonferroni correction and prequential plug-in methods}
\label{app:comparison-bonferroni}

We compare the performance of the AV test that uses the mixture martingale from
Section~\ref{sec:prior}---the proposal of this work summarized in Definition~\ref{def:adaptive-test-mart}---to a test that uses a
Bonferroni-type correction and to a test that uses a prequential plug-in
martingale as~(\ref{eq:preq-plug-in-mart}). Section~\ref{sec:bonf-corr}
describes the Bonferroni-type approach and our implementation, and
Section~\ref{sec:preq-plug-meth} describes the predictable estimator that we use
for~(\ref{eq:preq-plug-in-mart}) in the comparison.
Section~\ref{sec:bonfi-sim-results} shows the results of the comparisons.

\subsection{Bonferroni Correction}\label{sec:bonf-corr}

In building a test against the composite alternative
in~(\ref{eq:generic-alternative}), a first idea might be to use an already
existing fixed-sample test repeatedly over time and to use a Bonferroni-type
correction on top to account for repeated testing. We now describe the
implementation of this idea to construct and AV test based on $\mathrm{HC}_t$
from~(\ref{eq:hc-statisic}). For
$\alpha \in [0, 1]$, let $h(\alpha)$ the $\alpha$ quantile of $\mathrm{HC}_t$
under the null hypothesis; that is, $h(\alpha)$ satisfies
$\mathbf{P}_0\{\mathrm{HC}_t \geq h(\alpha)\} = \alpha$. Notice that, since
$\mathrm{HC}_t$ takes as inputs the null-standardized stream means
$(Z_{1, t}, \dots, Z_{K, t})$ with
$Z_{i, t} = (1 / \sqrt{t})\sum_{s\leq t}X_{i, s}$, and
$Z_{i, t}\sim \mathrm{Normal}(0, 1)$, the null distribution of
$\mathrm{HC}_t$ does not depend on $t$---and neither does $h(\alpha)$. Consider
$(w_t)_{t\geq 1}$ a probability distribution on $\{1, 2, \dots\}$ and define the
sequential test that stops at
$\tau_{\mathrm{HC}, w} = \inf\{t : \mathrm{HC}_t \leq h(w_t\alpha)\}$ and
rejects the null hypothesis at $\tau_{\mathrm{HC}, w}$ whenever
$\tau_{\mathrm{HC}, w} < \infty$. The fact that this test has type-I error
smaller than $\alpha$ is a consequence of a union bound; indeed, the type-I
error is
\begin{equation*}
  \mathbf{P}_0\{\tau_{\mathrm{HC}} < \infty\}
  =
  \mathbf{P}_0\{\exists t : \mathrm{HC}_t < h(w_t\alpha)\}
  \leq
  \sum_t   \mathbf{P}_0\{ \mathrm{HC}_t < h(w_t\alpha)\}
  =   \sum_t  w_t\alpha
  = \alpha.
\end{equation*}

The practical implementation of the sequentialized test that uses
$\mathrm{HC}_t$ and the Bonferroni correction as described requires two
ingredients: a choice of weights $(w_t)_{t\geq 1}$, and finding the quantiles
$h(w_t\alpha)$ for a specified level $\alpha$ and $t\geq 1$. Later, we discuss
challenges associated with these choices. An good choice of weights
$(w_t)_{t\geq 1}$ would put more mass at the times at which the signals become
detectable; that is, at the times $t\geq t^\star$, where $t^\star$ is as in
Theorem~\ref{thm:lr_achieves_boundary}. However, such a choice would require
knowledge of the true signal strength $\delta$, which is unknown. Thus, in
experiments, we use a default choice of weights, $w_t = \frac{1}{t(t + 1)}$. For
the second ingredient, the quantiles $h(w_t\alpha)$ for $t\geq 1$, we use Monte
Carlo simulation. Indeed, using a sample of $M$ independent null realizations
$\mathrm{HC}_{1, 1}, \dots, \mathrm{HC}_{1, M} $ of $\mathrm{HC}_1$, we estimate
$\alpha \mapsto h(\alpha)$ using the generalized inverse
$\alpha \mapsto \hat{h}_M(\alpha) := \inf\{x: \hat{F}_{\mathrm{HC}, M}(x) \geq
\alpha \}$ of the empirical cumulative distribution function
$\hat{F}_{\mathrm{HC}, M}(x) = (1 / M)\sum_{i\leq M}\indicator{\mathrm{HC}_{1,
    i} \leq x}$.

According to the previous discussion, we use weights $w_t = 1 / [t(t + 1)]$
and the estimand $\hat{h}_M$ with $M = 10^6$ to implement an approximate version
of the sequentialized $\mathrm{HC}_t$ test with a Bonferroni correction. These
choices come with certain challenges. Notice that any choice of weights $w_t$
that puts positive mass on all $t\geq 1$---and consequently also our choice
$w_t = 1 / [t(1 + t)]$---penalizes large times. Hence, any choice of $w_t$
with full support necessarily becomes more and more conservative when $t^\star$
is large. Additionally, in order to maintain an (approximate) type-I error
guarantee, it is necessary to count on good estimates of each $w_t\alpha$ null
quantiles of $\mathrm{HC}_{t}$, which become larger---farther in the tail of the
null distribution of $\mathrm{HC}_{t}$---as $t$ grows larger. Our current choice
of estimand, $\hat{h}_M(w_t\alpha)$, can be improved using, for instance,
importance-weighted sampling. However, the main goal of this section is to
present a comparison between this test and our approach. To our purposes, it is
enough to choose a setting that is, in principle, favorable to the
sequentialized test that uses $\mathrm{HC}_t$ and the Bonferroni correction, and
does not require further optimization.

\subsection{Prequential plug-in method}\label{sec:preq-plug-meth}

In order to specify an AV test based on~(\ref{eq:preq-plug-in-mart}), we need to
specify a sequence of estimators
$((\hat{\varepsilon}_t, \hat{\delta}_t))_{t \geq 1}$. Recall from the discussion
around~(\ref{eq:preq-plug-in-mart}) that, in order to preserve the martingale
property, the sequence needs to be predictable; that is, for each $t$, the pair
$(\hat{\varepsilon}_t, \hat{\delta}_t)$ of estimators need to depend only on
data up to and including $t - 1$, not on data at time $t$. We use a smoothed
maximum likelihood approach and use Expectation Maximization (EM)
algorithm~\citep{dempster_maximum_1977} (which is standard in fitting Gaussian
mixtures) to compute it~\citep[see also][Section~8.5.1]{hastie_elements_2009}.
The popularity of the EM algorithm is explained by its ease of implementation
and the fact that it is guaranteed to converge to a local maximum of the
likelihood function if ran long enough.
Below, we describe our implementation of the EM algorithm that, given
$(z_1, \dots, z_K)\in \mathbf{R}^K$, computes (approximate) maximum likelihood
estimates $\hat{\varepsilon}(z_1, \dots, z_K)$ and
$\hat{\delta}(z_1, \dots, z_K)$ of the parameter $\varepsilon^\star$ and
$\delta^\star$. The estimators that we use for $M_t$
in~(\ref{eq:preq-plug-in-mart}) are defined as
$(\hat{\varepsilon}_1, \hat{\delta}_1) = (0, 1)$---the first multiplicative
increment of $M_t$ is~1---and, for $t > 1$, we use
$\hat{\varepsilon}_t = \hat{\varepsilon}(Z_{1, t - 1}, \dots, Z_{K, t - 1})$ and
$\hat{\delta}_t = \hat{\delta}(Z_{1, t- 1}, \dots, Z_{K, t - 1})$, where, for
each $i\in [K]$ and each $t \geq 1$, the statistic
$Z_{i, t} = (1 / \sqrt{t})\sum_{s\leq t}X_{i, s}$ is the standardized
single-stream average.

We now describe our implementation of the EM algorithm without further
motivation---see \citep[][Section~8.5.1]{hastie_elements_2009} for a more
in-depth discussion. Let $\hat{\varepsilon}_0\in (0, 1)$ and $\hat{\delta}_0 > 0$ be the
initialization points for the EM algorithm, and let the $m_{\max}\in \mathbf{N}$
the maximum number of iterations we allow the EM algorithm to take. Given
$(z_1, \dots, z_K)$, $\hat{\varepsilon}_0$, $\hat{\delta}_0$ and $m_{\max}$, we
perform the following steps for $l = 1, \dots, m_{\max}$:
\begin{enumerate}
\item for $i\in [K]$, compute
  $\pi_{i, l} = \frac{\hat{\varepsilon}_{l - 1}\exp\paren{\hat{\delta}_{l - 1} z_i -
      \hat{\delta}_{l - 1}^2 / 2}}{1 - \hat{\varepsilon}_{l - 1} +
    \hat{\varepsilon}_{l - 1}\exp\paren{\hat{\delta}_{l - 1} z_i - \hat{\delta}_{l - 1}^2 / 2}}$;
\item compute $\hat{\delta}_{l} = \sum_{i\in [K]}\pi_{i, l} \  z_{i}$;
\item compute
  $\hat{\varepsilon}_l = \frac{1}{K}\sum_{i\in [K]}\pi_{i, l}$;
\item if $|\ln \hat{\varepsilon}_l - \ln \hat{\varepsilon}_{l - 1}| < 10^{-4}$
  and $|\hat{\delta}_{l}- \hat{\delta}_{l-1}| < 10^{-4}$, stop and return
  $(\hat{\varepsilon}_l,\hat{\delta}_l)$.
\end{enumerate}
If the algorithm fails to converge after $m_{\max}$ steps, we use
$(\hat{\varepsilon}, \hat{\delta}) = (0.0, 1.0)$---in this case the
multiplicative increment of $M_t$ takes the value $1$. In our experiments, we
use $\hat{\varepsilon}_0 = 1/ K$, $\hat{\delta}_0 = \sqrt{2\ln K}$, which are in the
parameter zone of interest, and $m_{\max} = 10^3$. In addition, the overall test
performance remained almost the same when using a more stringent stopping
criteria or a larger maximum number of iterations.

\subsection{Simulation Results}
\label{sec:bonfi-sim-results}
\begin{figure}[h!]
  \centering
  \includegraphics{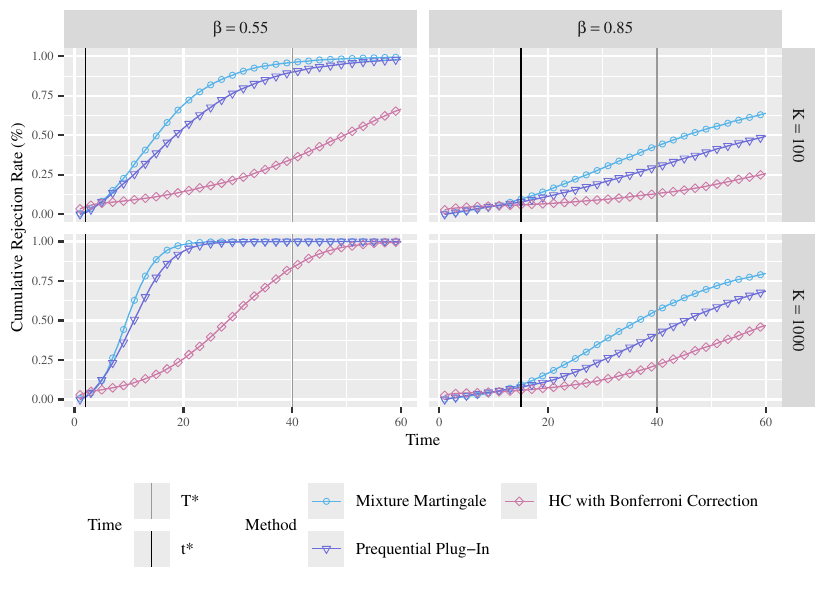}
  \caption{Cumulative rejection rate for the sequentialized $\mathrm{HC}_t$ with
    Bonferroni correction and the adaptive construction from
    Section~\ref{sec:prior}. The former test is stops more often at earlier
    times, but its performance is overall worse than the latter
    (see~Appendix~\ref{app:comparison-bonferroni} for a discussion of the results
    and Appendix~\ref{app:simulation-details} for more details about the
    simulation). Detection is possible between $T^\star$ and $t^\star$, when it
    is not always possible to estimate the anomalies. Only every other time point
    is marked to avoid overplotting. See more details in
    Section~\ref{sec:simulations}.}
  \label{fig:bonferroni-power}
\end{figure}

Under several alternatives in ${\cal H}_1$, Figure~\ref{fig:bonferroni-power}
shows estimates the cumulative rejection rates for three methods: the
sequentialized HC test described in Section~\ref{sec:bonf-corr}, the prequential
plug-in method from Section~\ref{sec:preq-plug-meth}---which monitors $M_t$
from~(\ref{eq:preq-plug-in-mart})---, and our proposal, the AV test that
monitors the mixture martingale $E_t(\mathbf{\Pi})$ from
Section~\ref{sec:prior}. Let $\tau^{\mathbf{\Pi}}$ be the random sample size of
AV test that uses $E_t(\mathbf{\Pi})$; $\tau_{\mathrm{HC}, w}$, that of the
sequentialized HC test; and $\tau_{\mathrm{preq}}$, that of AV test associated
to the prequential plugin test martingale. We show estimates of the cumulative
rejection rates $t \mapsto \mathbf{P}\{\tau^{\mathbf{\Pi}} \leq t\}$,
$t \mapsto \mathbf{P}\{\tau_{\mathrm{HC}, w} \leq t\}$, and
$t \mapsto \mathbf{P}\{\tau_{\mathrm{preq}} \leq t\}$. To this end, the same
procedure was used as the one used in Section~\ref{sec:simulations} for
generating Figure~\ref{fig:power} (see further details in
Appendix~\ref{app:simulation-details}). In Figure~\ref{fig:bonferroni-power}, it
can be observed that the $\mathrm{HC}$-Bonferroni is more aggressive (it stops
more frequently) at times $t < t^\star$ and its rejection rate is significantly
lower than that of the test based on either $E_t(\mathbf{\Pi})$ or $M_t$ at
later times. This higher rate at the beginning of the experiment can be
explained by the fact that, even under the null hypothesis, the
$\mathrm{HC}$-Bonferroni test is expected to reject with probability
(approximately) $\alpha / 2$ at the first time step. In contrast, the test based
on $E_t(\mathbf{\Pi})$ is more conservative under the null hypothesis at this
first time step (not shown). The test based on the prequential plug-in method
consistently underperforms our proposal.

\section{Simulation Details}
\label{app:simulation-details}

This section contains the details of the simulations in
Section~\ref{sec:simulations}.

\paragraph*{Test calibration} We fix $T^\star = 40$ and $N = 10^4$. In order to
calibrate the AV tests in each simulation scenario, that is, to estimate the
thresholds $A^\star$ and $A^{\mathbf{\Pi}}$ so that both the tests based on the
log-optimal likelihood ratio and the mixture likelihood ratio have type-I error
$\alpha$, we sample $N$ realizations of the processes
$t\mapsto E^\star_t, E_t(\mathbf{\Pi})$ for $t \leq 2T^\star$ under the null
hypothesis. Under the null hypothesis $\ln E^\star_t \to K \ln(1 - \varepsilon)$
almost surely as $t\to\infty$ and preliminary simulations show that this
convergence happens well before the threshold that was used, $2T^\star$. This
choice of threshold was made to shorten the computation times; the threshold $C$
(which was taken to be $5T^\star$ in simulations) can also be used with
identical results. A similarly fast convergence to its null asymptotic value is
observed for $t\mapsto \ln E_t(\mathbf{\Pi})$. For each of the $N$ realizations
of $t\mapsto E^\star_t, E_t(\mathbf{\Pi})$, the maxima over time,
$\max_{t} E^\star_t$ and $\max_t E_t(\mathbf{\Pi})$, are
computed, and both $A^\star$ and $A^{\mathbf{\Pi}}$ are estimated as the sample
$1 - \alpha$ empirical quantiles of the $N$ realizations of the maxima,
respectively.
\paragraph*{Details for Figure~\ref{fig:power}} We fix $T^\star = 40$ and
$N = 10^4$. Once the previous calibration is performed, for each
$\varepsilon^\star = K^{-\beta^\star}$ with $\beta^\star \in \{0.55, 0.85\}$ and
each $K\in \{100, 1000\}$, $N$ realizations (independent from the ones used to
calibrate the tests) of $t\mapsto E^\star_t, E_t(\mathbf{\Pi})$ are sampled for
$t \leq 2T^\star$ from which $N$ realizations $\tau^\star\wedge 2T^\star$ and
$\tau^{\mathbf{\Pi}}\wedge 2T^\star$ are obtained. In computing
$E_t(\mathbf{\Pi})$, we used $C = 5T^\star$, but similar results are observed
for a wide range of choices of $C$. Regarding the choice of length of each
sample path, $2T^\star$, it is arbitrary, but a large number is needed to ensure
the tests stop with high probability before the end of the simulation under the
alternative in order to estimate the desired cumulative rejection rates and the
quantiles. The cumulative rejection rates $F_\tau^\star(t)$ and
$F_\tau^{\mathbf{\Pi}}(t)$ are estimated using empirical cumulative distribution
functions for $t \leq 2T^\star$. The thresholds for the likelihood ratio test
and the test based on $\mathrm{HC}_t$ are obtained using the upper $5\%$
quantile of $10^4$ samples from the null distribution. The power of both
fixed-sample tests, at each sample size $t$, are estimated as the fraction of
rejections over $10^4$ realizations of the statistic under each alternative
distribution.

\paragraph*{Details for Table~\ref{tab:average-st}} Using the same sampling
procedure as for Figure~\ref{fig:power} for each test martingale $E_t$---the
well-specified likelihood ratio and the mixture---we sampled $N = 10^4$
realizations of the stopping time $\tau = \inf\{t : E_t \geq 1/\alpha\}$ under
each $\beta^\star \in \{0.55, 0.85\}$ and each $K\in \{100, 1000\}$. We
estimated $n(0.8)$, the smallest $n$ such that
$\mathbf{P}^\star\{\tau \leq n\} \geq 0.8\gamma_{\max}$ (see
Section~\ref{sec:quantile-comparison} for the computation of the maximum power
$\gamma_{\max}$) using the empirical distribution function of the stopping
times. Indeed, if $\tau_1, \dots, \tau_N$ are the $N$ realizations of the
stopping time $\tau$, then we use
$\hat{F}(n) = (1 / N)\sum_{i\in [N]}\indicator{\tau_i \leq n}$ and estimated
$n(0.8)$ with $\hat{n}(0.8)$, the smallest $n$ such that
$\hat{F}(n)\leq 0.8\gamma_{\max}$. Then, we estimate
$\mathbf{E}^\star[\tau \wedge n(0.8)]$ with
$(1 / N)\sum_{i\in [N]}\tau_i\wedge \hat{n}(0.8)$. To estimate the sample
size needed to reach $80\%$ of the maximum power $\gamma_{\max}$ for
fixed-sample tests, we repeated the same procedure as in the previous paragraph
and found the smallest sample size that resulted in estimated power equal or
larger to $80\%$ of the maximum power $\gamma_{\max}$.

\end{appendix}

\end{document}